\PassOptionsToPackage{figuresright}{rotating}
\PassOptionsToPackage{table}{xcolor}
\documentclass[final,5p,times,twocolumn]{elsarticle}
\usepackage[utf8]{inputenc}
\usepackage{amsmath}
\usepackage{algorithm}
\usepackage{algpseudocode}
\usepackage{svg}
\usepackage{amsfonts}
\usepackage[table]{xcolor}
\usepackage{soul}
\soulregister\ref7
\soulregister\textit7
\usepackage{fancyhdr}
\usepackage{blindtext}

\usepackage{pdflscape}
\usepackage{float}
\usepackage[colorinlistoftodos]{todonotes}
\usepackage{changepage}
\usepackage{nomencl}
\usepackage{graphicx}
\usepackage{rotating}
\makenomenclature
\usepackage{rotating}
\usepackage{multirow}
\usepackage{tabularx}
\usepackage{array}
\usepackage[figuresright]{rotating}
\usepackage{booktabs}
\usepackage{longtable}
\usepackage{hyphenat}
\usepackage{pdflscape}
\usepackage{longtable}
\usepackage{enumitem}
\usepackage{geometry}
\usepackage{siunitx}
\usepackage{multirow}
\usepackage{caption} 
\usepackage{enumitem}
\usepackage{array}
\usepackage{lscape}
\usepackage{microtype}
\usepackage{adjustbox}
\usepackage[english]{babel} 
\usepackage{microtype} 
\usepackage{caption}
\usepackage{xcolor}
\usepackage{ragged2e}
\usepackage{makecell}
\usepackage{amssymb} 
\usepackage{placeins}
\usepackage[font=normalsize,labelfont=bf]{caption}
\usepackage{booktabs}

\def\tsc#1{\csdef{#1}{\textsc{\lowercase{#1}}\xspace}}
\tsc{WGM}
\tsc{QE}
\tsc{EP}
\tsc{PMS}
\tsc{BEC}
\tsc{DE}
\begin{document}
\let\WriteBookmarks\relax
\def\floatpagepagefraction{1}
\def\textpagefraction{.001}

\title{Reliability as a Design Principle: A Systematic Review and Integrated Framework for Renewable-Based Microgrids}

\author[a]{Mohammed Zeehan Saleheen}
\ead{mohammed.saleheen@monash.edu}

\author[a]{Markus Wagner}
\ead{markus.wagner@monash.edu}

\author[b]{Reza Razzaghi}
\ead{reza.Razzaghi@monash.edu}

\author[a]{Hao~Wang\corref{cor1}}
\ead{hao.wang2@monash.edu}
\cortext[cor1]{Corresponding author: Hao Wang.}

\affiliation[a]{
      organization={Department of Data Science and AI, Faculty of IT, Monash University},
      city={Melbourne},
      country={Australia}
}
\affiliation[b]{
      organization={Department of Electrical and Computer Systems Engineering, Monash University},
      city={Melbourne},
      country={Australia}
}

\begin{abstract}
   Reliable operation is a central motivation for deploying renewable-based microgrids, yet reliability is frequently treated as a secondary outcome of cost-driven design rather than as a governing planning objective. Existing reviews remain fragmented across techno-economic sizing, optimization methods, and reliability assessment, limiting the translation of planning outcomes into operationally reliable systems. This paper presents a systematic rapid review that positions reliability as the central organizing principle for microgrid design. Specifically, this review systematically synthesizes recent literature to examine how planning assumptions, optimization formulations, operational flexibility mechanisms, and reliability assessment frameworks jointly shape reliability outcomes. The synthesis shows that reliability in renewable-based microgrids is governed primarily by chronological, time-coupled energy adequacy rather than installed capacity alone, with \textit{Dunkelflaute} events emerging as a key determinant of adequacy failure. Reliability outcomes are shaped by the joint interaction of resource portfolios, storage operating policies, and state trajectories, network features, and protection feasibility under inverter-dominated operation. The review further demonstrates that reliability indices inherited from conventional power systems are poorly suited for renewable-based microgrids, as they compress performance into single dimensions and obscure temporal, spatial, and service-critical risk concentrations. Across optimization practice, reliability is increasingly embedded through multi-objective and constrained formulations; however, persistent gaps remain in representing correlated renewable scarcity, mission-profile-dependent component reliability, and interruption valuation (e.g., value of lost load and customer damage functions) in a consistent and decision-relevant manner. Overall, this review consolidates planning factors, optimization approaches, reliability evaluation methods, and metric suitability into an integrated roadmap for reliability-centered microgrid planning, and outlines future directions toward state-aware, service-oriented planning and assessment frameworks.
\end{abstract}


\maketitle

\renewcommand{\nomname}{List of Acronyms and Abbreviations}
\printnomenclature

\section{Introduction}
\label{Introduction}

The global transition toward decarbonized energy systems has fundamentally reshaped the architecture and operating principles of electrical power networks. Driven by climate imperatives, policy commitments, and declining technology costs, renewable energy sources accounted for approximately 30\% of global electricity generation in 2023, with continued acceleration projected toward net-zero targets~\cite{IEA2024RenewablesElectricity}. While this transition is essential for environmental sustainability, it represents a structural departure from conventional power systems dominated by large, centrally dispatched synchronous generators to variable renewable energy resources. However, these sources introduce intermittent and weather-dependent generation profiles, reduced inertia and fault current contributions, and fundamentally different dynamic behaviors, altering the operating conditions under which power systems must ensure reliable electricity supply.

In this evolving energy landscape, ensuring reliability, the ability of a system to deliver electricity to customers within acceptable standards and in the quantities demanded~\cite{allan2013reliability}, has emerged as a central challenge for modern power systems~\cite{islam2024improving}. This challenge is further exacerbated by climate-driven extreme weather events~\cite{alhazmi2025resilient}, aging network infrastructure~\cite{kopsidas2017power}, and the rapidly increasing societal dependence on electricity for critical services~\cite{poudel2018critical}. These challenges are particularly acute in remote, rural, and edge-of-grid regions, where outages are more frequent and restoration times are longer. In such contexts, the limitations of traditional, centrally supplied distribution systems, originally designed for predictable, synchronous generation, become increasingly evident in maintaining reliable electricity supply.

Microgrids have emerged as a promising architectural paradigm for addressing these above-mentioned challenges in this evolving energy landscape~\cite{hittinger2015evaluating}. Defined by IEEE Std 2030.7 as a group of interconnected loads and distributed energy resources within clearly defined electrical boundaries that acts as a single controllable entity with respect to the grid, microgrids can operate in either grid-connected or islanded modes~\cite{danley2019defining}. This operational flexibility represents a fundamental distinction from the conventional distribution system, enabling reduced dependence on upstream network availability~\cite{rosales2019microgrids,gali2023adaptive}. The reliability benefits are particularly pronounced in renewable-based microgrids, where the strategic integration of distributed generation, energy storage, and intelligent control systems can maintain electricity supply during grid outages, offering a pathway to enhanced reliability for both grid-edge communities and remote areas lacking grid connection entirely. However, realizing these potential benefits requires systematic attention to reliability throughout the planning and design process, an objective that proves considerably more challenging than conventional distribution system planning due to the unique characteristics of microgrids.

The challenge of reliability-oriented microgrid design stems from multiple factors that distinguish these systems from traditional power networks. Unlike large interconnected power systems, where reliability is primarily ensured through network redundancy, firm generation reserves, and spatial aggregation of resources, microgrids exhibit a strong dependence between long-term planning decisions and short-term operational feasibility under uncertainty~\cite{pang2024long}. Decisions made during the planning stage, such as the sizing and placement of distributed generation, energy storage, and network assets, directly determine whether feasible operating states exist during periods of renewable scarcity or component outages. As a result, designs that satisfy conventional planning criteria or perform well under average conditions may still experience significant reliability degradation when subjected to realistic operating constraints and stochastic disturbances. This challenge is particularly acute in renewable-dominated microgrids, where solar and wind generation are inherently variable and energy-limited, and where storage dynamics couple past operational decisions to future supply adequacy. Under these conditions, reliability outcomes are governed by the joint realization of generation availability, load demand, storage state, and component failures, rather than by installed capacity or reserve margins alone. Consequently, deterministic planning approaches that decouple planning from operation, which are often adequate for traditional power systems~\cite{capasso2015new}, are insufficient to capture the true reliability behavior of renewable-based microgrids. 

In addition, microgrids are required to maintain acceptable reliability performance across two fundamentally different operating modes~\cite{ganjian2019seamless}, namely grid-connected operation and islanded operation, each of which imposes distinct technical and reliability requirements. During grid-connected operation, the upstream utility network typically provides voltage and frequency reference, reserve capacity, and contingency support, allowing local generation and storage to operate with greater flexibility. In contrast, during islanded operation the microgrid must independently balance generation and demand while simultaneously maintaining voltage stability, frequency regulation, and power quality within tight limits. The transition between these modes, whether planned or unplanned, introduces additional reliability risks that are not adequately represented by steady-state planning assumptions. Moreover, the widespread use of inverter-based interfaces for renewable generation and energy storage introduces failure modes and degradation mechanisms that differ fundamentally from those of conventional electromechanical equipment~\cite{peyghami2020incorporating,9648165}. For example, power electronic converters exhibit a higher sensitivity to thermal stress, control interactions, and coordination issues~\cite{ahmed2023dynamic}, and their failure behaviors are often abrupt rather than gradual~\cite{sandelic2022reliability}. However, most established reliability assessment methods were developed for conventional generation and network components, and therefore do not accurately capture the operational dependencies, dynamic behavior, and failure characteristics of inverter-dominated microgrids~\cite{zhang2020reliability}.

In addition, reliability-oriented microgrid design is further complicated by the limitations of conventional protection schemes under bidirectional power flow and above-mentioned inverter-dominated operation~\cite{basu2023dynamic}. Traditional protection systems were developed for radial networks with unidirectional power flow and high fault currents supplied by synchronous generators~\cite{kennedy2016review}. In contrast, microgrids with high penetration of distributed generation and energy storage frequently experience bidirectional power flows~\cite{nandi2023coordination}, variable fault current contributions~\cite{brearley2017review}, and operating conditions that change between grid-connected and islanded modes~\cite{ganjian2019seamless}. Moreover, inverter-based resources typically exhibit limited and controlled fault current injection~\cite{stankovic2021fault}, which can undermine the effectiveness of traditional overcurrent-based protection and disrupt coordination among protective devices. As a result, protection miscoordination, delayed fault isolation, or unnecessary disconnections can occur, directly affecting reliability outcomes. 

Furthermore, the cost-reliability trade-off in microgrid design is also challenging due to the non-smooth and highly nonlinear relationship between investment decisions and achieved reliability outcomes. In renewable-based microgrids, incremental investments in generation, energy storage, or network capacity may not result in proportional improvements in reliability~\cite{sakthivelnathan2024cost}, as reliability is governed by binding adequacy constraints related to energy balance, storage state of charge, and contingency coverage. As a result, modest changes in component sizing or placement may yield little reliability benefit until critical levels of redundancy or flexibility are reached, beyond which reliability can improve sharply. This threshold-driven behavior contrasts with traditional power systems, where reserve margins and network reinforcement typically lead to more gradual and predictable reliability gains. Consequently, the cost-reliability relationship in microgrids cannot be adequately represented by smooth trade-off functions, posing a challenge for reliability-oriented microgrid design.

Notably, despite substantial research attention to microgrid planning and design, the existing review literature exhibits several notable gaps that limit its utility for practitioners seeking to implement reliability-oriented design approaches. The current body of literature remains fragmented in how it frames reliability and, in many cases~\cite{weng2022reliability,adefarati2017reliability,riou2021multi}, positions it as a peripheral performance attribute rather than a central decision principle for design decisions. A large subset of reviews is dominated by techno-economic design and cost-optimal sizing~\cite{veilleux2020techno,thomas2016optimal,jeyaprabha2022probabilistic,elborlsy2025case,yadav2025hybrid,ali2025evaluating}, where reliability is either implicitly assumed through generic reserve margins or represented using a narrow set of customer-weighted interruption indices, without examining whether these indices remain valid for inverter-dominated and renewable-based operation. Similarly, several studies focus primarily on optimization algorithms and solution techniques, but treat reliability modeling as a simplified constraint set~\cite{amir2019reliability,huo2022reliability,sarfi2018economic}, rarely discussing the underlying reliability assumptions, uncertainty representations, and operational feasibility mechanisms that drive reliability outcomes in practice.

Moreover, existing reviews generally treat key design principles in isolation, without providing integrated frameworks that demonstrate how these factors interact in reliability-oriented microgrid design. For instance, some reviews primarily emphasize optimization formulations and solution techniques~\cite{dawoud2018hybrid,adewuyi2024appraisal,kumar2024review}, while others focus on reliability assessment methods in isolation~\cite{ou2023review,singh2023comprehensive,sandelic2022reliability}. Similarly, technical design factors are often examined separately from economic or social considerations, with limited discussion of how these dimensions jointly influence reliability-oriented microgrid design. In addition, reviews focused specifically on reliability assessment \cite{garip2022power,singh2023comprehensive,muhtadi2021distributed} typically discuss reliability indices and evaluation methods, but provide limited guidance on how evaluation results should inform planning decisions. These studies often emphasize component-level modeling, adequacy indices, or probabilistic simulation approaches, yet do not systematically connect these evaluation outputs to concrete planning decisions such as DER sizing, siting, network reinforcement, and operational strategy. Thus, these studies leave a disconnect between assessment and planning, which limits design interpretability and actionable guidance. For example, Sandelic et al. \cite{sandelic2022reliability} explicitly address reliability aspects in microgrid design and identify the critical gap regarding power electronics reliability, yet their treatment of planning factors remains limited given its focus on power electronics challenges. A further limitation is the inconsistent treatment and classification of reliability indices across reviews~\cite{garip2022power,aruna2021comprehensve,lopez2020reliability}, where indices are frequently mixed without clear definitions, boundaries, or justification for their applicability to microgrids. Hence, such classification schemes occasionally conflate distinct conceptual dimensions, providing limited guidance for practitioners attempting to select appropriate metrics for specific planning contexts.

The current review literature falls short of offering an integrated and reliability-oriented perspective on microgrid design. This gap motivates the need for a holistic review that consolidates key planning factors in reliability-oriented microgrid design, systematically examines optimization techniques for balancing reliability with competing objectives, and evaluates the suitability of existing reliability metrics and modeling approaches for modern renewable-based microgrids. Hence, this study aims to bridge these gaps by providing a state-of-the-art review that positions reliability as the central organizing principle guiding microgrid planning and design decisions. Adopting a rapid review methodology, this paper offers a unified and integrated perspective on reliable microgrid design by jointly examining key design factors, optimization techniques, and reliability assessment methods and by critically evaluating their suitability for modern renewable-based microgrids.

A rapid review is a structured evidence synthesis approach that employs a streamlined version of systematic review procedures to provide a timely and focused overview of the most relevant literature~\cite{garritty2024updated}. This methodological choice is motivated by the fast-evolving nature of microgrid research, where new technologies, operational frameworks, and planning approaches continue to emerge, making it essential to capture current challenges and state-of-the-art developments. While widely used in health sciences~\cite{RN1030}, rapid reviews are relatively unexplored in engineering fields. However, in dynamic domains such as microgrids, rapid reviews offer an effective balance between timeliness and methodological rigor by prioritizing high-quality and directly relevant studies. To ensure robustness and transparency in study selection, the AMSTAR2 (A Measurement Tool to Assess Systematic Reviews version 2) framework~\cite{shea2017amstar} is employed to assess the methodological quality and relevance of the included literature. Unlike traditional reviews, which are often time-intensive~\cite{RN1031} and may incorporate outdated evidence, these combined approaches enable a focused synthesis that is well aligned with the objectives of reliability-oriented microgrid design research.

Accordingly, a rapid review methodology is adopted to address the need for timely and relevant evidence synthesis in the rapidly evolving field of renewable-based microgrid design. As described in the preceding paragraphs, in recent years, reliability-oriented microgrid planning has progressed beyond conventional techno-economic sizing toward a more integrated and multi-dimensional design framework. While traditional systematic reviews aim to provide a comprehensive evidence base, their extended completion timelines, often spanning 12–24 months, limit their suitability for domains undergoing rapid methodological and technological change. Rapid reviews address this limitation by selectively streamlining review procedures, including prioritizing the most relevant and methodologically robust studies and applying a clearly defined chronological scope. Importantly, this approach does not compromise rigor; rather, it involves explicit and transparent trade-offs between comprehensiveness and timeliness. In this study, the applied chronological scope is designed to capture both foundational contributions and recent developments, ensuring that the synthesis remains historically grounded while reflecting current research directions. This is particularly important given the rapid evolution of reliability-oriented microgrid planning, driven by the increasing penetration of distributed renewable resources and the transition toward inverter-dominated systems. Overall, the adopted rapid review protocol provides a structured and transparent synthesis that is aligned with contemporary developments, while clearly defining the scope boundaries without compromising the internal validity of the review.

Building upon this foundation, the main contribution of this study lies in establishing reliability as a unifying analytical principle through which microgrid planning and design are systematically reviewed. Instead of treating reliability as a secondary performance metric, this review restructures existing studies based on how reliability considerations influence planning decisions, optimization formulations, and assessment methodologies. In particular, the study provides a structured examination of optimization techniques that explicitly trade off reliability against economic and operational objectives, and offers a critical evaluation of commonly used reliability indices with respect to their validity and interpretability in renewable-dominated microgrids. By integrating planning factors, optimization approaches, and reliability assessment methods within a unified review framework, this paper advances the current state of knowledge beyond fragmented surveys and provides a consolidated foundation for future methodological development in reliability-oriented microgrid design. The paper is organized as follows: Section~\ref{rapid review methodology} covers the methodology used to conduct the rapid review. A comprehensive guideline outlining the key factors of reliable microgrid design is presented in Section~\ref{Sec:3}. Section~\ref{sec:4} summarizes optimization techniques applied to microgrid design. Section~\ref{sec:5} analyzes reliability evaluation methods and metrics currently used for microgrids, highlighting their limitations. Finally, Section~\ref{6} provides further discussion and integrated synthesis, while Section~\ref{7} presents the key findings, recommendations, limitations, and conclusion.

\section{Rapid Review Methodology: Synthesizing Microgrid Literature}
\label{rapid review methodology}

Recently, rapid review methodology has emerged as a pragmatic approach for synthesizing evidence in fast-evolving research domains. Based on the foundation of traditional systematic review, a rapid review applies structured and transparent procedures while narrowing the scope of search and screening to enable timely synthesis of the most relevant literature~\cite{haby2024best}. Compared to traditional systematic reviews, which often require a year or more to complete~\cite{beller2013systematic}, rapid reviews are typically conducted within shorter time frames (six to eight months), allowing emerging technological developments and research trends to be captured more effectively~\cite{ganann2010expediting}. This characteristic is particularly relevant for microgrid research, where rapid advances in renewable integration, power electronics, and reliability strategies can render older evidence less representative of current practice.

\begin{figure*}[htbp]
    \centering
    \includegraphics[width=0.6\textwidth]{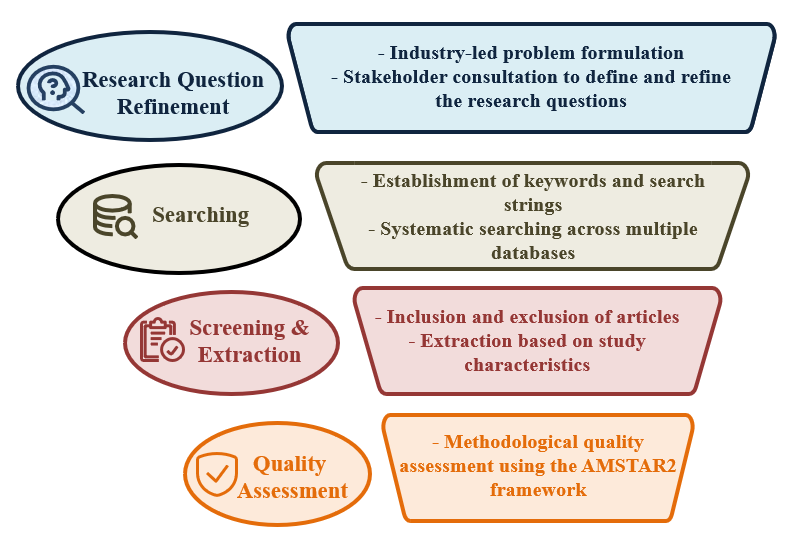}
    \caption{Rapid review methodology.}
    \label{Fig.1}
\end{figure*}

In this context, rapid reviews emphasize the synthesis of high-quality secondary literature and authoritative sources, including academic reviews, industry reports, and policy documents, supported by clearly defined inclusion and exclusion criteria~\cite{syed2021rapid}. While rapid review methodologies are well established in the health sciences, their application within the microgrid domain remains limited; to date, only a single study has been reported \cite{syed2021rapid}, published in 2021. Accordingly, this study applies a rapid review framework to examine reliability-oriented microgrid design and planning. The methodological workflow, where each step is recurrent and is interrelated with the subsequent step, is outlined in Figure~\ref{Fig.1}. An in-depth discussion of these steps is presented in  the following subsections.

\subsection{Research Questions}
\label{2.1}
The research questions are formulated to structure the review and to capture the key technical and methodological developments in reliability-oriented microgrid design within a rapidly evolving research landscape. Following preliminary literature scoping and engagement with industry-relevant challenges, the refined research questions guiding this review are as follows:
\begin{enumerate}
    \item What are the key factors that govern the planning and design of microgrids to achieve improved reliability?
    \item What are the state-of-the-art optimization techniques and their relative merits for addressing the reliability-oriented microgrid design problem?
    \item Which reliability metrics are commonly employed for microgrid reliability assessment, and how suitable are these metrics for modern renewable-rich microgrids?
\end{enumerate}

While these research questions provide a structured framework for conducting the rapid review, the analysis is not limited to addressing them in isolation. Instead, the review synthesizes insights across these questions to develop an integrated understanding of reliability-oriented microgrid design, including the interactions between planning factors, optimization formulations, and reliability assessment approaches.

\subsection{Search Criteria and Filtering Method}
\label{2.2}

\begin{figure*}[htbp]
    \centering
    \includegraphics[width=0.7\textwidth]{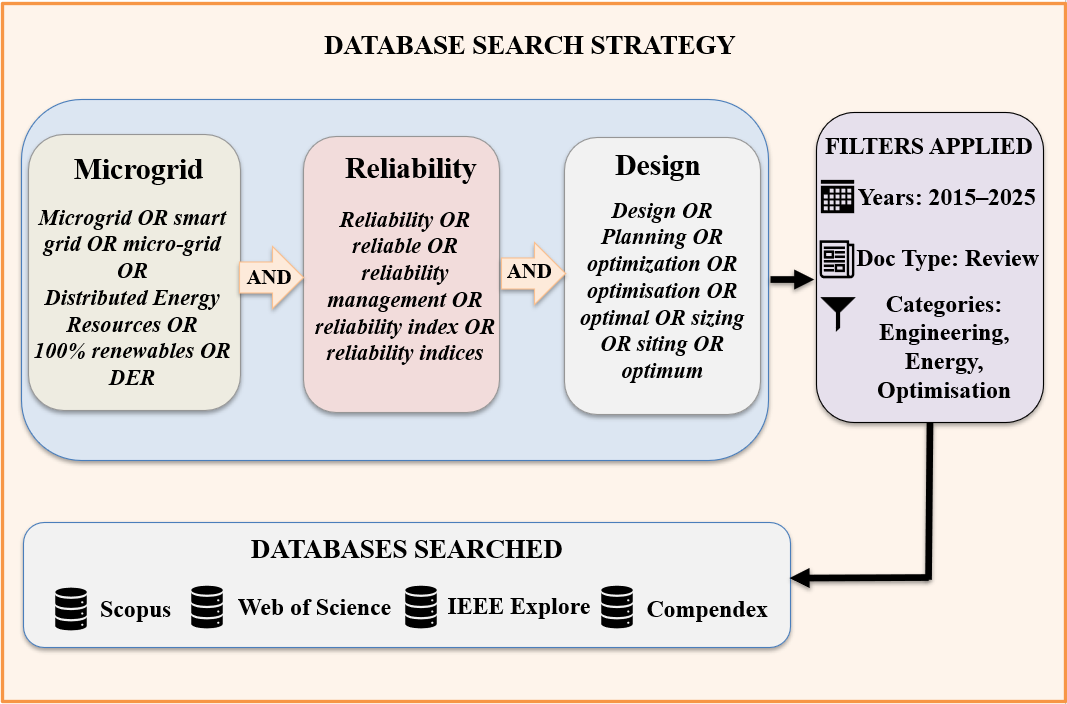}
    \caption{Database search strategy.}
    \label{Fig.2}
\end{figure*}

Guided by the defined research questions, a list of keywords and associated topic-related terms was compiled to construct the search strings as shown in Figure~\ref{Fig.2}. The literature search was conducted on October 8, 2025, across major scientific databases, including Scopus, Web of Science, Compendex, and IEEE Xplore. The search was restricted to publications from 2015 to 2025, with early-access articles from 2026 included, in order to capture recent and relevant developments in microgrid design and reliability-oriented research. To complement the database search, additional searches were performed using Google Scholar and Google Search to identify relevant industry reports and academic publications not indexed in the primary databases. Study selection was guided by predefined eligibility and scope criteria, detailed in Section~\ref{2.3}, with particular emphasis placed on peer-reviewed review articles to ensure a comprehensive synthesis of existing knowledge. The search was limited to English-language publications and filtered using titles, abstracts, and keywords. Database-specific search strings and filtering criteria are reported in Appendix~A to ensure transparency, full reproducibility of the literature search.

\subsection{Screening, Eligibility Criteria and Extraction}
\label{2.3}

In this phase, records retrieved from the database searches were exported to reference management software, where duplicate entries were identified and removed. The remaining articles were screened at the title and abstract level to assess their relevance to the scope of the review. Studies that satisfied the initial screening were subsequently subjected to full-text assessment based on the following eligibility criteria:

\begin{enumerate}
    \item Review articles that provide a comprehensive or systematic overview of microgrid planning and design.
    \item Review studies that examine a broad range of optimization approaches, including conventional optimization methods and emerging AI-based techniques, together with their applications and reported limitations in microgrid design.
    \item Review articles that explicitly address reliability-related aspects of microgrids, including reliability improvement strategies, assessment methodologies, and reliability indices.
\end{enumerate}

\begin{figure*}[htbp] 
    \centering 
    \includegraphics[width=0.6\linewidth]{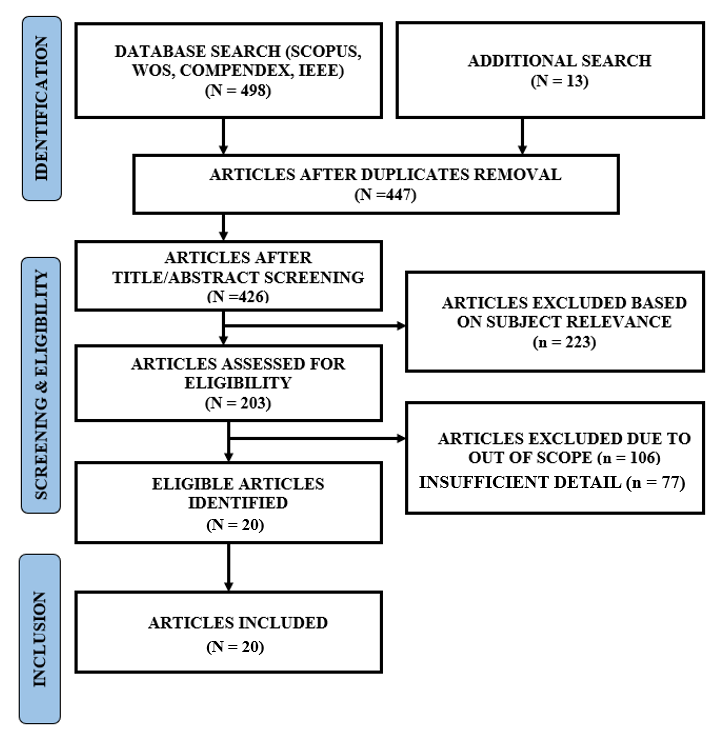} 
    \caption{Search and screening as per PRISMA flow diagrams.} 
    \label{Fig.3}  
\end{figure*}

 From 498 references initially identified, 20 review articles were identified as eligible and included following the screening and selection process guided by the PRISMA framework, as illustrated in Figure~\ref{Fig.3}. Given that the majority of the included reviews did not report statistical or quantitative meta-analyses, a qualitative data extraction and synthesis approach was adopted. For each included review, key bibliographic and methodological attributes were systematically extracted, including publication year, review typology, scope, number of primary studies reviewed, country of origin, funding sources, and any declared conflicts of interest. 

 Despite the transparency of the PRISMA-guided workflow, it is important to acknowledge potential sources of selection bias inherent to this rapid review. First, restricting the search to English-language publications may have excluded relevant studies published in other languages, particularly from regions with extensive microgrid research and deployment activity, such as China, Russia, and Latin America. Second, the primary focus on peer-reviewed review articles, while ensuring a high-level synthesis of existing knowledge, may omit important primary studies or grey literature (e.g., industry technical reports, regulatory documents, and standardization body publications) not captured through the database search. Third, the reliance on four major databases (Scopus, Web of Science, Compendex, and IEEE Xplore), while comprehensive for the engineering domain, may not index all relevant venues, particularly interdisciplinary journals spanning social science, policy, or environmental domains that occasionally address microgrid reliability. Fourth, the focus on recent literature improves timeliness but may reduce coverage of older foundational work unless that work is recovered through citation tracking. 

 To mitigate these potential biases, the database search was complemented by Google Scholar and Google searches. Furthermore, a systematic bidirectional snowball citation-tracking procedure was applied to all 20 included reviews, identifying relevant primary studies that contextualize and validate the review-level findings, thereby reducing the risk of omitted evidence and supporting the robustness of key claims. Nevertheless, the findings should be interpreted with awareness of the inherent trade-off between timeliness and comprehensiveness associated with the rapid review approach.
 
 Building upon this mitigation strategy, 59 underlying primary research studies were incorporated through the bidirectional snowball citation-tracking procedure described above. Selection followed a transparent hierarchical criteria: (i) studies cited by three or more of the 20 included reviews as foundational or methodologically significant work in reliability-oriented microgrid planning, design, or optimization; (ii) studies presenting novel optimization formulations, reliability assessment frameworks, or planning approaches — supported by empirical validation or case study evidence, that address specific technical gaps not covered by multiply-cited works; and (iii) recent publications (2021 - 2025) demonstrating state-of-the-art advances in AI-driven methods, probabilistic planning, or reliability standard that have not yet accumulated multiple review citations due to publication recency. Thematic data extraction subsequently focused on reliability-oriented microgrid planning frameworks, multi-objective optimization approaches addressing cost-reliability trade-offs, reliability-influencing design and operational factors, reliability assessment methodologies and indices, and reported strategies (including AI-driven methods) for reliability enhancement under technical and economic constraints. These themes were synthesized to address the defined research questions and to identify recurring patterns, research gaps, and methodological limitations across the reviewed literature.

\subsection{Quality Assessment of Selected Studies}
\label{2.4}

Assessing the methodological quality and transparency of included reviews is essential for interpreting the evidentiary strength of a rapid review, particularly given the heterogeneous nature of the secondary literature considered. Accordingly, the included review articles were evaluated using AMSTAR 2 (A Measurement Tool to Assess Systematic Reviews, version 2),~\cite{shea2017amstar}, a comprehensive and validated instrument for appraising the methodological rigor of review level evidence . Although originally developed for healthcare research, AMSTAR 2 is applied here in a limited and pragmatic sense as a structured framework for examining review-process reporting, including scope definition, search transparency, study selection, and discussion of evidentiary limitations. Hence, AMSTAR 2 in this study, not adopted as a discipline-specific benchmark for engineering review quality, nor as a proxy for technical novelty or analytical sophistication; rather, as a pragmatic framework for characterizing methodological transparency and the traceability of evidence synthesis across the included literature.

A key aspect of this  adaptation involves reinterpreting the PICO framework (Population, Intervention, Comparator, and Outcome) referenced in Questions~1 and~8 to fit the engineering context. In this context, Population refers to the problem or system under investigation, such as specific microgrid configurations, network topologies, or application contexts; Intervention denotes the methodological approach, including planning frameworks, optimization techniques, or design approaches; Comparator refers to baseline models, benchmark methods, or alternative approaches for comparisons; and Outcome represents the performance-related metrics, including reliability, cost, and other evaluation indicators. In addition, Questions~11,~12, and~15 were excluded, as they pertain to quantitative meta-analysis Such statistical aggregation is not applicable to the present evidence base due to substantial heterogeneity across studies in system configurations, modeling assumptions, temporal resolution, optimization techniques, and performance metrics. Consequently, this review adopts a qualitative narrative synthesis. Table~\ref{AMSTAR2 Table} summarizes the AMSTAR2 assessment results using color-coded indicators, where green (=1) denotes ``yes'', yellow (=0.5) indicates ``partial or unclear'', and red (=0) represents ``no''. It is important to note that these scores should be interpreted solely as relative indicators of reporting completeness and methodological transparency, rather than as absolute measures of technical merit or engineering significance. Their purpose is to provide context regarding how clearly each review communicates its evidentiary basis, rather than to rank their scholarly contribution.

Among the 20 included review articles, four are classified as systematic reviews and sixteen as narrative or scoping reviews. Systematic reviews typically adopt more explicit and reproducible procedures (e.g., structured search strategies and screening logic), whereas narrative/scoping reviews in engineering often prioritize conceptual synthesis and methodological taxonomy. The AMSTAR2 assessment reflects this distinction: a smaller subset of studies demonstrates stronger reporting transparency, particularly in scope definition (Q1), justification of study types (Q3), descriptive synthesis (Q8), and discussion of heterogeneity (Q14), while most studies report these aspects partially or implicitly. This pattern reflects common review conventions in the engineering literature. However, this should not be construed as evidence of poor technical quality or limited scholarly contribution. Instead, they reflect the propensity of engineering reviews to prioritize conceptual synthesis, methodological taxonomy, and design insight over systematic-review protocols.

Finally, it is also crucial to highlight that, consistent with the above role of AMSTAR2, the assessment results were used to inform the interpretation of the evidence base during the synthesis. Specifically, the results influenced the synthesis in two ways. First, reviews demonstrating higher methodological transparency, particularly in scope definition (Q1), study type justification (Q3), descriptive synthesis (Q8), and discussion of heterogeneity (Q14), were given greater interpretive weight when identifying convergent findings across the literature, as their results are supported by more transparent and reproducible procedures. For example, when multiple reviews converge on a finding, this convergence carries greater evidential strength when supported by methodologically transparent reviews. Second, findings derived primarily from reviews with lower transparency were systematically cross checked against the 59 snowballed primary studies before inclusion in the synthesis. This ensures that the conclusions of this rapid review are not overly reliant on evidence with limited methodological reporting. Overall, this approach allows AMSTAR2 to function as a quality weighting and cross validation mechanism rather than a strict inclusion criterion, consistent with its role in this study as a tool for characterizing methodological transparency rather than ranking scholarly merit.

\newcommand{\cellgreen}{\cellcolor{green!25}} 
\newcommand{\cellyellow}{\cellcolor{yellow!50}} 
\newcommand{\cellred}{\cellcolor{red!25}}

\begin{table*}[ht]
\centering
\caption{AMSTAR2 assessment checklist}
\label{AMSTAR2 Table}
\setlength\tabcolsep{8pt}
\begin{tabular}{
    |>{\centering\arraybackslash}m{2cm}
    |>{\centering\arraybackslash}m{0.5cm}
    |>{\centering\arraybackslash}m{0.5cm}
    |>{\centering\arraybackslash}m{0.5cm}
    |>{\centering\arraybackslash}m{0.5cm}
    |>{\centering\arraybackslash}m{0.5cm}
    |>{\centering\arraybackslash}m{0.5cm}
    |>{\centering\arraybackslash}m{0.5cm}
    |>{\centering\arraybackslash}m{0.5cm}
    |>{\centering\arraybackslash}m{0.5cm}
    |>{\centering\arraybackslash}m{0.5cm}
    |>{\centering\arraybackslash}m{0.5cm}
    |>{\centering\arraybackslash}m{0.5cm}
    |>{\centering\arraybackslash}m{0.5cm}
    |>{\centering\arraybackslash}m{0.8cm}|
}
    \hline
        \multirow{2}{*}{Authors} & \multicolumn{13}{c|}{Quality Assessment (AMSTAR2)} & \multirow{2}{*}{\shortstack{Total\\Points}} \\ \cline{2-14} 
        & Q1 & Q2 & Q3 & Q4 & Q5 & Q6 & Q7 & Q8 & Q9 & Q10 & Q13 & Q14 & Q16 & ~ \\ \hline
        
        Ahmadi et al., 2026,~\cite{ahmadi2026comprehensive} & \cellgreen 1 & \cellyellow 0.5 & \cellgreen 1 & \cellgreen 1 & \cellred 0 & \cellred 0 & \cellyellow 0.5 & \cellgreen 1 & \cellred 0 & \cellred 0 & \cellred 0 & \cellgreen 1 & \cellgreen 1 & 7 \\ \hline
        Hadi etal., 2025,~\cite{hadi2025artificial} &  \cellyellow 0.5 & \cellred 0 & \cellred 0 & \cellred 0 & \cellred 0 & \cellred 0 & \cellred 0 & \cellgreen 1 & \cellred 0 & \cellred 0 & \cellred 0 & \cellyellow 0.5 & \cellgreen 1 & 3 \\ \hline
        Ghanbarzadeh et al., 2025,~\cite{ghanbarzadeh2025addressing} & \cellgreen 1 & \cellyellow 0.5 & \cellgreen 1 & \cellgreen 1 & \cellred 0 & \cellred 0 & \cellred 0 & \cellgreen 1 & \cellred 0 & \cellred 0 & \cellred 0 & \cellgreen 1 & \cellgreen 1 & 6.5 \\ \hline
        Malika et al., 2025,~\cite{malika2025critical} & \cellgreen 1 & \cellyellow 0.5 & \cellgreen 1 & \cellgreen 1 & 
        \cellred 0 & \cellred 0 & \cellyellow 0.5 & \cellgreen 1 & 
        \cellred 0 & \cellred 0 & \cellred 0 & \cellgreen 1 & \cellgreen 1 & 
        7 \\ \hline
        Akter et al., 2024,~\cite{akter2024review} & \cellgreen 1 & \cellred 0 & \cellgreen 1 & \cellred 0 & \cellred 0 & \cellred 0 & \cellred 0 & \cellgreen 1 & \cellred 0 & \cellred 0 & \cellred 0 & \cellgreen 1 & \cellgreen 1 & 5 \\ \hline
        M.Thiruna et al., 2023,~\cite{thirunavukkarasu2023comprehensive} & \cellyellow 0.5 & \cellred 0 & \cellred 0 & \cellred 0 & \cellred 0 & \cellred 0 & \cellred 0 & \cellgreen 1 & \cellred 0 & \cellred 0 & \cellred 0 & \cellyellow 0.5 & \cellgreen 1 & 3 \\ \hline
        Alasali et al., 2023,~\cite{alasali2023powering} & \cellyellow 0.5 & \cellred 0 & \cellgreen 1 & \cellred 0 & \cellred 0 & \cellred 0 & \cellred 0 & \cellgreen 1 & \cellred 0 & \cellred 0 & \cellred 0 & \cellyellow 0.5 & \cellgreen 1 & 4 \\ \hline
        De Mel et al., 2022,~\cite{de2022balancing} & \cellgreen 1 & \cellyellow 0.5 & \cellgreen 1 & \cellgreen 1 & \cellred 0 & \cellred 0 & \cellred 0 & \cellgreen 1 & \cellred 0 & \cellred 0 & \cellred 0 & \cellgreen 1 & \cellgreen 1 & 6.5 \\ \hline
        Zhang et al., 2022,~\cite{zhang2022systematic} & \cellgreen 1 & \cellred 0 & \cellgreen 1 & \cellgreen 1 & \cellred 0 & \cellred 0 & \cellred 0 & \cellgreen 1 & \cellred 0 & \cellred 0 & \cellred 0 & \cellgreen 1 & \cellgreen 1 & 6 \\ \hline
        Meera et al., 2022,~\cite{meera2022reliability} & \cellyellow 0.5 & \cellred 0 & \cellred 0 & \cellred 0 & \cellred 0 & \cellred 0 & \cellred 0 & \cellgreen 1 & \cellred 0 & \cellred 0 & \cellred 0 & \cellyellow 0.5 & \cellgreen 1 & 3 \\ \hline
        Polleux et al., 2022,~\cite{polleux2022overview} & \cellgreen 1 & \cellred 0 & \cellgreen 1 & \cellred 0 & \cellred 0 & \cellred 0 & \cellred 0 & \cellgreen 1 & \cellred 0 & \cellred 0 & \cellred 0 & \cellyellow 0.5 & \cellgreen 1 & 4.5 \\ \hline
        Sandelic et al., 2022,~\cite{sandelic2022reliability} & \cellyellow 0.5 & \cellred 0 & \cellgreen 1 & \cellred 0 & \cellred 0 & \cellred 0 & \cellred 0 & \cellgreen 1 & \cellred 0 & \cellred 0 & \cellred 0 & \cellyellow 0.5 & \cellgreen 1 & 4 \\ \hline
        Nazir et al., 2022,~\cite{nazir2021system} & \cellyellow 0.5 & \cellred 0 & \cellred 0 & \cellred 0 & \cellred 0 & \cellred 0 & \cellred 0 & \cellgreen 1 & \cellred 0 & \cellred 0 & \cellred 0 & \cellyellow 0.5 & \cellgreen 1 & 3 \\ \hline
        Liang et al., 2022,~\cite{liang2022planning} & \cellyellow 0.5 & \cellred 0 & \cellred 0 & \cellred 0 & \cellred 0 & \cellred 0 & \cellred 0 & \cellgreen 1 & \cellred 0 & \cellred 0 & \cellred 0 & \cellyellow 0.5 & \cellyellow 0.5 & 2.5 \\ \hline
        Lopez-Prado et al., 2020,~\cite{lopez2020reliability} & \cellyellow 0.5 & \cellred 0 & \cellyellow 0.5 & \cellyellow 0.5 & \cellred 0 & \cellred 0 & \cellred 0 & \cellgreen 1 & \cellred 0 & \cellred 0 & \cellred 0 & \cellgreen 1 & \cellgreen 1 & 4.5 \\ \hline
        Peyghami et al., 2020,~\cite{peyghami2020overview} & \cellgreen 1 & \cellred 0 & \cellred 0 & \cellred 0 & \cellred 0 & \cellred 0 & \cellred 0 & \cellgreen 1 & \cellred 0 & \cellred 0 & \cellred 0 & \cellgreen 1 & \cellgreen 1 & 4 \\ \hline
        Escalera et al., 2018,~\cite{escalera2018survey} & \cellyellow 0.5 & \cellred 0 & \cellgreen 1 & \cellred 0 & \cellred 0 & \cellred 0 & \cellred 0 & \cellgreen 1 & \cellred 0 & \cellred 0 & \cellred 0 & \cellgreen 1 & \cellgreen 1 & 4.5 \\ \hline
        Ehsan et al., 2018,~\cite{ehsan2018optimal} & \cellyellow 0.5 & \cellred 0 & \cellred 0 & \cellred 0 & \cellred 0 & \cellred 0 & \cellred 0 & \cellgreen 1 & \cellred 0 & \cellred 0 & \cellred 0 & \cellgreen 1 & \cellyellow 0.5 & 3 \\ \hline
        Prakash et al., 2016,~\cite{prakash2016optimal} & \cellyellow 0.5 & \cellred 0 & \cellred 0 & \cellred 0 & \cellred 0 & \cellred 0 & \cellred 0 & \cellgreen 1 & \cellred 0 & \cellred 0 & \cellred 0 & \cellyellow 0.5 & \cellgreen 1  & 3 \\ \hline
        Gamarra et al., 2015,~\cite{gamarra2015computational} & \cellyellow 0.5 & \cellred 0 & \cellgreen 1 & \cellred 0 & \cellred 0 & \cellred 0 & \cellred 0 & \cellgreen 1 & \cellred 0 & \cellred 0 & \cellred 0 & \cellyellow 0.5 & \cellred 0 & 3 \\ \hline

\end{tabular}
\normalsize
\end{table*}

\subsection{Study Characteristics}
\label{2.5}

The 20 reviewed articles, published between 2015 and early 2026, collectively examined 1,773 sources and reflect contributions from diverse geographical regions. These reviews form the core synthesis corpus of this study and are integrated throughout subsequent sections. Complementing this, 59 primary technical studies were curated via forward and backward citation tracking (Section~\ref{2.3}) to develop a methodological taxonomy of reliable microgrid design. These snowballed studies provide a structured framework to guide researchers and practitioners through the diverse landscape of reliability-oriented microgrid design. As shown in Table 2, the reviewed literature spans three core dimensions: design and planning, optimization methodologies, and reliability assessment frameworks. While planning-oriented studies dominate, many works address optimization and explicitly evaluate reliability, highlighting the field’s multidimensional nature. Notably, eight articles~\cite{ahmadi2026comprehensive,sandelic2022reliability,malika2025critical, meera2022reliability, polleux2022overview,alasali2023powering, ghanbarzadeh2025addressing,nazir2021system} provide comprehensive analyses of reliability in microgrid design, reflecting its growing importance as a central planning criterion. However, very few studies have addressed the key factors that guide the design of reliable microgrids. A clear dichotomy also emerges between long-term infrastructure planning and short-term operational optimization: studies such as Sandelic et al.~\cite{sandelic2022reliability} and Ghanbarzadeh et al.~\cite{ghanbarzadeh2025addressing} adopt macro-temporal perspectives focused on resource sizing and expansion planning, whereas others, including Hadi et al.~\cite{hadi2025artificial} and Zhang et al.~\cite{zhang2022systematic} emphasize short-term operational control, energy management, and scheduling.

Overall, the included studies indicate that designing a reliable microgrid requires a multifaceted framework informed by diverse research perspectives and theoretical viewpoints. For example, Zhang et al.~\cite{zhang2022systematic} and Prakash et al.~\cite{prakash2016optimal} approached microgrid optimization from different perspectives. Theoretical considerations are emphasized in~\cite{zhang2022systematic}, particularly the trade-off between model complexity and solution accuracy, while~\cite{prakash2016optimal} adopted a more application-oriented perspective by examining optimal sizing and siting strategies for distributed generators within the network. Building on this practical planning viewpoint, Malika et al.~\cite{malika2025critical} further investigated sizing and siting methodologies by highlighting demand uncertainty, conditions monitoring requirements, and infrastructure integration challenges. Complementing these works, Thirunavukkarasu et al.~\cite{thirunavukkarasu2023comprehensive} reviewed the use of commercial software tools in microgrid optimization studies, emphasizing the challenges associated with translating theoretical formulations into implementable solutions. The reviewed literature demonstrates a clear temporal evolution in methodological preferences. Earlier reviews, including Prakash et al.~\cite{prakash2016optimal} (2016), Ehsan et al.~\cite{ehsan2018optimal} (2018), and Escalera et al.~\cite{escalera2018survey} (2018), emphasize mathematical programming and metaheuristic approaches, reflecting the dominant paradigm of that period. In contrast, reviews published after 2022 increasingly incorporate AI-driven methodologies, with~\cite{hadi2025artificial,ahmadi2026comprehensive,malika2025critical} devoting substantial attention to machine learning, deep learning, and reinforcement learning applications.

In addition to optimization formulation and techniques, ensuring reliability is also crucial. Accordingly, Sandelic et al.~\cite{sandelic2022reliability} focused on the wear and tear of power electronics in microgrids, emphasizing that even the most rigorously formulated theoretical models are susceptible to failure if the reliability of the associated physical components is not considered. In contrast, Ghanbarzadeh et al.~\cite{ghanbarzadeh2025addressing} focused on the reliability implications of high renewable penetration from a system planning perspective. Their introduction of the "\textit{Dunkelflaute}" phenomenon, which refers to extended periods of minimal or no wind and solar generation, represents a critical reliability challenge that is absent from earlier reviews focused on component-level reliability. This macro-level perspective on system adequacy further complements the micro-level focus of Peyghami et al.~\cite{peyghami2020overview} on power electronics reliability, together providing a more complete picture of the reliability challenges spanning from component to system levels.

As reliability challenges in microgrids continue to grow in complexity, artificial intelligence (AI) has recently emerged as a complementary approach to support the design of reliable microgrids, as mentioned earlier. Hadi et al.~\cite{hadi2025artificial} present a structured analysis of AI applications throughout the entire life cycle of a microgrid, encompassing design, control, and maintenance. Their work identifies critical research gaps, particularly the lack of an end-to-end perspective on AI in microgrids and the limited discussion on AI for microgrid sizing. However, Ahmadi et al.~\cite{ahmadi2026comprehensive} complemented this by providing a unified framework for real-time microgrid management and fault detection, systematically evaluating machine, deep, and reinforcement learning techniques across various tasks. Their emphasis on hybrid models that combine AI with optimization strategies addresses practical deployment challenges including scalability, computational complexity, and adaptability.

Despite the growing role of AI in supporting microgrid design and operation, rigorous reliability evaluation frameworks are essential for assessing system performance and validating design decisions. Regarding reliability evaluation methodologies, López-Prado et al.~\cite{lopez2020reliability} reported a near-equal distribution of use of Monte Carlo Simulation (MCS), and analytical approaches to model real-world scenarios. Escalera et al.~\cite{escalera2018survey} similarly acknowledged the computational advantages of analytical techniques while recognizing the superiority of MCS in capturing stochastic behavior. In terms of reliability metrics, customer-oriented interruption indices remain the most widely adopted across the reviewed studies~\cite{lopez2020reliability}, while energy-based metrics are employed for economic assessments~\cite{ghanbarzadeh2025addressing}. Nazir et al.~\cite{nazir2021system} further introduced two additional performance indicators and examined architectural and modeling aspects of microgrids, emphasizing their implications for distribution-level reliability.

When reliability targets cannot be fully achieved through planning and operational measures alone, protection strategies play a critical role in mitigating faults and limiting their impact on overall system performance. In this context, a key insight emerging from the comparative analysis is the ongoing debate surrounding the role of artificial intelligence in microgrid protection. Ahmadi et al.~\cite{ahmadi2026comprehensive} and Hadi et al.~\cite{hadi2025artificial} characterize AI as a potentially transformative approach capable of addressing limitations inherent to conventional protection schemes. At the same time, both studies acknowledge substantial implementation challenges. Ahmadi et al. identify the simulation-to-reality gap as a fundamental barrier to real-world deployment, while Hadi et al. note that current AI-based protection solutions remain fragmented across control and maintenance stages. Complementing this perspective, Meera et al.~\cite{meera2022reliability} demonstrate that traditional overcurrent protection schemes can suffer from reduced sensitivity and coordination failures when fault current levels differ markedly between grid-connected and islanded operating modes. This protection-reliability nexus is further examined by Alasali et al.~\cite{alasali2023powering}, who highlight similar coordination limitations and advocate communication-assisted protection schemes as a means to preserve reliability under dynamic network topologies.

In summary, across the reviewed corpus, a clear and analytically significant pattern emerges: while the literature increasingly recognizes reliability as a central planning concern, it remains methodologically fragmented in how that concern is represented and operationalized. Planning-oriented reviews tend to emphasize resource adequacy, optimal siting, and long-term topology decisions, treating reliability primarily as a sizing/siting constraint rather than a governing design objective. Optimization-oriented reviews, whereas, focus more strongly on algorithmic tractability and solution quality, with reliability typically reduced to a penalty term or threshold condition within the objective function, rarely examining whether the underlying reliability representation reflects the operational realities of inverter-dominated, storage-dependent systems. Reliability-focused reviews, in contrast, provide deeper treatment of indices and failure mechanisms but frequently remain disconnected from concrete planning decisions, leaving an assessment-planning gap that limits their actionability for practitioners. Notably, no single reviewed article spans all three strands in a unified and mutually consistent framework. The central analytical gap is therefore not a shortage of literature, but the absence of integration among these strands, an absence that risks producing designs that are algorithmically efficient, adequately sized, and formally assessed, yet still unreliable in operation because the three concerns were never resolved jointly. This fragmentation directly motivates the integrated, reliability-as-design-principle framing adopted in the present review. The subsequent sections systematically examine these dimensions in order to address the research question posed in this review.

\clearpage 
\onecolumn 

\begin{longtable}{|>{\centering\arraybackslash}p{1.5cm}|>{\centering\arraybackslash}p{1.5cm}|>{\centering\arraybackslash}p{0.9cm}|>{\centering\arraybackslash}p{2.4cm}|>{\centering\arraybackslash}p{1.68cm}|>{\centering\arraybackslash}p{1.2cm}|>{\centering\arraybackslash}p{2.0cm}|>{\centering\arraybackslash}p{1.5cm}|>{\centering\arraybackslash}p{1.5cm}|}
\caption{Characteristics of the 20 included review articles} \\ 
\hline
\textbf{First Author., Year} & \textbf{Country} & \textbf{Review Type} & \textbf{Study Scope} & \textbf{Topics Covered} & \textbf{Cited studies} & \textbf{Funding} & \textbf{Conflicts of Interest} & \textbf{Snowballed Primary Studies} \\ \hline
\endfirsthead 

\multicolumn{9}{c}%
{{\bfseries \tablename\ \thetable{} -- continued from previous page}} \\
\hline 
\textbf{First Author., Year} & \textbf{Country} & \textbf{Review Type} & \textbf{Study Scope} & \textbf{Topics Covered} & \textbf{Cited studies} & \textbf{Funding} & \textbf{Conflicts of Interest} & \textbf{Snowballed Primary Studies} \\ \hline
\endhead 

\hline \multicolumn{9}{|r|}{{Continued on next page}} \\ \hline
\endfoot 

\hline
\endlastfoot 

Ahmadi et al., 2026,~\cite{ahmadi2026comprehensive} & Canada & SLR & Reviews AI-driven techniques for microgrid reliability & AI, Reliability & 318 & No funding acknowledged & None &~\cite{islam2024improving,gali2023adaptive,khalid2024smart,hu2023towards,chen2021multi}  \\ \hline
Malika et al., 2025,~\cite{malika2025critical} & India, Czech Republic & SLR & Reviews optimal placement and sizing of DGs in microgrid & Planning, Optimization & 234 & EU (REFRESH Project) & None &~\cite{adefarati2017reliability,opathella2020milp,agajie2022impact,xie2024optimal,jaleel2021reliability,nick2017optimal,elkadeem2019optimal,yaghoubi2021optimal}  \\ \hline
Hadi et al., 2025,~\cite{hadi2025artificial} & France, Italy, Spain & NR & Reviews AI applications in microgrid design & Planning, Reliability & 187 & ISEN Yncréa Ouest & None &~\cite{gali2023adaptive,chrif2025techno,ramli2018optimal,zia2019energy}  \\ \hline
Ghanbarza deh et al., 2025,~\cite{ghanbarzadeh2025addressing} & Australia & SLR & Reviews reliability challenges in planning under high renewables & Planning, Reliability & 163 & Edith Cowan University & None &~\cite{islam2024improving,adefarati2017reliability,adefarati2019reliability,heylen2018review,firouzi2022reliability,mayer2023probabilistic}  \\ \hline
Akter et al., 2024,~\cite{akter2024review}  & Canada, Bangladesh, Australia & NR & Reviews scopes and trends of microgrid optimization & Optimization & 200 & Canada National Sciences and Engineering Research Council & None &~\cite{gali2023adaptive,nojavan2017efficient,harasis2020reliable,kamal2023optimal,mashayekh2017mixed,jafari2020optimal}  \\ \hline
M.Thiruna et al., 2023,~\cite{thirunavukkarasu2023comprehensive} & India & NR & Reviews optimization techniques & Optimization & 305 & No funding acknowledged & None &~\cite{moghaddam2019designing,ghorbani2018optimizing}  \\ \hline
Alasali et al., 2023,~\cite{alasali2023powering}  & Libya, Jordan, Saudi Arabia & NR & Reviews protection techniques for reliable microgrid & Protection, Reliability & 227 & King Khalid University, Saudi Arabia  & None &~\cite{ganjian2019seamless,sardari2018enhancement,ghotbi2020design} \\ \hline
De Mel et al., 2022~\cite{de2022balancing} & UK & SLR & Discusses the trade-off between model accuracy and complexity in microgrids and OPF. & Optimization & 129 & No funding acknowledged & None &~\cite{hittinger2015evaluating,mashayekh2017mixed,huang2017optimal,karmellos2019multi}\\ \hline
Zhang et al., 2022~\cite{zhang2022systematic}  & Australia & SLR & Reviews integration challenges and solutions for renewable DGs and energy storage. & Planning, Optimization & 103 & Murdoch University & None &~\cite{adefarati2017reliability,nick2017optimal,huang2017optimal,lazzeroni2019optimal,wang2019joint} \\ \hline
Nazir et al., 2022~\cite{nazir2021system}  & India, USA & NR & Discusses microgrid architecture and its impact on reliability. & Planning, Reliability & 100 & No funding acknowledged & None &~\cite{zia2019energy} \\ \hline
Sandelic et al., 2022~\cite{sandelic2022reliability}  & Denmark & NR & Examines power electronics' role in microgrid reliability. & Power Electronics, Reliability & 148 & Villum Foundation & None &~\cite{islam2024improving,hittinger2015evaluating,rosales2019microgrids,ganjian2019seamless,jafari2020optimal,peyghami2020incorporating,adefarati2017reliability,amir2019reliability,chen2021multi,adefarati2019reliability,nojavan2017efficient,arefifar2013optimum,arefifar2014dg,khodaei2014microgrid,mitra2016optimal,farzin2017role,baghaee2016reliability,wang2014optimal,cortes2017microgrid,che2016optimal,cagnano2020microgrids} \\ \hline
Meera et al., 2022~\cite{meera2022reliability} & India & NR & Reviews the impact of renewable DG on network reliability. & Planning, Reliability & 126 & No funding acknowledged & None &~\cite{adefarati2017reliability,amir2019reliability,elkadeem2019optimal,yaghoubi2021optimal,adefarati2019reliability,arefifar2013optimum,farzin2017role,pirouzi2022hybrid,arefifar2014dg,battu2015dg,canales2021cost,mitra2016optimal,farzin2016reliability,baghaee2016reliability,xu2015evaluation,de2016modeling} \\ \hline
Polleux et al., 2022~\cite{polleux2022overview}  & France & NR & Discusses challenges in integrating renewables into microgrids. & Planning, Reliability & 192 & Total Energies S.E & Financial interests/personal relationships &~\cite{rosales2019microgrids,ganjian2019seamless,khodaei2014microgrid,cagnano2020microgrids} \\ \hline
Liang et al., 2022,~\cite{liang2022planning}  & Canada & NR & Reviews the role of microgrid reconfiguration & Reliability, Reconfiguration & 92 & University of Saskatchewan & None &~\cite{rosales2019microgrids,arefifar2013optimum,khodaei2014microgrid,mitra2016optimal,mohamed2018efficient,yazdavar2020optimal,cortes2017microgrid,che2016optimal} \\ \hline 
Lopez-Prado et al.,  2020~\cite{lopez2020reliability}  & Colombia & NR & Reviews microgrid reliability metrics and methods. & Reliability & 156 & No funding acknowledged & None &~\cite{adefarati2017reliability,arefifar2013optimum,arefifar2014dg,farzin2017role,de2016modeling,nikmehr2016reliability,wang2014optimal,farzin2016reliability,xu2015evaluation} \\ \hline
Peyghami et al., 2020~\cite{peyghami2020overview} & Denmark, Netherlands & NR & Discusses reliability challenges in modern power electronic-based microgrids & Reliability, Power Electronics & 68 & Villum Fonden, Denmark & None &~\cite{peyghami2020incorporating}\\ \hline
Escalera et al., 2018~\cite{escalera2018survey} & Spain, Ireland & NR & Surveys methodologies for assessing reliability. & Reliability & 112 & Spanish Ministry of Economy and Competitiveness project RESmart  & None &~\cite{adefarati2017reliability,arefifar2013optimum,arefifar2014dg,mitra2016optimal,baghaee2016reliability,wang2014optimal,xu2015evaluation,farzin2016reliability,nikmehr2016reliability,awad2014optimal}\\ \hline
Ehsan et al., 2018~\cite{ehsan2018optimal}  & China & NR & Reviews analytical techniques for optimal planning & Planning, Optimization & 97 & National High Technology Research and Development Program, China & None &~\cite{awad2014optimal} \\ \hline
Prakash et al., 2016~\cite{prakash2016optimal}  & India & NR & Reviews optimal sizing and siting techniques & Planning, Optimization & 112 & No funding acknowledged & None &~\cite{battu2015dg,dehghanian2013optimal}\\ \hline
Gamarra et al., 2015~\cite{gamarra2015computational}  & Spain, Denmark & NR & Reviews optimization techniques in microgrid planning. & Planning, Optimization & 125 & No funding acknowledged & None &~\cite{hittinger2015evaluating,khodaei2014microgrid} \\ \hline
\end{longtable}

\twocolumn

\section{Key Factors in Designing a Reliable Microgrid}
\label{Sec:3}

This rapid review defines microgrids as localized aggregation of distributed energy resources (DERs) \nomenclature[A]{DER}{distributed energy resource} such as solar panels, batteries, and generators that act as a unified power system capable of operating autonomously or in conjunction with the main grid. Accordingly, identifying clear design objectives upfront is pivotal when planning such a complex system to ensure the anticipated advantages or mitigate potential challenges. These objectives can range from reducing environmental impact~\cite{zhong2022optimal} and minimizing costs~\cite{uddin2023techno} to improving energy dispatch capabilities~\cite{shen2023multi}, integrating more renewables~\cite{hasan2023critical} or enhancing reliability~\cite{garg2021design} and resiliency~\cite{gargari2021preventive}. It is important to note that reliability and resilience, while related, are conventionally distinguished by the nature of the events they address. Reliability refers to the ability of a microgrid to continuously supply electricity under normal operating conditions, whereas resilience concerns its capacity to withstand and recover from high impact, low probability exogenous events such as extreme weather, natural disasters, or cyber-physical disruptions. However, this boundary becomes less distinct in renewable-based microgrids. Prolonged scarcity events, such as \textit{Dunkelflaute}, can be interpreted as reliability constraints related to adequacy, rather than resilience challenges, because they result from normal stochastic variability rather than external shocks. Accordingly, this review focuses on cost-optimal reliable microgrid design, treating adequacy under extended scarcity as a reliability concern intrinsic to renewable-dominated systems. Resilience to exogenous disruptions, including natural disasters, equipment destruction, and cyber attacks, remains outside the scope of this review.

While most reviews included in this article extensively explore the microgrid design process, which includes various steps such as benchmark compilation, element modeling, case study development, etc., they largely overlook the key factors that should guide overall planning and decision-making. This oversight represents a critical gap, as these key factors, from technical considerations to economic implications, ultimately determine the success of any microgrid project throughout its life cycle. 

Notably, only~\cite{sandelic2022reliability} and~\cite{liang2022planning} offered a limited discussion on some factors and criteria that are considered in the design process. However, Sandelic et al.~\cite{sandelic2022reliability} cited two notable papers~\cite{8746836,farhangi2019microgrid} that could be of interest while defining key design factors for reliability-oriented microgrids. Drawing insights from these studies, the following subsections explore various key factors for a reliable microgrid design.

\subsection{Technical Factors}
\label{3.1}
\subsubsection{Load Analysis}
\label{3.1.1}
Every consumer is different when it comes to power consumption; however, all consumers require a continuous and reliable electricity supply. This load variation can be observed in many ways, such as the type of load (critical or non-critical), sector (residential, commercial, or industrial), and specific operational demands. Therefore, understanding the consumer demand profile is a fundamental factor in microgrid planning~\cite{naderipour2020optimal}. Lopez-Prado et al.~\cite{lopez2020reliability} highlight the critical role of accurate load forecasting while designing microgrids for reliability, as inadequate forecasting risks over/under-sizing of generation and storage components. For instance, in sunny regions, precise load analysis enables the efficient planning of storage systems to capture surplus solar energy during peak hours for dispatching during high-demand periods, thereby maintaining reliability.

Similarly, incorporating dynamic load models into the design is crucial to ensure reliable microgrid operation, especially in real-world scenarios where load conditions fluctuate due to various factors~\cite{liang2022planning}. This concern is further reinforced by Malika et al.~\cite{malika2025critical}, who argue that emerging load behaviors, such as electric vehicle charging and nonlinear demand profiles, significantly complicate the optimal placement and sizing of distributed generation and storage resources in microgrids. Furthermore, Sandelic et al.~\cite{sandelic2022reliability} addressed the challenges associated with data acquisition for load forecasting models, as the quality and availability of data can affect the reliability of load forecasting models. Hadi et al.~\cite{hadi2025artificial} advance this discussion by demonstrating that AI driven load forecasting approaches can mitigate such data quality challenges by learning complex patterns from historical data. In particular, deep learning models, such as long short-term memory (LSTM) networks, exhibit superior performance in capturing temporal dependencies that conventional forecasting techniques often fail to represent. This represents a methodological advancement compared to the purely statistical approaches emphasized in other reviews~\cite{kondaiah2022review,khodaei2014microgrid}, as AI based models are capable of capturing complex and nonlinear relationships in load behavior. However, Ahmadi et al.~\cite{ahmadi2026comprehensive} note that the challenges related to computational burden and the gap between simulated models and real world operating conditions remain significant barriers to large scale deployment. Furthermore, the appropriate selection of load data resolution has also been identified as a key planning consideration by Sandelic et al.~\cite{sandelic2022reliability}, suggesting that short-term planning should rely on hourly or daily load profiles, while long-term planning requires the use of yearly data to capture broader demand trends.

\subsubsection{Generation Mix: Balancing Conventional and Renewable Sources}
\label{3.1.2}

Determining the optimal generation mix is equally crucial as load analysis for long-term microgrid reliability~\cite{sandelic2022reliability}, as it enables the system to withstand different uncertainties, including changing weather patterns~\cite{dashtaki2023optimal}, volatile fuel prices~\cite{ruiz2017optimal}, and grid outages~\cite{alramlawi2024chance}. DERs like solar and wind power introduce both opportunities and challenges in this context. While they offer low-loss localized generation, their inherent variability~\cite{zhang2022systematic, liang2022planning} challenges reliable supply. Hence, a balanced mix of conventional and renewable generation sources, where variable renewable outputs are supported by more stable sources like diesel generators, is crucial. However, the ideal combination requires proper optimization techniques, which are detailed in Section 4. 

While the ideal generation combination is crucial for microgrid reliability, it is also important to integrate the nonlinear characteristics of these generation sources during the design phase, as overlooking or simplifying these nonlinearities can lead to improper system sizing and compromise real-world performance reliability~\cite{de2022balancing}. Therefore, rigorously modeling nonlinearities of generation sources should be a key consideration for planners to ensure reliability. Furthermore, Prakash et al.~\cite{prakash2016optimal} highlighted the strategic placement of these generation sources, while~\cite{nazir2021system,ehsan2018optimal} recommended considering scalability and flexibility in the design of the generation mix, allowing easy modifications to meet evolving energy demands.

Despite these well-established design principles, recent studies suggest that even carefully optimized and flexible generation mixes may fall short when confronted with prolonged renewable scarcity events. In particular, Ghanbarzadeh et al.~\cite{ghanbarzadeh2025addressing} argued that existing approaches inadequately account for the \textit{Dunkelflaute} phenomenon, defined as extended periods of simultaneously low wind and solar generation. Such events can persist for days or even weeks in certain regions and introduce high-impact adequacy risks that may not be fully addressed through generation-mix optimization alone. Addressing \textit{Dunkelflaute} events  therefore requires explicit representation of correlated, persistent renewable scarcity within planning frameworks, extending beyond conventional assessments based on average conditions. In practice, this can be implemented by constructing \textit{Dunkelflaute}-aware scenarios that preserve persistence and cross-resource dependence. For example, planners may cluster historical or reanalysis weather into regime states (e.g., synoptic patterns) and sample multi-day sequences using observed transition frequencies, or alternatively run multi-year sequential simulations to retain autocorrelation in wind/solar and their joint dependence with demand. Such regime- or sequence-based scenario sets can then be embedded in reliability evaluation and in planning optimization by enforcing survivability-style constraints (e.g., minimum autonomy hours at high-percentile scarcity) alongside cost objectives.

This perspective also redefines the role of energy storage and backup resources within the portfolio. Whereas studies such as Polleux et al.~\cite{polleux2022overview} and Escalera et al.~\cite{escalera2018survey} primarily consider energy storage systems for load shifting and frequency regulation, Ghanbarzadeh et al.~\cite{ghanbarzadeh2025addressing} emphasize their function as strategic reserves capable of sustaining supply during prolonged renewable droughts. Accordingly, storage duration requirements should be derived from the tail of the renewable scarcity distribution rather than average shortfall characteristics. For example, the 95th or 99th percentile duration of consecutive periods in which renewable generation falls below a context-specific threshold (e.g., 10--20\% of rated capacity) can be used to infer the minimum discharge duration required to bridge extended scarcity without curtailing critical loads. 

From a modeling standpoint, sequential Monte Carlo simulation used for reliability evaluation ideally spans multi-year horizons (e.g., a decade or more where data permit) to capture interannual variability in extreme scarcity events, whereas single-year or "typical meteorological year" approaches can under-represent their frequency and persistence. Beyond sizing, storage operating policies should incorporate forecast-informed reserve targets that raise state-of-charge set-points when extended low-renewable periods are anticipated, rather than optimizing purely for economic arbitrage that may deplete reserves immediately before scarcity. Complementing this perspective, Malika et al.~\cite{malika2025critical} emphasize that mitigating long-duration adequacy challenges motivates evaluating emerging long-duration storage technologies (e.g., flow batteries and hydrogen-based systems) alongside conventional options due to their ability to provide extended discharge durations.

Taken together, it is evident that reliability-oriented generation planning in renewable-based microgrids extends well beyond installed capacity. What ultimately governs supply security is whether the resource portfolio, in conjunction with storage dispatch and operational policy, can sustain demand continuously across time under correlated renewable scarcity. This temporal dimension is critical: capacity that appears sufficient in aggregate may prove inadequate during prolonged periods of simultaneously low wind and solar generation, when storage reserves are depleted and no single resource can compensate independently. Accordingly, capacity expansion decisions should be evaluated not only by their nominal power contribution, but by their ability to maintain feasible operating conditions across successive time periods under adequacy stress. This reframing has direct implications for the optimization frameworks used in microgrid planning, as methods that operate on aggregated or average representations of renewable availability are structurally unable to capture these dynamics, motivating a shift toward temporally resolved, uncertainty-aware planning approaches examined in the following sections.

\subsubsection{Grid Connectivity }
\label{3.1.3}

Grid connectivity is another crucial factor in reliable microgrid design, serving as a backup by connecting to the main utility grid during local generation shortfalls or maintenance~\cite{guichi2021optimal}. However, the decision between grid-connected and off-grid configurations should be tailored based on a holistic analysis of various elements~\cite{sandelic2022reliability}, including infrastructural readiness, economic feasibility, regional energy market dynamics, anticipated load demands, environmental implications, and the broader socio-cultural aspirations of the served communities~\cite{polleux2022overview}. Although the general trend suggests off-grid configurations for remote areas and grid-connected systems for urban regions~\cite{sandelic2022reliability}, planners should not take it as a strict norm. Each region presents its unique challenges and opportunities which require a thorough  evaluation to determine the most suitable approach. Furthermore, implementing grid connectivity requires proper coordination and control mechanisms to ensure seamless integration and operation between the microgrid and the main grid~\cite{zhu2021optimal}. Power electronics devices play a vital role in this regard, as they can significantly improve system reliability by facilitating the connection between the distribution network and the microgrid ~\cite{liang2022planning}. Furthermore, grid connectivity has been shown to improve reliability while minimizing costs, as demonstrated by Meera et al.~\cite{meera2022reliability}, who reported that maximum reliability with minimum cost is achieved when a microgrid is connected to an upstream network. 

\subsubsection{Energy Storage System}
\label{3.1.4}

Similar to grid connectivity, energy storage systems (ESS) \nomenclature[A]{ESS}{Energy Storage System}play a crucial role in ensuring the reliability of microgrids and should be a key consideration for planners during the design phase. Microgrids with ESS can store surplus energy during low demand and release it during peak hours or outages, enhancing both reliability and resilience, as highlighted by Polleux et al.~\cite{polleux2022overview}. The authors emphasized that ESS could serve as a buffer to compensate for unplanned power fluctuations and provide system services in the event of unplanned events. Similarly, Escalera et al.~\cite{escalera2018survey} mentioned the versatility of ESS in load shifting, congestion easing, frequency control, voltage regulation, and electricity trading. However, they also emphasized that the planning stage must consider the location, size, and specific functional details (type, lifetime, and ramps to charge/discharge) of the ESS. However, Sandelic et al.~\cite{sandelic2022reliability} expressed concern regarding the high upfront cost of ESS, as it affects the overall economics of the system. Extending this discussion beyond cost considerations, Ghanbarzadeh et al.~\cite{ghanbarzadeh2025addressing} emphasized that absolute reliability cannot be guaranteed due to component failures and uncertainties in renewable generation and demand, positioning ESS as a fundamental enabler of reliability rather than a discretionary asset in renewable-dominated microgrids.

Although ESS are increasingly recognized as indispensable reliability enablers in renewable-dominated microgrids, their effective contribution critically depends on how they are modeled and evaluated over their operational lifetime. In this regard, Malika et al.~\cite{malika2025critical} highlighted that a key aspect often overlooked in microgrid planning is the explicit incorporation of degradation behavior and replacement of ESS cycles in life-cycle cost assessments. They also argued that static assessments of ESS economics fail to capture the dynamic nature of battery degradation, which varies significantly depending on the operating conditions, depth of discharge patterns, and temperature profiles. Therefore, microgrid planners should move beyond simple capacity optimization to integrated ESS in planning frameworks and explicitly consider degradation behavior, lifecycle economics, and long-term reliability impacts.

In summary, the reviewed literature suggests several key design implications for ESS in microgrids. First, ESS sizing should be guided not only by average energy balancing requirements but also by extreme operating conditions, ensuring sufficient capacity and duration to maintain reliability during prolonged supply deficits. Second, beyond sizing, spatial placement of ESS within the network should be strategically optimized to alleviate congestion and enhance local reliability. Third, the selection of storage technology and its operational characteristics (e.g., charge–discharge rates, efficiency, and lifetime) must be aligned with the specific services required, including peak shaving, frequency regulation, and resilience support. Finally, planning frameworks should explicitly incorporate degradation dynamics and replacement strategies, recognizing ESS as a long-term, evolving asset rather than a static component. Collectively, these considerations highlight that effective ESS deployment requires an integrated, lifecycle-oriented design approach that balances reliability, operational flexibility, and economic performance.

\subsection{Reliability \& Security Factors}
\label{3.2}
\subsubsection{Protection Measures}
\label{3.2.1}

Another critical factor in the design of reliable microgrids is the implementation of robust protection measures. The primary objective of protection is to ensure the safety of the equipment while maintaining the efficiency and reliability of the electrical network~\cite{alasali2023powering}. However, traditional overcurrent protection schemes, commonly used in standard distribution grids~\cite{nazir2021system}, are often inadequate for microgrids due to bidirectional power flow caused by high penetration of distributed generation (DG) \nomenclature[A]{DG}{Distributed Generation}penetration~\cite{meera2022reliability}. In addition, challenges such as false tripping~\cite{larik2022comprehensive}, fluctuating fault currents~\cite{augustine2020fault}, unintended islanding~\cite{ahangar2020review}, and out-of-synchronism reclosing~\cite{patnaik2020ac} complicate protection coordination and often necessitate modifications to network topology. These characteristics make protection a non-trivial design consideration in microgrids.

From a planning perspective, these challenges imply that protection system performance directly constrains the feasible operating envelope of reliability-oriented microgrid design. Critical functionalities such as successful islanding, rapid fault isolation, effective post-fault reconfiguration, and reliable transitions between grid-connected and islanded modes depend on protection capabilities that must be defined at the design stage. Consequently, protection-related decisions, including device placement and coordination parameters (e.g., switch locations, recloser settings, sectionalizer placement, and communication infrastructure requirements), should be treated as co-optimized planning variables alongside generation sizing and network topology, rather than as downstream operational considerations. Failure to integrate protection at the planning stage can result in system designs that appear adequate from a generation–storage perspective but are operationally unprotected, leading to degraded reliability in real-world deployment.

Despite its importance, protection has often been overlooked in microgrid design studies, with many works implicitly assuming that protection devices operate reliably under all conditions~\cite{meera2022reliability}. This assumption neglects several factors that can significantly influence system reliability, including protection device failures, backup protection requirements, maintenance outages, and the possibility of concurrent failures. In addition, undesired activation of protection devices, exposure to varying operating environments, and incorrect protection system operation introduce further uncertainty that is rarely captured in planning models. These limitations highlight the need for a more comprehensive representation of protection behavior within reliability-oriented microgrid planning.

To address these challenges, several studies have proposed enhancements to conventional protection schemes. For instance, Escalera et al.~\cite{escalera2018survey} advocated the use of advanced control strategies for protection elements, demonstrating their ability to reduce both the frequency and duration of interruptions. They also recommended the deployment of automatic reclosers and tele-controlled switches to improve system resilience. Similarly, Lopez-Prado et al.~\cite{lopez2020reliability} emphasized the role of fault detection and diagnosis systems (FDDS)\nomenclature[A]{FDDS}{Fault Detection and Diagnosis Systems} and remote-controlled disconnectors, while Alasali et al.~\cite{alasali2023powering} proposed decentralized backup protection and the integration of cloud-based architectures. Although these approaches offer significant improvements in protection performance, their practical implementation is often associated with increased costs and communication complexity, posing challenges for large-scale deployment.

More recently, research has shifted toward integrated and adaptive protection paradigms that align more closely with planning requirements. Ahmadi et al.~\cite{ahmadi2026comprehensive} proposed a unified protection framework in which AI-based fault detection systems are tightly integrated with energy management and optimization layers. In this approach, diagnostic outputs are quantified as fault-specific risk indicators that directly inform optimization constraints or risk-weighted cost functions, while operational set-points provide feedback to dynamically adjust detection thresholds. This bidirectional interaction enables a coordinated treatment of protection and operation, providing a structured mechanism to manage protection complexity and reduce unnecessary protection actions, thereby mitigating both reliability and economic impacts.

The application of artificial intelligence to microgrid protection has also been explored from a predictive maintenance perspective. Hadi et al.~\cite{hadi2025artificial} framed fault detection and classification as proactive processes based on real-time monitoring of voltage, current, frequency, and harmonic distortion, enabling early fault identification and improved power quality management. Complementarily, Ahmadi et al.~\cite{ahmadi2026comprehensive} classified AI-driven detection methods according to fault types—including high-impedance, short-circuit, open-circuit, and arc faults—highlighting their effectiveness in scenarios where conventional protection schemes struggle due to bidirectional power flow and variable fault currents. These developments demonstrate the potential of AI-enhanced protection to address the intrinsic uncertainties of microgrid operation.

Overall, it is clear that protection should not be treated as a purely operational add-on but as an integral component of microgrid planning. Incorporating protection constraints, and coordination requirements directly into planning models is essential to ensure that designed systems are not only cost-effective but also operationally feasible and reliable. Hence, microgrid planners should adopt AI-enhanced and adaptive protection frameworks at the planning stage to ensure reliable operation under high renewable penetration and uncertain fault conditions while limiting system complexity and coordination costs.

\subsubsection{Power Electronics Reliability}
\label{3.2.2}

Power electronics reliability represents a critical yet frequently overlooked factor in microgrid planning, as recently highlighted by Sandelic et al.~\cite{sandelic2022reliability}. Modern microgrids heavily rely on power electronic components to interface distributed energy resources, regulate voltage and frequency, and enable flexible operation under both grid-connected and islanded modes. Consequently, the reliability of these components directly influences the overall availability and operational robustness of microgrid systems. However, these components are prone to frequent failures, influenced by factors such as power loading, ambient temperature, humidity, and thermal cycling~\cite{peyghami2020overview,liu2022impact}. Therefore, neglecting the reliability impact of power electronics can lead to inaccurate microgrid planning, resulting in insufficient generation capacity, unpredictable outages, and higher lifecycle costs than anticipated during the design phase.

A key challenge arises from the fact that power electronics failure rates are not constant over their operational lifetime. As reported by Sandelic et al.~\cite{sandelic2022reliability}, failure behavior typically follows a bathtub-shaped profile, with relatively stable rates during the useful life phase and a rapidly increasing failure probability during the wear-out phase. Conventional microgrid planning models often assume constant failure rates, implicitly neglecting wear-out effects and stress-induced degradation. This simplification can result in inaccurate estimates of long-term reliability, replacement needs, and lifecycle costs, particularly for converter-dominated microgrids with high renewable penetration. However, to address this limitation, a more accurate approach could involve shifting from constant failure rate assumptions to the integration of physics-of-failure-based modeling within microgrid planning. This approach  investigates failure mechanisms and their root causes by establishing explicit relationships between operational mission profiles and component degradation mechanisms. In practice, this involves translating time-varying electrical and thermal stressors, such as junction temperature, thermal cycling amplitude, and switching frequency, into degradation rates using established reliability models, such as temperature acceleration relationships and formulations based on thermal cycling~\cite{peyghami2020overview}. These models enable the derivation of stress-dependent failure rates or mean time to failure estimates that reflect realistic operating conditions.

However, a major practical barrier to the implementation of physics-of-failure-based reliability models in microgrid planning is data availability. Mission-profile data are often proprietary to converter manufacturers or simply unavailable for specific models deployed in field installations. Furthermore, long-term field failure data for inverter-dominated microgrids remain scarce, as many such systems have operational histories shorter than the wear-out phase of their power electronic components (typically 10–15 years), limiting the calibration and validation of physics-based degradation models against real-world failure observations. The challenge is again compounded by the diversity of converter topologies, control strategies, and environmental conditions encountered across microgrid deployments, which makes it difficult to generalize failure data from one installation to another. While emerging data-driven and AI-based approaches may help augment limited datasets through synthetic stress scenarios, their validation against field data remains an open challenge. Consequently, caution is required when applying physics-of-failure models calibrated primarily on synthetic data to planning-scale decisions.

Despite these challenges, microgrid planners should integrate degradation-aware reliability models for power electronics into their planning frameworks, as this represents a paradigm shift from conventional static reliability assumptions to dynamic, physics-informed decision-making, which offers substantial practical benefits. By incorporating stress-dependent failure behavior, planners can then develop mission profile-based derating strategies that proactively limit converter loading or adjust capacity allocation to mitigate accelerated degradation under high-stress conditions. Furthermore, time-varying failure rates should be incorporated into reliability-constrained planning formulations to inform maintenance scheduling, component replacement decisions, and accurate expected energy not served. Similarly, degradation dynamics should also be explicitly reflected in lifecycle cost assessments to enable balanced trade-offs between upfront investment costs and long-term reliability performance. Importantly, accounting for power electronics reliability at the planning stage has implications beyond component longevity. As degradation alters the converter behavior over time, fault response characteristics and protection coordination can deviate from the initial assumptions, reinforcing the need to jointly consider reliability and protection during microgrid design. Consequently, planners who adopt the aforementioned methodologies will be better positioned to design systems that are both reliable and cost-effective, so as to maintain performance integrity over extended operational lifetimes.

\subsection{Economic Factors}
\label{3.3}

\subsubsection{Total Cost of Ownership (TCO)}\nomenclature[A]{TCO}{Total Cost of Ownership}
\label{3.3.1}

The total cost of ownership (TCO) provides a comprehensive lifecycle perspective that covers all costs from inception to decommissioning, thereby enabling planners to evaluate the true economic implications of microgrid design decisions. In the reviewed literature, TCO is consistently decomposed into capital expenditures (CAPEX) and operational expenditures (OPEX), although the level of component granularity and cost attribution varies substantially~\cite{thirunavukkarasu2023comprehensive, sandelic2022reliability, meera2022reliability}. CAPEX comprises equipment procurement, installation, grid interconnection (if grid-connected), and preliminary assessments, whereas OPEX includes fuel costs (for dispatchable sources), scheduled and unscheduled maintenance, control system updates, and personnel expenses. Therefore, when aiming for high reliability, it is vital to weigh these costs against the benefits of long-term reliability. 

Again, a comparative assessment of the economic indicators employed in the literature reveals notable methodological diversity. Levelized Cost of Energy (LCOE), defined as the ratio of the total discounted lifecycle costs to the total energy generated, remains the most widely used metric for cross-technology comparison~\cite{thirunavukkarasu2023comprehensive, meera2022reliability}. However, LCOE alone does not adequately reflect the value of reliability, as it is insensitive to supply interruptions, redundancy investments, or resilience improvement measures. To address this limitation, several studies have adopted more comprehensive cost indicators, including net present cost (NPC) and annual levelized cost (ALC), which aggregate annualized capital, operation, maintenance, replacement, and fuel costs over the project horizon~\cite{thirunavukkarasu2023comprehensive}. 

Notably, the interdependencies between the TCO and other planning parameters are particularly critical in reliability-oriented microgrid design. For example, the composition of the generation mix significantly influences both capital and operational cost structures, with renewable-dominated configurations exhibiting higher upfront costs but lower operational expenditures, while dispatchable generation reduces capital requirements while incurring ongoing fuel costs~\cite{gamarra2015computational}. Geographic location further influences resource availability, construction costs, environmental exposure, and fuel logistics, while microgrid scale introduces economies of scale for larger systems but might require additional network reinforcement and protection investments. Most importantly, high-reliability targets substantially influence TCO through redundancy requirements, protection system complexity, reserve margins, and backup generation provisions. In addition, ignoring replacement costs when project durations exceed component lifetimes, particularly for batteries with 8-15 year lifespans, can also lead to unrealistic planning outcomes~\cite{malika2025critical}. Such omissions might result in oversized storage installations, suboptimal coordination between distributed generation and energy storage devices, and systematic underestimation of TCO, which can jeopardize the financial sustainability of the project. These interdependencies highlight the fundamental challenge of balancing cost and reliability, which requires holistic planning frameworks that simultaneously optimize both objectives.

\subsubsection{Interruption Cost}
\label{3.3.2}

\begin{figure}[!t] 
    \centering 
    \includegraphics[width=1\linewidth]{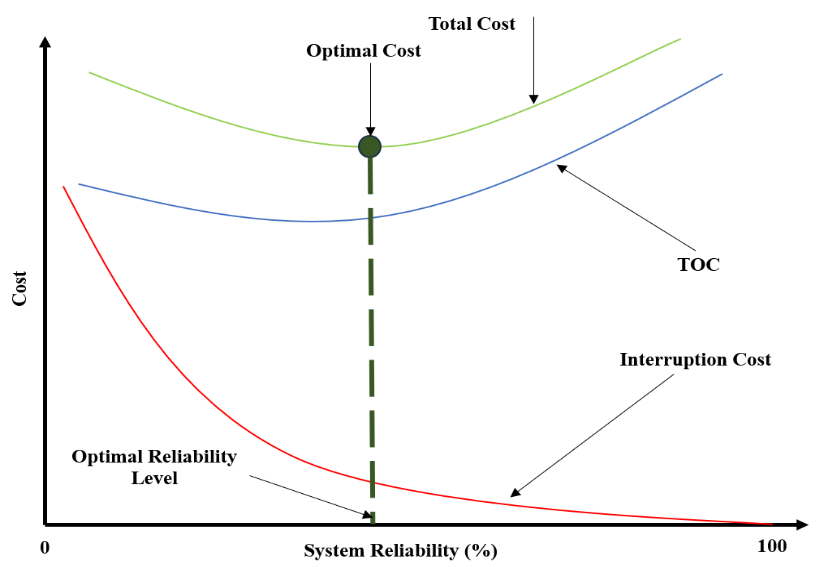} 
    \captionsetup{font=normalsize, justification=centering}
    \caption{Illustration of the tradeoff between reliability and cost.} 
    \label{Fig.4}  
\end{figure}

Interruption cost quantifies the economic impact of power supply disruptions on end-users and represents a critical input for reliability-oriented microgrid design. However, the magnitude of interruption costs varies substantially across customer types, creating heterogeneous reliability requirements within microgrid service territories. For commercial and industrial applications, even short-duration outages can result in substantial revenue losses~\cite{polleux2022overview}, operational disruptions, equipment damage, and contractual penalties, whereas residential customers generally experience comparatively lower direct economic impacts~\cite{morrissey2018cost,macmillan2023shedding}, with consequences more closely related to comfort, convenience, and quality of service. In particular, interruption cost exhibits an inverse relationship with TCO, defining the fundamental economic trade-off in reliability-driven microgrid design, as shown in Figure~\ref{Fig.4}. As planners invest more in reliability measures such as redundant generation units, advanced protection systems, and backup storage capacity, the TCO naturally increases. Consequently, these investments simultaneously minimize interruption costs by reducing the frequency and duration of outages. This inverse relationship establishes an economically optimal design point at which the marginal cost of reliability improvement is equal to the marginal benefit of reduced interruption costs. Design choices that fall below this optimum represent underinvestment, resulting in inadequate system reliability and excessive outage-related losses that outweigh capital savings. In contrast, design choices that exceed this optimum represent overinvestment, yielding diminishing economic returns despite continued reliability improvement.

This challenge of balancing the TCO against the interruption costs requires a standardized method to economically express the impact of reliability. However, while the reviewed literature consistently acknowledges interruption costs as a critical design parameter, it has not explicitly defined how these interruption costs should be quantified in practice. In this regard, the Value of Lost Load (VOLL)~\cite{gorman2022quest,ovaere2019detailed} \nomenclature[A]{VOLL}{Value of Lost Load} and Customer Damage Functions (CDF)~\cite{ravaghi2021probabilistic,rickerson2024value} \nomenclature[A]{CDF}{Customer Damage Function} emerge as essential frameworks, allowing planners to convert reliability into measurable financial impacts that can be directly compared against infrastructure investment costs. VOLL represents the economic value that customers assign to reliability, or equivalently, the cost incurred per unit of unserved energy during supply interruptions. It is typically expressed in $\text{\$/kWh}$ or $\text{\$/MWh}$, and reflects the willingness of customers to pay to avoid outages. In contrast, CDFs capture the nonlinear relationship between outage characteristics (duration, frequency, timing) and their associated economic impacts, recognizing that outage costs may increase disproportionately owing to longer outage durations, adverse timing, seasonal effects, and repeated occurrences. Hence, microgrid planners must incorporate explicit financial quantification of reliability, such as VOLL and CDFs, into reliability-oriented microgrid design to avoid both overinvestment in unnecessary reliability measures and underinvestment that leaves customers exposed to expensive outages.

\subsection{Social Factors}
\label{3.4}

As microgrid systems become increasingly decentralized and community-centric, social factors are no longer peripheral considerations but warrant explicit integration into microgrid planning~\cite{scheubel2017modeling, shandiz2020resilience}. The reviewed literature commonly treats customer satisfaction as the primary social dimension~\cite{gamarra2015computational}, typically represented through reliability, power quality, and environmental performance. However, social considerations extend beyond these technical metrics to include community engagement, trust building, and stakeholder participation in defining acceptable service levels. Importantly, understanding community needs and aspirations can provide insights that purely technical data cannot capture.

A key implication for reliability-oriented microgrid design is that “acceptable reliability'' is context dependent rather than uniform across communities. Different communities exhibit varying tolerance for power interruptions and differing expectations of reliable services, particularly with respect to critical facilities, vulnerable populations, and local economic activities. Consequently, reliability targets should not be treated as fixed engineering values, but should instead be informed through participatory processes that reflect local priorities, risk perceptions, and acceptable trade-offs among cost, reliability, and renewable integration.

In parallel with community-specific reliability preferences, socioeconomic heterogeneity further shapes reliability expectations by influencing willingness to pay and the distribution of reliability benefits. In this regard, Thirunavukkarasu et al.~\cite{thirunavukkarasu2023comprehensive} introduced social indicators such as the Human Development Index and job creation metrics to capture the broader socioeconomic implications of microgrid deployment. These indicators highlight how reliability improvements support essential services, local economic activity, and social development, thereby justifying reliability targets beyond purely techno-economic criteria. Together, these insights suggest that the incorporation of social and socioeconomic considerations can support context-aware, equitable, and broadly accepted reliability-oriented microgrid designs.

Finally, to systematically integrate technical, economic, and social considerations, the availability of standardized design guidelines is essential for robust microgrid planning. Building on the framework originally proposed in~\cite{8746836}, Sandelic et al.~\cite{sandelic2022reliability} presented a structured design and planning procedure that provides a comprehensive foundation for microgrid development. This framework begins with a clear definition of planning objectives, which are explicitly categorized into economic and reliability dimensions.

The design process then proceeds through successive stages, including system configuration, electrical design, and automation, with detailed guidance provided at each stage. The subsequent steps involve the formulation and comparison of alternative microgrid configurations using technical, economic, and environmental performance criteria, enabling an informed evaluation of the trade-offs between competing design objectives. The final design is selected based on an integrated assessment that balances these factors according to the predefined planning goals.

\begin{figure}[!t] 
    \centering 
    \includegraphics[width=0.95\linewidth]{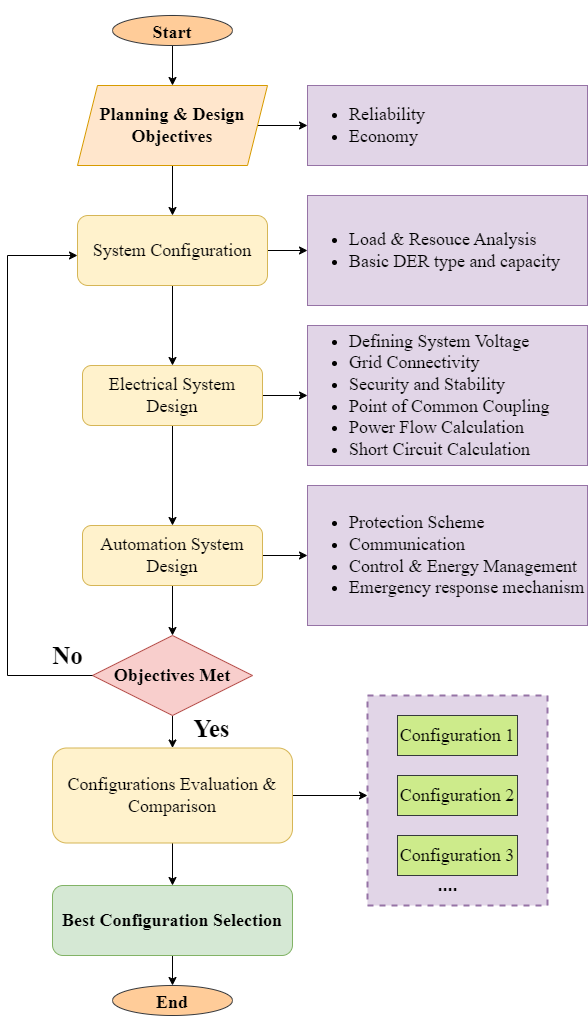} 
    \captionsetup{font=small}
    \captionsetup{justification=centering}
    \caption{Microgrid design and planning procedure (modified for reliability-oriented Design~\cite{sandelic2022reliability,8746836}).}
    \label{Fig.5}  
\end{figure}

This standardized procedure serves as a fundamental blueprint for microgrid planning, clarifying key objectives, decision points, and evaluation metrics. Figure~\ref{Fig.5} illustrates the corresponding microgrid planning and design flowchart, which was adapted to emphasize reliability considerations. While planners may tailor or extend specific steps to suit local conditions or emerging technologies, adherence to this structured guideline provides a consistent and transparent foundation for designing reliable, cost-effective, and socially responsive microgrid systems in the future.  

\section{Optimization Techniques in Reliable Microgrid Design}
\label{sec:4}
As discussed in Section 3, designing reliable microgrids is inherently complex and requires simultaneous consideration of technical, economic, and social factors. Again, planning objectives such as cost minimization and reliability maximization often exhibit conflicting relationships, requiring careful trade-offs. This is fundamentally an optimization problem~\cite{obara2018study,malika2025critical,gamarra2015computational}, where optimization techniques serve as the mathematical framework to mediate these trade-offs, addressing diverse planning decisions, including generation mix selection~\cite{dallavalle2023improved}, component sizing~\cite{wu2020stochastic}, optimal siting~\cite{cetinay2017optimal}, and scheduling~\cite{qing2024energy}.  However, the selection of an appropriate optimization technique depends critically on the mathematical structure of the problem, uncertainty characteristics, computational constraints, and the nature of the decision variables involved. Therefore, this section aims to provide a systematic taxonomy of the optimization models and solution techniques employed in reliability-aware microgrid planning.

Gamarra et al.~\cite{gamarra2015computational} formally defined optimization as the process of finding the best solution from a set of feasible alternatives with respect to specified criteria, subject to the defined constraints. This process involves an iterative interaction between an optimizer and a system model, wherein the optimizer adjusts the decision variables based on the model's evaluation of the objective functions and constraint satisfaction until the convergence criteria are met. In the context of reliable microgrid planning, decision variables encompass component capacities, discrete placement decisions, operational set points, protection device locations, and network topology configuration. The formulation of these optimization problems, comprising objective functions, constraints, and decision variables, constitutes the foundational step that determines both the quality of the solution and the computational tractability~\cite{malika2025critical,ehsan2018optimal}. 

\subsection{Optimization Models}
\label{4.1}

The formulation of optimization models represents the critical first step in microgrid planning, as the mathematical structure fundamentally determines the solution space, computational requirements, and, ultimately, the quality of planning decisions. A poorly selected model can result in suboptimal solutions that compromise microgrid performance, operational efficiency, and most critically, reliability. Therefore, understanding the characteristics, strengths, and limitations of different optimization models is essential for planners seeking to achieve specific reliability objectives.

The planning time horizon fundamentally shapes how reliability considerations are incorporated into optimization models~\cite{malika2025critical,gamarra2015computational}.This classification distinguishes between long-term strategic decisions and short-term tactical adjustments, each with distinct implications for the reliability of the system. Offline (static) optimization models address long-term strategic decisions such as system topology, generation mix selection, and capacity planning. These models utilize historical data and long-term forecasts to determine the optimal configurations that will serve the microgrid over its planning horizon, which typically spans 10 to 25 years. For example, when a planner must decide whether to install a 500 kW solar array with 200 kWh battery storage or a 300 kW solar array with 400 kWh storage, an offline model evaluates these alternatives based on projected load growth, renewable resource availability, and reliability targets. However, static formulations inherently assume that averaged or representative conditions adequately capture system behavior, potentially overlooking critical reliability events that occur during specific operational periods, such as coincident peak demand and low renewable generation. 

In contrast, online (dynamic) optimization models operate in real time or near real time, continuously adjusting to immediate changes in load demand, generation availability, and system conditions. From a reliability perspective, these models are essential for operational decisions, such as load-shedding sequences during contingencies, real-time dispatch adjustments when a generator goes offline, or activation of backup resources during islanding events. For instance, when an upstream grid fault forces a microgrid into islanded operation, an online optimization model must rapidly determine which loads to curtail and how to redispatch available generation to maintain frequency stability—decisions that directly impact reliability.

Again, multi-period (dynamic) formulations bridge the gap between static and online approaches by optimizing over multiple time steps while capturing inter-temporal dependencies~\cite{malika2025critical}. This is particularly important for reliability-oriented planning involving energy storage systems, where the state-of-charge at any given hour depends on all previous charging and discharging decisions. In a microgrid serving critical loads, a multi-period model ensures that storage reserves are maintained at sufficient levels before anticipated high-risk periods (such as forecasted storms or peak demand hours), rather than allowing storage to be depleted during low-risk periods and leaving the system vulnerable when reliability is most critical.

Notably, the treatment of uncertainty represents perhaps the most critical differentiator among optimization models for reliability-oriented planning, since reliability itself is fundamentally a probabilistic concept~\cite{malika2025critical,yan2024integrated}. How a model handles uncertainty in renewable generation, load demand, component failures, and grid availability directly determines whether the resulting design performs reliably under real-world conditions. Deterministic models use fixed parameter values without probabilistic elements, assuming perfect knowledge of all input parameters~\cite{abunima2022two}. For example, a deterministic model might assume that solar generation always equals 20\% of the rated capacity during the winter months and that the load demand follows a fixed daily profile. Although computationally efficient and straightforward to implement, this simplification often neglects uncertainties critical to reliability. Hence, a microgrid sized using deterministic assumptions may perform adequately under "average" conditions but fail during the very scenarios that matter most for reliability, such as a week of cloudy weather coinciding with an unexpected load increase.

On the other hand, stochastic optimization models explicitly incorporate uncertainty through probability distributions and optimize the expected value of objectives across multiple scenarios~\cite{malika2025critical}. Instead of assuming a fixed renewable output, for example, 20\% of the installed capacity, a stochastic formulation may evaluate a set of scenarios spanning a wider range of capacity factors, such as 5\% - 35\%, with each scenario weighted by its likelihood of occurrence. For reliability-oriented planning, this approach enables the model to size components such that the reliability indices remain below acceptable thresholds across the full spectrum of plausible operating conditions. For instance, a stochastic model may indicate that achieving a 99\% supply reliability target requires approximately 15\% additional battery capacity relative to the deterministic solution that fails to capture low-probability but high-impact scenarios where renewable shortfalls coincide with peak demand. However, the limitation of stochastic optimization lies in its computational burden, particularly when many uncertain parameters must be considered simultaneously, and a large number of scenarios are required to ensure adequate representation.

Conversely, robust optimization models adopt a worst-case perspective, ensuring that the solution remains feasible and satisfy predefined reliability requirements under the most adverse conditions within a defined uncertainty set~\cite{yang2021robust,malika2025critical}. Unlike stochastic optimization, which optimizes expected performance, robust optimization prioritizes guaranteed system performance against extreme but plausible combinations of uncertain parameters. In reliability-oriented microgrid planning, this approach ensures that the generation and storage capacities are sufficient to withstand simultaneous unfavorable conditions, such as prolonged periods of low renewable availability coinciding with peak demand. The principal advantage of robust optimization lies in its ability to provide explicit performance guarantees, which makes it particularly well-suited for safety-critical or mission-critical applications in which service interruptions are unacceptable. However, this conservatism often leads to higher capital requirements than stochastic approaches because designs must accommodate extreme scenarios. Consequently, the selection between stochastic and robust optimization reflects a fundamental trade-off between economic efficiency based on expected conditions and reliability assurance under worst-case scenarios.

Chance-constrained optimization models provide a pragmatic middle ground that is particularly well-suited for reliability-oriented planning~\cite{malika2025critical}. Rather than enforcing constraints deterministically under all conditions or optimizing performance expectations, chance-constrained formulations require that critical constraints be satisfied with a specified probability. For example, planners may impose requirements that voltage limits be respected at least 95\% of the time or that load curtailment remain below 5\% of demand with 99\% confidence. This formulation aligns closely with how reliability targets are defined in practice, as planners do not seek to eliminate all interruptions but instead aim to guarantee predefined reliability levels, such as achieving 99.99\% supply availability. It is also important to note that chance-constrained optimization is increasingly recognized as the natural framework for reliability-constrained planning~\cite{hemmati2020optimal,huo2022reliability}, although it requires careful specification of acceptable violation probabilities and may present computational challenges depending on the underlying uncertainty distributions.

The number of objectives considered further differentiates optimization models into single-objective and multi-objective formulations~\cite{prakash2016optimal}. Single-objective models focus on optimizing a single criterion, such as minimizing cost or maximizing reliability~\cite{fioriti2022multiple}. While computationally efficient, such formulations often overlook other relevant planning objectives, limiting their ability to capture the multidimensional nature of reliability-oriented microgrid design.

In contrast, multi-objective optimization models simultaneously address multiple, often competing objectives, allowing systematic exploration of trade-offs between cost, reliability, environmental performance, and other planning criteria~\cite{sachs2016multi}. These formulations provide a more comprehensive decision framework but typically incur higher computational complexity and often require advanced solution techniques. Moreover, multi-objective optimization is commonly formulated using either weighted-sum formulations or Pareto-based approaches. The weighted-sum method aggregates multiple objectives into a single function using predefined importance weights, whereas Pareto-based methods explicitly identify sets of non-dominated solutions, providing a richer representation of cost-reliability trade-offs at the expense of increased computational effort~\cite{thirunavukkarasu2023comprehensive}.

It is important to note that, despite the differences among various models, all optimization models share a common mathematical foundation. In this context, it is essential to distinguish the three fundamental components of any optimization model: objective functions, constraints, and decision variables. The objective function quantifies what optimization seeks to achieve, such as minimizing cost, maximizing reliability, reducing environmental impact, or improving power quality~\cite{malika2025critical,gamarra2015computational,ehsan2018optimal}. For reliability-oriented microgrid planning, objective functions are designed to balance service continuation against economic and operational constraints under uncertainty. Complementing the objective function, constraints define the boundaries of feasible solutions by enforcing physical laws, operational limits, and regulatory requirements~\cite{malika2025critical, ehsan2018optimal}. Without appropriate constraints, an optimization model may yield solutions that are mathematically optimal but physically infeasible or operationally unacceptable. Although constraints such as power balance, voltage limits, and thermal ratings are generic to microgrid planning, reliability-oriented formulations must incorporate explicit reliability constraints. These are commonly implemented by imposing threshold limits on reliability indices or constraining the probability of load curtailment and component unavailability through chance-based or scenario-dependent reliability requirements. Lastly, decision variables define the design and operational choices determined by the optimization model and directly shape the reliability of the system. In reliability-oriented microgrid planning, these typically include sizing variables for generation and energy storage capacities, placement variables that determine the location of critical components or protective devices, and operational variables that govern dispatch, storage charging and discharging, and load management over time~\cite{prakash2016optimal, malika2025critical}. Power flow variables further capture the resulting network state, such as bus voltages and line currents, ensuring that reliability requirements are met under feasible electrical operating conditions. Table 3 presents an integrated overview of representative decision variables, constraints, and objective functions commonly employed in microgrid planning studies. It serves as a unified reference framework that supports microgrid planning from multiple domain perspectives, including sizing, operation, reliability, protection, economics, and environmental considerations. However, a detailed explanation of each formulation is beyond the scope of this review, as the primary focus is on the optimization techniques employed for reliability-oriented microgrid design. Readers seeking in-depth discussions of individual constraints and their physical interpretations are referred to the cited literature~\cite{malika2025critical,gamarra2015computational,ehsan2018optimal}. The following subsections critically examine the optimization techniques adopted in the literature to address microgrid planning problem.

\begin{sidewaystable*}[p]
\centering
\caption{Integrated overview of representative formulations for microgrid planning across the literature~\cite{malika2025critical,gamarra2015computational,ehsan2018optimal}.}
\label{tab:optimization_overview}

\rotatebox{0}{%
\begin{minipage}{0.9\textheight}
\centering
\begin{tabular}{@{}>{\raggedright\arraybackslash}p{3cm}
                p{5.2cm}p{6cm}p{6cm}@{}}
\toprule
\textbf{Domain} & \textbf{Decision Variables} & \textbf{Constraints} & \textbf{Objective Functions} \\
\midrule

\makecell[l]{\textbf{SIZING}\\\textbf{\&}\\\textbf{PLACEMENT}}
& $\bullet$ $P_{\text{DG},i}$: DG active power capacity at bus $i$
\newline $\bullet$ $Q_{\text{DG},i}$: DG reactive power capacity at bus $i$
\newline $\bullet$ $E_{\text{ESS},i}$: ESS energy capacity (kWh) at bus $i$
\newline $\bullet$ $P_{\text{ESS},i}$: ESS power rating (kW)
\newline $\bullet$ $N_{\text{DG}}$, $N_{\text{ESS}}$: Number of units to install
\newline $\bullet$ $x_i \in \{0,1\}$: Binary placement decision
\newline $\bullet$ $y_i \in \{0,1\}$: ESS installation decision
& $\bullet$ $P_{\text{DG},\min} \leq P_{\text{DG},i} \leq P_{\text{DG},\max}$
\newline $\bullet$ $E_{\text{ESS},\min} \leq E_{\text{ESS},i} \leq E_{\text{ESS},\max}$
\newline $\bullet$ $\sum x_i \leq N_{\text{DG},\max}$
\newline $\bullet$ $\sum y_i \leq N_{\text{ESS},\max}$
\newline $\bullet$ Total investment $\leq \text{Budget}_{\max}$
\newline $\bullet$ Space/land availability at each bus
& $\bullet$ Min Total Investment Cost
\newline $\bullet$ Min Levelized Cost of Energy (LCOE)
\newline $\bullet$ Min Total System Cost (TSC)
\newline $\bullet$ Min Net Present Cost (NPC)
\newline $\bullet$ Max DG/RES penetration level
\\ \midrule

\makecell[l]{\textbf{TECHNICAL}\\\textbf{\&}\\\textbf{OPERATIONAL}}
& $\bullet$ $P_{ij}(t)$: Active power flow in line $i$--$j$
\newline $\bullet$ $Q_{ij}(t)$: Reactive power flow in line $i$--$j$
\newline $\bullet$ $V_i(t)$: Voltage magnitude at bus $i$
\newline $\bullet$ $I_{ij}(t)$: Magnitude of the current flowing through line $i$--$j$ at time $t$
\newline $\bullet$ $\delta_i(t)$: Voltage angle at bus $i$
\newline $\bullet$ $P_{\text{DG}_{dis},i}(t)$: DG dispatch at time $t$
\newline $\bullet$ $\text{SOC}_j(t)$: State of charge of ESS $j$
\newline $\bullet$ $P_{\text{ch},j}(t)$: Charging power of ESS $j$
\newline $\bullet$ $P_{\text{dis},j}(t)$: Discharging power of ESS $j$
\newline $\bullet$ $P_{\text{grid}}(t)$: Power exchange with main grid
& $\bullet$ $\sum P_{\text{gen}}(t) = \sum P_{\text{load}}(t) + P_{\text{loss}}(t)$ (power balance)
\newline $\bullet$ $\sum Q_{\text{gen}}(t) = \sum Q_{\text{load}}(t) + Q_{\text{loss}}(t)$ (reactive balance)
\newline $\bullet$ $V_{\min} \leq V_i(t) \leq V_{\max}$
\newline $\bullet$ $|I_{ij}(t)| \leq I_{ij}^{\max}$ (thermal limit)
\newline $\bullet$ $\text{SOC}_{\min} \leq \text{SOC}_j(t) \leq \text{SOC}_{\max}$
\newline $\bullet$ $\lvert \Delta P_{\text{DG},i}(t) \rvert \leq \text{RR}_{\max}$ (ramp rate)
\newline $\bullet$ $\text{PF}_i \geq \text{PF}_{\min}$ (power factor)
& $\bullet$ Min Power Loss: $\sum R_{ij} I_{ij}^2$
\newline $\bullet$ Min Voltage Deviation: $\sum (V_i - V_{\text{nom}})^2$
\newline $\bullet$ Min Line Loading Index
\newline $\bullet$ Max Voltage Stability Index (VSI)
\newline $\bullet$ Min Total Operating Cost
\\ \midrule

\makecell[l]{\textbf{RELIABILITY,}\\\textbf{ADEQUACY}\\\textbf{\&}\\\textbf{PROTECTION}}
& $\bullet$ $x_{\text{sw},k} \in \{0,1\}$: Switch placement at location $k$
\newline $\bullet$ $x_{\text{rec},k} \in \{0,1\}$: Recloser installation
\newline $\bullet$ $x_{\text{sect},k} \in \{0,1\}$: Sectionalizer placement
\newline $\bullet$ $x_{\text{tie},mn} \in \{0,1\}$: Tie-line between feeders $m$--$n$
\newline $\bullet$ $x_{\text{backup},i} \in \{0,1\}$: Backup DG installation
\newline $\bullet$ $z_{\text{restore},i}(t) \in \{0,1\}$: Restoration sequence
\newline $\bullet$ $N_{\text{sections}}$: Number of feeder sections
& $\bullet$ SAIFI $\leq \text{SAIFI}_{\text{target}}$
\newline $\bullet$ SAIDI $\leq \text{SAIDI}_{\text{target}}$
\newline $\bullet$ CAIDI $\leq \text{CAIDI}_{\text{target}}$
\newline $\bullet$ LOLE $\leq 0.1$ days/year
\newline $\bullet$ EENS $\leq \varepsilon\%$ of annual demand
\newline $\bullet$ ASAI $\geq 0.9999$
\newline $\bullet$ Reserve margin $\geq \text{RM}_{\min}$
\newline $\bullet$ $\mathbb{P}(V_{\min} \leq V_i \leq V_{\max}) \geq 1-\varepsilon$ [chance]
\newline $\bullet$ $\mathbb{P}(P_{\text{gen}} \geq P_{\text{load}}) \geq 1-\varepsilon$ [chance]
& $\bullet$ Min SAIFI: $\sum(\lambda_i N_i)/N_{\text{total}}$
\newline $\bullet$ Min SAIDI: $\sum(U_i N_i)/N_{\text{total}}$
\newline $\bullet$ Min EENS: $\sum \mathbb{P}(s) E_{\text{curtailed}}(s)$
\newline $\bullet$ Min LOLE: $\sum \mathbb{P}(C < L)\,T$
\newline $\bullet$ Max EIR: $1-\text{EENS}/E_{\text{demand}}$
\newline $\bullet$ Min Cost of Unserved Energy
\newline $\bullet$ Min Customer Interruption Cost

\\ \midrule
\bottomrule
\end{tabular}
\end{minipage}}
\end{sidewaystable*}

\begin{sidewaystable*}[p]
\ContinuedFloat
\caption[]{Integrated overview of representative formulations for microgrid planning across the literature~\cite{malika2025critical,gamarra2015computational,ehsan2018optimal}}

\rotatebox{0}{%
\begin{minipage}{0.9\textheight}
\centering
\begin{tabular}{@{}>{\raggedright\arraybackslash}p{3cm}
                p{5.2cm}p{6cm}p{6cm}@{}}
\toprule
\textbf{Domain} & \textbf{Decision Variables} & \textbf{Constraints} & \textbf{Objective Functions} \\
\midrule

\textbf{ECONOMIC}
& $\bullet$ Sizing variables ($P_{\text{DG}}$, $E_{\text{ESS}}$, $x_i$, $y_j$)
\newline $\bullet$ Operational variables ($P_{\text{DG}_{dis},i}(t)$, $\text{SOC}_j(t)$,$P_{\text{ch},j}(t)$,$P_{\text{dis},j}(t)$)
\newline $\bullet$ Protection variables (same as Reliability \& Protection domain)
\newline \newline \textit{Note: Economic domain uses variables from other domains; costs are calculated outputs.}
& $\bullet$ Total investment $\leq \text{Budget}_{\max}$
\newline $\bullet$ Annual O\&M cost $\leq C_{\text{oper},\max}$
\newline $\bullet$ Payback period $\leq \text{PB}_{\text{target}}$
\newline $\bullet$ IRR $\geq \text{IRR}_{\min}$
\newline $\bullet$ NPV $\geq 0$
& $\bullet$ Min Total Cost: $\sum C_{\text{inv}} + \sum C_{\text{oper}} + \sum C_{\text{maint}}$
\newline $\bullet$ Min Net Present Cost (NPC)
\newline $\bullet$ Min Levelized Cost of Energy (LCOE)
\newline $\bullet$ Max Net Present Value (NPV)
\newline $\bullet$ Max Benefit-Cost Ratio (BCR)
\newline $\bullet$ Min Total Ownership Cost (TOC)
\\ \midrule

\textbf{ENVIRONMENTAL}
& $\bullet$ Sizing variables ($P_{\text{DG}}$, type selection)
\newline $\bullet$ Operational variables ($P_{\text{DG}}(t)$ dispatch)
\newline \newline \textit{Note: Emissions are calculated from dispatch decisions; RES fraction from sizing decisions.}
& $\bullet$ $\sum E_{\text{CO}_2} \leq \text{Emission}_{\text{limit}}$
\newline $\bullet$ $\sum E_{\text{NO}_x} \leq \text{NO}_x\ \text{limit}$
\newline $\bullet$ $E_{\text{renewable}}/E_{\text{total}} \geq \text{RES}_{\text{target}}$
\newline $\bullet$ Carbon intensity $\leq \text{CI}_{\max}$ (kg/kWh)
\newline $\bullet$ Compliance with environmental regulations
& $\bullet$ Min Total Emissions: $\sum P_{\text{DG}}(t)\,E_{\text{rate}}$
\newline $\bullet$ Min Carbon Footprint
\newline $\bullet$ Min Greenhouse Gas (GHG) Cost
\newline $\bullet$ Max Renewable Energy Fraction
\newline $\bullet$ Min Life Cycle Emissions
\\

\bottomrule
\end{tabular}
\end{minipage}}
\end{sidewaystable*}

\subsection{Analytical Techniques}
\label{4.2}

Analytical techniques offer a mathematically rigorous approach to microgrid design and optimization, relying on explicit mathematical formulations and power system theory to capture system behavior in a computationally efficient way. Based on closed-form formulations and theoretical analyses, these methods derive solutions from power flow equations, sensitivity analyses, and probabilistic representations of historical data~\cite{thirunavukkarasu2023comprehensive,centeno2023computationally}. Their principal advantages include computational efficiency~\cite{ehsan2018optimal}, accuracy~\cite{prakash2016optimal}, and strong physical interpretability~\cite{thirunavukkarasu2023comprehensive}, which make them particularly suitable for small-scale microgrids, where deterministic and tractable solutions are required. From a reliability perspective, analytical techniques allow planners to explicitly relate design variables and operating conditions to reliability indices without relying on extensive simulation or heuristic search procedures.

Ehsan et al.~\cite{ehsan2018optimal} proposed a comprehensive taxonomy of analytical formulations for microgrid planning, identifying six principal approaches: the exact loss formula, loss sensitivity factor, branch current loss formula, equivalent current injection, branch power flow loss formula, and phase feeder current injection. Each formulation provides a distinct mathematical representation of power flows and losses, enabling planners to target specific optimization objectives. Common decision variables across these formulations include the size, location, and power factor of distributed generation units, optimized subject to load characteristics and network constraints. Importantly, the branch power flow loss formula incorporates stochastic representations of renewable generation, making it particularly relevant for reliability-oriented planning in renewable-dominated microgrids where uncertainty plays a central role in system performance. 

Again, complementing formulation-based perspectives, Prakash et al.~\cite{prakash2016optimal} categorized analytical techniques according to their functional role in microgrid planning. Eigenvalue-based analysis focuses on system stability and dynamic behavior under varying loads and generation conditions, supporting reliable operation during disturbances and transients. Index methods quantify the relative changes in system performance and reliability metrics, enabling vulnerability assessment, benchmarking, and prioritization of network improvements. Sensitivity-based methods evaluate how variations in design variables influence target outputs such as losses or reliability indices, thereby reducing the search space and supporting efficient sizing and siting decisions. Point estimation methods address the uncertainty in renewable generation and load by approximating probabilistic inputs through representative points, allowing reliability indices to be estimated with a significantly lower computational burden than the full Monte Carlo simulation.

However, despite their advantages, analytical techniques exhibit inherent limitations as system complexity increases. Malika et al.~\cite{malika2025critical} positioned analytical approaches within a broader taxonomy that includes numerical, evolutionary, and advanced computational methods, clearly delineating their applicability. Their analysis showed that analytical techniques perform well for systems with limited dimensionality and moderate uncertainty but face challenges when extended to high-resolution temporal modeling, nonlinear component behavior, or complex interactions among multiple sources of uncertainty. These challenges are particularly evident in reliability-oriented planning, where discrete events such as component failures or network reconfigurations—introduce combinatorial complexities that are difficult to solve through closed-form equations alone.

To address these limitations, recent studies have explored extensions that retain the strengths of analytical methods while enhancing their flexibility. Hadi et al.~\cite{hadi2025artificial} argue that the growing complexity of modern microgrids necessitates more adaptive computational approaches, while emphasizing the continued importance of analytical formulations in structuring optimization problems and enforcing physical feasibility. Ahmadi et al.~\cite{ahmadi2026comprehensive} further demonstrate that hybrid frameworks, which combine analytical preprocessing techniques such as signal decomposition or feature extraction with data-driven models, can improve robustness to noise and uncertainty while preserving interpretability. From a reliability standpoint, these developments suggest that analytical techniques remain valuable either as standalone tools for tractable planning problems or as structural components within hybrid optimization frameworks.

Overall, it is evident that analytical techniques remain fundamental to reliability-oriented microgrid planning due to their computational efficiency, transparency, and strong physical grounding. However, their limited scalability restricts their effectiveness in highly complex or uncertainty-intensive planning scenarios. Consequently, analytical techniques are increasingly embedded within integrated optimization frameworks, where they provide structural rigor while complementary methods address system complexity and uncertainty.

\subsection{Mathematical Optimization Techniques}
\label{4.3}

Mathematical optimization techniques constitute one of the most established and theoretically grounded approaches for microgrid design and reliability-oriented planning. These methods identify optimal system configurations and operational strategies by minimizing or maximizing a well-defined objective function subject to a set of physical, operational, and reliability constraints~\cite{gamarra2015computational,prakash2016optimal,liu2019development}. Their structured formulation provides designers with explicit control over system modeling, enabling a transparent representation of trade-offs among capital investment, reliability indices, and technical performance~\cite{de2022balancing,abdmouleh2017review}. While the optimization models introduced in section~\ref{4.1} define the overarching planning objectives (e.g., handling uncertainty through Robust or Stochastic frameworks), the choice of mathematical technique determines the computational feasibility and physical accuracy of the resulting microgrid design. In reliability-oriented planning, these techniques effectively serve as the computational engines that solve the underlying mathematical models while enforcing system constraints. Consequently, mathematical optimization constitutes a cornerstone of reliability-oriented planning, particularly when rigorous guarantees of constraint satisfaction and solution optimality are required.

Gamarra and Guerrero~\cite{gamarra2015computational} classify mathematical optimization techniques applied to microgrid planning into several categories, including linear programming (LP) \nomenclature[A]{LP}{Linear Programming}, mixed integer linear programming (MILP)\nomenclature[A]{MILP}{Mixed-integer Linear Programming}, mixed integer nonlinear programming (MINLP)\nomenclature[A]{MINLP}{Mixed-integer Non Linear Programming}, dynamic programming (DP)\nomenclature[A]{DP}{Dynamic Programming}, and sequential quadratic programming (SQP) \nomenclature[A]{SQP}{Sequential Quadratic Programming}. Among these, LP, MILP, and MINLP remain the most widely used formulations in reliability-oriented microgrid planning as they provide an effective trade-off between modeling flexibility and computational tractability. Linear programming has been extensively employed in early microgrid studies where both objective functions and constraints can be expressed in linear form~\cite{zhang2022systematic}. Its primary advantage lies in computational efficiency and guaranteed convergence, making it suitable for optimizing resource allocation, such as distributed generation and storage sizing, under simplified reliability constraints~\cite{sandelic2022reliability}. However, the inherent linearity of LP limits its ability to capture discrete decisions and nonlinear behaviors that are essential in realistic microgrid systems. To address these limitations, MILP has emerged as the dominant mathematical optimization technique in microgrid planning~\cite{wu2021milp,li2017microgrid,marocco2021milp}.

By incorporating integer variables, MILP enables explicit modeling of discrete design and operational decisions, including unit commitment, switch status, protection device placement, and energy storage charging or discharging states. This capability is particularly important for reliability-oriented planning, in which system redundancy, contingency handling, and operational flexibility must be accurately represented. In reliability-oriented planning, MILP is typically used to integrate reliability metrics as either hard constraints or economic penalties within the objective function. Consequently, MILP has been widely adopted for generation mix selection, capacity sizing, and reliability-constrained planning problems~\cite{sandelic2022reliability,ghanbarzadeh2025addressing}. Moreover, MILP is frequently used to solve scenario-based formulations, where each scenario is represented as a set of linear constraints. However, as the temporal horizon expands to a full year (8,760 h) to capture seasonal reliability variations, MILP faces a "curse of dimensionality" owing to the exponential growth of binary variables. In addition, MILP formulations typically rely on linear approximations of system behavior and therefore struggle to represent nonlinear phenomena, such as power flow physics, battery degradation, and inverter dynamics.

MINLP is an advanced optimization technique that extends MILP by allowing both nonlinear relationships and discrete decision variables within the same optimization framework~\cite{elsido2017two}. This flexibility enables MINLP to capture complex interactions among microgrid components, such as nonlinear battery aging, voltage-dependent power flows, and renewable generation characteristics. In reliability-oriented planning, MINLP is particularly useful for jointly optimizing cost, reliability indices, and operational performance while accounting for component degradation and nonlinear constraints. Although numerous articles reviewed in this rapid review used MINLP for microgrid design, it is also important to acknowledge the scalability limitations of MINLP, as the increased complexity of large-scale problems can render MINLP computationally challenging. Moreover, owing to the nonconvex nature of the power flow equations~\cite{de2022balancing}, MINLP often struggles with local optima and extreme computational costs. To bridge this gap, recent studies have focused on convexification techniques, such as second-order cone programming (SOCP)~\cite{gil2021mixed} and semidefinite programming (SDP)~\cite{molzahn2019survey}. These methods relax or transform non-convex power flow equations into convex forms, allowing planners to retain the accuracy of nonlinear behavior while achieving the computational robustness and global optimality associated with linear solvers.

Additionally, DP and SQP represent alternative mathematical optimization approaches that are well-suited to specific reliability-oriented applications. Dynamic programming decomposes time-coupled problems into smaller subproblems, making it effective for scheduling, unit commitment, and energy management tasks under reliability constraints~\cite{de2022balancing,das2022approximate}. Sequential quadratic programming, on the other hand, iteratively solves nonlinear constrained optimization problems and is particularly useful when modeling continuous nonlinear behavior of microgrid components~\cite{prakash2016optimal}. However, while both methods are effective for solving time-coupled or nonlinear subproblems, their scalability becomes limited when reliability must be evaluated across multiple time periods, contingencies, or uncertainty realizations.  This scalability limitation has motivated the adoption of decomposition-based mathematical optimization techniques, particularly Benders Decomposition~\cite{lin2019optimal} and Column-and-Constraint Generation (C\&CG)~\cite{wu2020integrating}, in reliability-oriented microgrid planning. 

These decomposition methods offer two complementary advantages for reliability-oriented planning. First, they improve computational tractability by decomposing large-scale problems into coordinated master and subproblems, each addressing a manageable subset of decision variables and constraints. Second, and more critically for reliability planning, they provide a systematic framework to explicitly handle uncertainty by iteratively identifying adverse or worst-case operating scenarios. Within this hierarchical structure, the master problem optimizes strategic planning decisions such as generation capacity, storage sizing, and network topology, while subordinate subproblems evaluate operational feasibility and quantify reliability violations under stochastic or adversarial uncertainty realizations. This iterative exchange ensures that reliability requirements are satisfied not only in average representative scenarios but also under worst-case or high-impact conditions that critically determine system reliability. Despite their conceptual suitability for reliability-oriented microgrid planning, decomposition-based optimization techniques have received limited attention in the reviewed literature. Apart from a brief application of Benders Decomposition mentioned in a single study~\cite{zhang2022systematic}, none of the reviewed articles discuss these techniques, and their role in enforcing reliability under stochastic or robust uncertainty. This limited coverage highlights an important methodological gap and suggests that decomposition-based mathematical optimization techniques remain an underexplored yet promising direction for scalable and uncertainty-aware reliability-oriented microgrid planning.

\subsection{Heuristic \& Metaheuristic Techniques}
\label{4.4}

As microgrid planning problems grow in scale, dimensionality, and computational complexity, traditional analytical and mathematical optimization techniques may struggle to deliver feasible or timely solutions. This challenge is particularly evident in large-scale mixed-integer formulations and non-deterministic polynomial-time hard problems, where exact solution methods can exhibit slow convergence or fail to identify feasible solutions within practical computational limits. Such limitations are especially critical in reliability-oriented microgrid planning, where timely decision-making and the ability to explore complex design spaces are essential. Consequently, heuristic and metaheuristic optimization techniques have gained increasing attention as practical alternatives for solving computationally intensive microgrid planning problems~\cite{gamarra2015computational,akter2024review, thirunavukkarasu2023comprehensive} owing to their quicker~\cite{alrashidi2020metaheuristic} and approximate solutions that are often good enough for practical use.

Heuristic methods employ problem-specific rules or search strategies to rapidly identify high-quality solutions with relatively low computational overhead. Common examples in microgrid planning include greedy algorithms~\cite{dimovski2023holistic}, tabu search~\cite{vincent2023sustainable}, and simulated annealing (SA)~\cite{zhang2018optimization}. These techniques are attractive because of their simplicity and speed. However, their reliance on local search mechanisms renders them susceptible to entrapment in locally optimal solutions, particularly in highly nonconvex and multi-modal reliability planning problems. In contrast, metaheuristic techniques extend beyond problem-specific heuristics by providing general-purpose, stochastic search frameworks that balance exploration and exploitation of the solution space~\cite{zhang2018optimization,malika2025critical}. Inspired by natural, biological, or social processes, metaheuristics iteratively improve candidate solutions and can handle complex, nonlinear, and multi-objective microgrid planning formulations involving discrete design choices and reliability constraints~\cite{ehsan2018optimal}. Their flexibility and robustness have led to their widespread adoption in reliability-oriented microgrid designs, where exhaustive searches are infeasible.

\nomenclature[A]{SA}{Simulated Annealing}

Among the broad family of metaheuristic techniques, genetic algorithms (GA) and particle swarm optimization (PSO) have emerged as the most widely applied methods in microgrid planning studies~\cite{leonori2020optimization,hossain2019energy,raghavan2020optimization}. A recent review by Akter et al.~\cite{akter2024review} reports that PSO and GA account for approximately 25\% and 10\% of metaheuristic applications, respectively, highlighting their prominence in the literature. Genetic algorithms are based on principles of natural selection and evolutionary genetics, employing operators such as selection, crossover, and mutation to evolve populations of candidate solutions~\cite{prakash2016optimal}. Their ability to explore large and discontinuous solution spaces makes them well-suited for multi-objective microgrid planning problems, such as simultaneous optimization of cost, reliability, and renewable penetration. In reliability-oriented contexts, GA has been extensively used for optimal sizing and placement of distributed energy resources and energy storage systems~\cite{alhamali2017determination,yaghoubi2021optimal}. Nevertheless, performance of GA is sensitive to parameter tuning and population size, and the repeated evaluation of complex reliability metrics can result in a substantial computational burden.

Particle swarm optimization, on the other hand, draws inspiration from collective social behavior observed in biological systems~\cite{ehsan2018optimal}. One of the advantages of PSO over GA is its ability to improve solution quality with fewer iterations~\cite{prakash2016optimal}. It operates by initializing a swarm of particles in the solution space, each of which represents a potential solution. These particles then iteratively adjust their positions based on their own experiences and those of their neighbors, aiming to find the global optimum. PSO has been successfully applied to the siting and sizing of distributed generators, compensators, and storage systems in reliability-constrained microgrid planning~\cite{karunarathne2020optimal,sharafi2014multi}. However, PSO is prone to premature convergence, particularly in high-dimensional reliability problems, and its performance can degrade as system size and constraint complexity increase.

\nomenclature[A]{GA}{Genetic Algorithm}
\nomenclature[A]{PSO}{Particle Swarm Optimization}

Beyond GA and PSO, numerous bio-inspired metaheuristics have been explored in microgrid planning, including grey wolf optimizer~\cite{nimma2018grey}, harmony search~\cite{sun2024harmony}, differential evolution~\cite{ramli2018optimal}, artificial bee colony~\cite{saeed2021two}, ant colony optimization~\cite{guven2022design}, cuckoo search~\cite{singh2022energy}, crow search~\cite{thirunavukkarasu2023comprehensive}, and whale optimization algorithms~\cite{liu2022improved}. Although these methods demonstrate promising performance in specific applications, their comparative advantages are often problem-dependent, and detailed algorithmic analysis is beyond the scope of this review. For a comprehensive understanding of the various metaheuristic approaches, interested readers are encouraged to refer to the cited literature~\cite{akter2024review,thirunavukkarasu2023comprehensive,prakash2016optimal,gamarra2015computational}. 

A comparative synthesis of the heuristic and metaheuristic literature reveals a pattern of both convergence and persistent disagreement. There is broad agreement that GA and PSO outperform classical heuristics for multi-objective, combinatorial microgrid planning problems. In particular, PSO is frequently reported to achieve comparable solution quality with fewer iterations, which makes it attractive for large scale applications where computational efficiency is critical. However, this agreement does not extend to issues of robustness. GA-based studies frequently report a high sensitivity to parameter tuning, particularly population size and the configuration of crossover and mutation operators. In contrast, research on PSO consistently identifies premature convergence as a key limitation, especially in high dimensional reliability-oriented problems where the search space is complex and multimodal. These contrasting observations suggest that neither method offers a universally reliable performance, and that their effectiveness remains problem dependent. 

Bio inspired optimization methods, on the other hand, have been proposed as alternatives that may overcome some of these limitations. While these methods often report superior performance in specific case studies, the evidence base remains fragmented. Most contributions rely on isolated applications without employing standardized benchmark problems or consistent evaluation criteria. As a result, the reported performance improvements are difficult to generalize. The absence of systematic, controlled comparisons across reliability focused benchmarks represents a significant methodological gap in the literature. Without such a framework, the selection of appropriate metaheuristic techniques lacks a principled foundation, and the rapid proliferation of new variants risks prioritizing novelty over cumulative scientific insight. A further point of divergence concerns computational cost. In reliability oriented planning, the need to repeatedly evaluate candidate solutions through Monte Carlo simulation introduces substantial computational burden. Under these conditions, even relatively efficient metaheuristic algorithms may become impractical for large scale problems.

Overall, heuristic and metaheuristic optimization techniques play a vital role in reliability-oriented microgrid planning by enabling efficient exploration of complex, high-dimensional design spaces where exact methods become impractical. Although these methods do not guarantee global optimality, their flexibility, scalability, and adaptability render them indispensable for real-world reliability studies. However, the literature increasingly suggests that the greatest potential of these methods lies not in their standalone application but in hybrid optimization frameworks, as discussed in Section 4.6.

\subsection{AI-based Techniques}
\label{4.5}

Artificial intelligence (AI)-based techniques represent the most recent and rapidly evolving methodological direction in reliability-oriented microgrid planning~\cite{mohammadi2022review,zahraoui2024ai}. Unlike analytical, mathematical, or metaheuristic approaches that rely on predefined formulations and explicit models, AI techniques learn complex nonlinear patterns and decision policies directly from multi-modal data. This data-driven capability has enabled AI methods to address challenges that are increasingly difficult to capture through conventional optimization, including high-dimensional uncertainty, nonlinear component behavior, and complex interactions between demand, renewable generation, and network constraints~\cite{hadi2025artificial,ahmadi2026comprehensive,malika2025critical}.

Recent reviews provide a structured taxonomy of AI techniques applicable to microgrid planning. Hadi et al.~\cite{hadi2025artificial} distinguish three principal categories: machine learning, deep learning, and reinforcement learning. Machine learning techniques include supervised learning for forecasting and classification tasks, unsupervised learning for clustering and anomaly detection to identify early stage component failures, and semi-supervised learning for applications with limited labeled data. Deep learning extends these capabilities through multi-layer neural architectures such as convolutional neural networks for spatial feature extraction, recurrent and long short-term memory networks for sequential and temporal data, and transformer-based models for attention-driven learning. Reinforcement learning, in contrast, formulates planning as a sequential decision-making process. Unlike static optimization, RL agents learn optimal dispatch and protection policies through continuous interaction with a simulated environment, making them highly suited for real-time reliability management during islanding events~\cite{malika2025critical}.

Ahmadi et al.~\cite{ahmadi2026comprehensive} complemented this taxonomy by categorizing AI techniques according to their application domains in microgrid systems, including fault detection, reliability assessment, and optimization. Within the context of reliability-oriented microgrid planning, AI techniques are primarily applied at the strategic and operational planning levels of microgrid planning. Hadi et al.~\cite{hadi2025artificial} emphasize that most existing AI applications target tertiary-level decisions, including generation sizing, storage allocation, energy management strategies, and long-term scheduling under uncertainty. In these applications, AI models are used to forecast load demand, renewable generation availability, and failure likelihoods, thereby improving the quality of planning inputs and reducing reliance on overly conservative assumptions. Malika et al.~\cite{malika2025critical} further demonstrated that machine learning models can be used to generate plausible future scenarios for weather conditions, demand evolution and grid outages. These AI-generated scenarios are subsequently embedded into optimization frameworks, enabling planners to evaluate reliability performance across a broader and more realistic uncertainty space than is possible with historical averages or simplified probabilistic assumptions. This capability is particularly valuable for reliability-oriented planning, in which rare but high-impact events critically influence design decisions.

Another critical advancement in reliability planning is the use of Generative Adversarial Networks (GANs)~\cite{tang2021securing,ichinomiya2025generative} and Variational Autoencoders (VAEs)~\cite{kaur2021variational,jangilwar2025cso} to address the scarcity of historical failure data. These models generate high-fidelity, synthetic "stress-test" scenarios such as rare multi-component failures or extreme weather transients, that go beyond simple probabilistic distributions. Embedding these AI-generated scenarios into optimization frameworks allows planners to evaluate system reliability against "black swan" events that are often missed by traditional planning tools.

Despite their potential, AI-based techniques face limitations when applied to reliability-oriented microgrid planning. A central challenge identified by Ahmadi et al.~\cite{ahmadi2026comprehensive} is the simulation-to-reality gap, wherein AI models trained on simulated or historical datasets exhibit degraded performance when deployed in real-world systems. Contributing factors include insufficient representation of extreme operating conditions in training data, violations of model assumptions due to physical nonlinearities, computational latency, and limited interpretability of learned decision policies. These issues are particularly problematic in reliability planning, where incorrect predictions or opaque decisions can lead to under-designed systems or unacceptable risk exposure. To mitigate this, recent research has pivoted toward Physics-Informed Machine Learning (PIML)~\cite{xu2015evaluation,wu2024physics}. By embedding physical laws (such as Kirchhoff’s laws or swing equations) directly into the neural network's loss function, PIML ensures that the AI's predictions do not violate the fundamental principles of power system physics. This "physics-aware" approach provides the structural guardrails necessary to trust AI decisions in safety-critical microgrid operations. Complementing this, Hadi et al.~\cite{hadi2025artificial} further proposed a digital twin framework as a potential mitigation strategy, enabling extensive simulation-based testing and validation of AI models before deployment. However, the development and maintenance of accurate digital twins introduces additional modeling complexity and resource requirements, and discrepancies between the digital twin and the physical system remain a concern.

Recognizing these challenges, recent literature increasingly positions AI not as a standalone replacement for traditional planning methods, but as a complementary component within hybrid optimization frameworks. Ahmadi et al.~\cite{ahmadi2026comprehensive} identified several hybrid paradigms, including data-driven hybrid models that combine preprocessing techniques with AI predictors, optimization-based hybrids that integrate AI with mathematical programming or metaheuristics, and ensemble AI frameworks that combine multiple learning models to improve robustness and generalization. In reliability-oriented microgrid planning, these hybrid approaches allow AI to enhance scenario generation, uncertainty representation, and decision support, while established optimization techniques enforce physical feasibility, reliability constraints, and economic optimality. This integration helps address key limitations of purely data-driven models, including lack of guarantees on constraint satisfaction and limited transparency.

\subsection{Hybrid Techniques}
\label{4.6}

As mentioned briefly in earlier sections, hybrid optimization techniques have emerged as a pragmatic response to the limitations of standalone analytical, mathematical, heuristic, and AI-based methods in reliability-oriented microgrid planning. These approaches integrate two or more optimization paradigms within a unified framework to simultaneously exploit complementary strengths, such as global exploration, constraint handling, computational efficiency, and modeling accuracy~\cite{thirunavukkarasu2023comprehensive,de2022balancing,kumar2018optimizing}. By doing so, hybrid techniques aim to balance solution quality, convergence speed, and scalability when addressing complex planning problems involving discrete design decisions, nonlinear system behavior, and reliability constraints under uncertainty~\cite{zhang2022systematic}.

As discussed earlier, single-method optimization approaches often entail inherent trade-offs, particularly with respect to modeling flexibility, solution optimality, and computational tractability. Mathematical programming techniques provide optimality guarantees and rigorous constraint enforcement but may suffer from scalability issues in large or uncertainty-rich planning problems. Metaheuristic methods offer flexibility and global search capability but lack feasibility guarantees and may converge prematurely. AI-based approaches excel at learning complex patterns from data but face challenges related to interpretability, robustness, and deployment. Hence, Prakash et al.~\cite{prakash2016optimal} emphasized that hybrid optimization frameworks can address these shortcomings by combining methods in a targeted manner, enabling faster convergence and higher solution quality than any single technique.

A common class of hybrid techniques integrates metaheuristics with mathematical optimization. For example, coupling particle swarm optimization with mixed-integer linear programming leverages the global search capability of PSO while retaining the structured constraint handling and discrete decision modeling of MILP. Kim et al.~\cite{kim2020milp} demonstrate that using MILP to generate high-quality initial solutions for PSO allows the swarm to explore promising regions of the solution space, accelerating convergence and improving reliability-oriented design outcomes. Similarly, Almadhor et al.~\cite{almadhor2021hybrid} combine PSO with the bat algorithm to optimize distributed energy resource placement and sizing, exploiting PSO’s rapid convergence while using bat-algorithm dynamics to regulate exploration and reduce power losses.

Beyond PSO-MILP combinations, several studies report hybrid metaheuristic structures such as GA-PSO, GA-ANN, and simulated annealing combined with neural networks~\cite{malika2025critical}. These combinations typically use one algorithm to perform coarse global exploration and another to refine local solutions, improving robustness against local optima while maintaining acceptable computational effort. Gamarra et al.~\cite{gamarra2015computational} similarly advocate sequential hybridization, where one heuristic identifies an initial feasible design and a second heuristic performs targeted refinement to enhance reliability performance.

Again, recent reviews increasingly highlight the integration of AI techniques within hybrid optimization frameworks as a particularly promising direction for reliability-oriented microgrid planning. Malika et al.~\cite{malika2025critical} demonstrated that machine learning models can be used to generate realistic future scenarios for demand, weather, and outages, which are subsequently embedded within metaheuristic or mathematical optimization routines. This separation of scenario generation from optimization reduces computational burden while enabling planners to evaluate reliability performance across a richer uncertainty space. Ahmadi et al.~\cite{ahmadi2026comprehensive} proposed layered hybrid architectures in which AI-driven methods generate initial solutions, contingency suggestions, or uncertainty representations, while mathematical optimization models such as OPF or MILP act as feasibility filters to enforce physical and reliability constraints. This structure preserves interpretability and protection margins while recovering much of the performance advantage offered by data-driven intelligence. Importantly, within planning contexts, such architectures allow AI to enhance decision support without replacing the optimization backbone that ensures reliability compliance. Furthermore, the emergence of Digital Twin-driven hybridization allows for continuous validation of these hybrid models, providing the "virtual commissioning" necessary for reliability-oriented microgrids.

Despite their advantages, hybrid optimization techniques introduce additional complexity in model formulation, algorithm coordination, and parameter tuning. Singh et al.~\cite{singh2017review} note that many hybrid approaches have been applied primarily to systems with static loads or limited technology diversity, while their application to planning problems involving multiple distributed generation types and dynamic load profiles remains limited due to increased computational burden. Scalability also becomes a concern as hybrid frameworks grow in dimensionality and algorithmic interdependence. Parallel computing and high-performance computing resources have been suggested as potential enablers for addressing these challenges~\cite{gamarra2015computational}, though their adoption in planning-oriented studies remains limited.

Overall, hybrid optimization techniques are consistently positioned as an effective strategy for reliability-oriented microgrid planning, across the reviewed literature. Rather than replacing existing methods, hybrid frameworks combine mathematical rigor, heuristic flexibility, and data-driven insight to better address uncertainty, nonlinearity, and discrete decision spaces. Notably, this review identifies a critical gap in the existing literature. Despite the widespread adoption of hybrid optimization methods, no systematic framework exists to guide their selection for reliability-oriented microgrid planning. As a result, planners lack standardized guidance on choosing appropriate hybrid configurations based on microgrid scale, data availability, computational resources, and reliability targets. Addressing this architecture selection gap is essential for advancing hybrid optimization from largely ad hoc academic implementations to a standardized industrial planning tool.

\subsection{Software tools}
\label{4.7}

Beyond algorithmic optimization techniques, software platforms also play a decisive role in implementing microgrid planning methodologies by embedding modeling assumptions, optimization logic, and evaluation metrics into deployable decision-support environments~\cite{lenhart2021microgrid,thirunavukkarasu2023comprehensive,gamarra2015computational}. These tools differ not only in their underlying solvers but also in how explicitly they represent network topology, outages, uncertainty, and reliability trade-offs. Consequently, the choice of software tool can materially influence planning outcomes, particularly in reliability-oriented studies where outage behavior, load prioritization, and network feasibility are central concerns.

A broad ecosystem of tools has been applied across the reviewed literature, including HOMER~\cite{bouendeu2023systematic, al2023feasibility, ahmad2018techno}, GridLab-D~\cite{wang2022microgrid, ferrari2019assessment}, RETScreen~\cite{zaro2023design, rafique2018developing}, PVSYST~\cite{bouchekara2021decomposition, huda2024techno}, DER-CAM~\cite{jung2017optimal}, OpenDSS~\cite{jain2023simulation}, and the Sandia Microgrid Design Toolkit~\cite{eddy2017sandia} Each platform addresses specific stages of the planning workflow. For example, PVSyst and HelioScope are commonly used for detailed photovoltaic yield estimation~\cite{govindarajan2024analysing}, RETScreen for high-level financial screening~\cite{paradongan2024techno}, and GridLAB-D or OpenDSS for distribution-level power flow and voltage analysis~\cite{appiah2024survey}. However, no single tool fully integrates techno-economic optimization, electrical feasibility, and reliability assessment within a unified framework. As a result, many studies rely on multi-software workflows in which capacity optimization, reliability evaluation, and electrical validation are performed sequentially rather than jointly.

Among all available platforms, HOMER has emerged as the dominant industry and academic standard, particularly for early stage microgrid feasibility analyses~\cite{thirunavukkarasu2023comprehensive,gamarra2015computational,de2022balancing}. More than a thousand technical articles consistently cited HOMER as the most widely adopted and continuously used software~\cite{lambert2006micropower}. This is due to its accessibility, extensive technology libraries, and ability to rapidly enumerate and compare candidate system configurations based on net present cost. Moreover, HOMER’s chronological simulation over representative years enables basic assessment of supply adequacy under varying resource conditions~\cite{ringkjob2018review,babu2023sensitivity,uwineza2021feasibility}, making it well suited for cost-driven design problems. However, it is important to note that HOMER’s internal modeling structure is inherently cost-centric, which limits its applicability in reliability-oriented microgrid planning. Moreover, reliability considerations are not embedded within the optimization objective, but are evaluated indirectly through unmet load statistics after optimization. Additionally, HOMER does not incorporate standard reliability indices, outage costs, contingency modeling, or value-of-lost-load formulations. Furthermore, its single-node abstraction precludes the explicit modeling of networked microgrids, voltage constraints, and feeder-level power flows. As reliability requirements become more stringent, these simplifications necessitate complementary analyses using external tools, fragmenting the planning workflow and obscuring the interactions between cost, reliability, and network feasibility.

To address these limitations, a number of newer planning platforms have been proposed that place greater emphasis on representing reliability-related requirements alongside least-cost design, including XENDEE~\cite{pecenak2019efficient}. XENDEE is built around a mixed-integer linear programming (MILP) formulation and can be configured to consider multiple objectives and constraints (e.g., cost, emissions, redundancy, and reliability-related requirements) within a single optimization problem~\cite{mathiesen2021techno}. In contrast to legacy tools where reliability is evaluated only after a cost-optimal configuration is selected, outage scenarios and critical-load requirements can be represented within the optimization formulation, enabling reliability requirements to directly influence investment and operational decisions~\cite{Xendee2025_online}. In addition, the platform supports network-represented microgrid models and can be coupled with electrical analyses so that operational feasibility considerations (e.g., voltage limits, thermal constraints, and losses) can be incorporated during design rather than assessed solely as an ex-post validation step~\cite{Xendee2025_online}. It also supports explicit modeling of islanded operation duration, component redundancy, battery degradation, and load growth, making it particularly suited to mission-critical and community-scale microgrids where service continuity is paramount.

Despite these capabilities, XENDEE remains under-represented in the academic literature, primarily because it was introduced only recently (approximately 2018) and has therefore not yet been extensively examined in academic microgrid planning studies. Only a limited number of studies have applied this platform~\cite{mathiesen2021techno,stadler2019planning,pecenak2020impact}, and existing comparisons with HOMER are cursory and primarily cost-focused~\cite{musfiq2024comparison}. Notably, no reviewed study conducted a systematic, feature-based comparison of microgrid planning software from a reliability-oriented perspective. This gap is significant, as nearly all existing software comparisons conclude by default that HOMER is the most capable tool without evaluating emerging platforms against reliability-specific criteria such as outage modeling, network feasibility, and resilience optimization. Table~\ref{Table.4} provides a structured comparison of HOMER and XENDEE across features that are directly relevant to reliability-oriented microgrid planning.

\begin{table*}[!t]
\centering
\caption{Comparison of HOMER Pro and XENDEE for reliability-oriented microgrid planning.}
\label{Table.4}
\renewcommand{\arraystretch}{1.2}
\begin{tabular}{p{4.2cm} p{5.4cm} p{6.0cm}}
\toprule
\textbf{Feature} & \textbf{HOMER Pro} & \textbf{XENDEE} \\
\midrule
Optimization algorithm 
& Blackbox algorithm 
& Mixed-Integer Linear Programming (MILP) \\

Objective structure 
& Single-objective (net present cost) 
& Multi-objective (cost, reliability, emissions, resilience) \\

Reliability metrics in objective 
& Not supported 
& Explicitly supported (outage cost, critical load constraints) \\

Outage modeling 
& Limited, post-analysis only 
& Integrated outage windows and islanding constraints \\

Value of Lost Load (VOLL) 
& Not supported 
& Explicitly supported \\

Network topology modeling 
& Single-node, aggregated representation 
& Multi-node, GIS-based representation \\

Power flow analysis 
& Not supported 
& Supported via OpenDSS integration \\

Voltage and thermal constraints 
& Not modeled 
& Explicitly enforced \\

Battery degradation modeling 
& Simplified lifetime approximation 
& Explicit multi-year degradation modeling \\

Load growth modeling 
& Limited or static 
& Explicit multi-period load growth modeling \\

Redundancy and N--1 logic 
& Not supported 
& Explicitly supported \\

Typical application domain 
& Feasibility screening and cost minimization 
& Reliability- and resilience-driven microgrid planning \\

Academic adoption 
& Extensive and well-established 
& Limited but growing \\
\bottomrule
\end{tabular}
\end{table*}

Finally, by reviewing the full spectrum of optimization techniques employed in reliability-oriented microgrid planning, a qualitative understanding of their theoretical foundations, strengths, and limitations was established. However, qualitative synthesis alone does not reveal how these techniques are adopted in practice, nor does it capture emerging methodological trends across the literature. To address this gap, Table~\ref{tab:technique_distribution} presents a quantitative meta-analysis that systematically examines the frequency, temporal evolution, and co-occurrence of analytical, mathematical optimization, heuristic, AI-based, and hybrid techniques reported in the studies. The percentage distribution reported in Table~\ref{tab:technique_distribution} was derived through a systematic data extraction and categorization protocol designed to reflect the relative prevalence of optimization techniques reported in the reviewed articles rather than subjective weighting. Primary study counts were identified through a structured examination of technique-specific tables, reference lists, and methodological discussions within each review, with cross-validation applied to identify overlapping citations and minimize both omissions and duplications. Each primary study was assigned to a single, mutually exclusive category based on its dominant algorithmic foundation, ensuring internal consistency and that the aggregated shares summed to 100\%. For studies employing hybrid formulations, classification was guided by the dominant or novel methodological contribution as explicitly emphasized by the source article, thereby avoiding double-counting. Consequently, the reported distribution represents a normalized frequency-based synthesis of methodological adoption in microgrid planning, capturing prevailing research trends while accounting for overlaps across reviews.

The distribution reported in Table~\ref{tab:technique_distribution} reveals clear methodological preferences in reliability-oriented microgrid planning. Heuristic and Metaheuristic algorithms dominate the literature at 38.4\%, reflecting their flexibility in addressing nonconvex, multi-objective, and combinatorial planning problems. It was also found that  genetic Algorithm (26.9\% of category) and particle swarm optimization (24.6\%) collectively represent more than half of metaheuristic applications. This predominance reflects the natural suitability of population-based algorithms to handle discrete placement decisions and continuous sizing variables in microgrid planning problems. This widespread adoption of heuristics and meta heuristics can be attributed to their ability to incorporate complex reliability objectives without requiring strict mathematical formulations, making them well suited for large-scale and highly nonlinear planning problems. In contrast, mathematical programming techniques constitute the second largest share at 27.2\%, underscoring their continued importance in applications where exact optimality, constraint satisfaction, and model interpretability are essential. Their prevalence in optimal power flow, operational scheduling, and generation expansion planning highlights their suitability for structured problems with explicit reliability constraints. However, their comparatively lower share reflects known limitations related to scalability and nonconvexity, particularly when detailed reliability modeling, storage degradation, or multi-period uncertainty is introduced.

AI-based and machine learning techniques account for 17.1\% of the reviewed studies, indicating a rapidly growing but still maturing research direction. Their adoption is primarily motivated by the need to address uncertainty, forecasting, and fault-related tasks that are difficult to model deterministically. The lower proportion relative to metaheuristics reflects ongoing challenges associated with interpretability, data requirements, and robustness in planning-oriented applications. Again, hybrid optimization approaches represent 9.2\% of the literature, suggesting increasing recognition that no single method is sufficient for microgrid planning. These approaches aim to combine global search capability, computational efficiency, and constraint enforcement. Their limited adoption to date is attributed to higher implementation complexity and the absence of standardized hybrid design frameworks.

Analytical and classical methods comprise the smallest share at 5.7\%, consistent with their declining applicability to modern microgrids characterized by high renewable penetration and uncertainty. While still valuable for sensitivity analysis and preliminary assessment, their limited expressiveness constrains their role in comprehensive reliability-oriented planning. Overall, the distribution reflects a clear methodological shift from analytically tractable and exact formulations toward more flexible and hybridized approaches. This transition mirrors the growing complexity of reliability requirements in modern microgrids and highlights the need for integrated frameworks that balance optimality, scalability, and uncertainty representation.

\begin{table*}[!t]
\centering
\caption{Quantitative distribution of optimization techniques in microgrid planning.}
\label{tab:technique_distribution}
\begin{tabular}{l c p{6.8cm}}
\hline
\textbf{Category} & \textbf{Percentage} & \textbf{Primary application domains} \\
\hline
Mathematical  & 27.2\% & DG sizing and siting, optimal power flow, unit commitment, expansion planning, operational scheduling \\
Heuristic \& Metaheuristic & 38.4\% & DG sizing and siting, capacity optimization, multi-objective planning \\
AI / machine learning & 17.1\% & Load forecasting, fault detection, state estimation, predictive planning \\
Hybrid optimization & 11.6\% & Complex multi-objective problems, convergence enhancement, uncertainty handling \\
Analytical / classical & 5.7\% & Sensitivity analysis, power flow studies, stability assessment \\
\hline
\end{tabular}
\end{table*}

Finally, the preceding analysis of optimization techniques suggests that no single optimization formulation is universally optimal; rather, the suitability of an optimization formulation for reliability-oriented microgrid planning is fundamentally determined by how reliability is represented within the planning problem. For long-term planning problems that involve discrete investment decisions, multi-period operational scheduling, and explicit reliability constraints, mixed-integer linear programming (MILP) constitutes the most appropriate baseline formulation. Its capacity to encode design-feasibility logic, chronological energy balance, unit commitment, and contingency constraints within a single tractable framework, while providing provable optimality guarantees, makes it particularly well-suited to this class of planning problems.

When uncertainty is a dominant planning driver, as is characteristically the case in renewable-rich microgrids, stochastic programming and chance-constrained formulations become essential. These approaches ensure solution feasibility across the full distribution of plausible operating conditions, rather than under average or representative scenarios alone. Where planning objectives must additionally hedge against worst-case realizations or ambiguity in the uncertainty set, robust and distributionally robust optimization formulations are more appropriate, as they explicitly prioritize resilience assurance over expected-value performance.

In network-constrained microgrids where electrical fidelity is central to the reliability assessment, convex relaxations of the AC power flow equations, most notably second-order cone programming (SOCP) formulations, or suitably structured mixed-integer nonlinear programming (MINLP) models, are required to accurately represent power flows, voltage limits, and converter behavior. Although these formulations impose a higher computational burden relative to their linear counterparts, this cost is justified when network feasibility is a binding reliability requirement.

For large-scale or highly nonlinear design problems where exact mathematical programming methods face scalability limitations, hybrid frameworks that integrate mathematical programming with metaheuristic search algorithms or AI-based components offer a pragmatic improvement in computational tractability. Such methods are particularly valuable for scenario generation, surrogate modeling, and solution acceleration. It should be noted, however, that the reliability value of any hybrid approach remains bounded by the fidelity of the underlying planning formulation; algorithmic sophistication cannot compensate for inadequate reliability representation.

In summary, the reviewed literature consistently demonstrates that formulation quality, encompassing the adequate representation of chronological operation, uncertainty, contingency coverage, and network feasibility, exerts a substantially greater influence on planning outcome credibility than the choice of solution algorithm itself. This finding underscores the primacy of rigorous problem formulation in reliability-oriented microgrid planning research and practice.

\section{Evaluating Reliability: Methods and Metrics}
\label{sec:5}

While optimization techniques are central to the design of reliable microgrids, the credibility of any design ultimately depends on how reliability is evaluated. Therefore, reliability assessment serves not only as a validation step for optimized solutions but also as a diagnostic tool that informs topology refinement, component reinforcement, and protection coordination. Consequently, reliability evaluation is an integral component of reliability-oriented microgrid planning. Traditionally, reliability in power systems has been defined as the ability to supply electricity continuously within acceptable voltage and frequency limits with minimal service interruptions~\cite{escalera2018survey,alotaibi2016efficient}. This definition has historically focused on the impacts of customer interruption and the performance of distribution system operators during fault events. However, modern microgrids extend the reliability concept beyond continuity of supply to include adaptive operation under uncertainty, high penetration of renewables, islanding capability, and self-healing through fault isolation and reconfiguration~\cite{lopez2020reliability}. These expanded requirements require a more sophisticated and context-aware reliability evaluation. This section reviews the current methods and metrics for reliability assessment and their alignment with evolving microgrid systems.

\subsection{Reliability Evaluation Methods}
\label{5.1}
Reliability evaluation methods are designed to measure the ability of the system to adequately, and securely meet the required supply demands with minimal disruption~\cite{escalera2018survey,nazir2021system}. These methods aim to quantify the likelihood of system failure or inadequate performance considering various factors such as component outages~\cite{al2019reliability}, load variations~\cite{xiang2016power}, and system configurations~\cite{ramezanzadeh2021reliability}. Existing studies broadly classify these methods into analytical and simulation-based approaches, which are detailed in the following sections.

\subsubsection{Analytical Methods}
\label{5.1.1}

Analytical methods employ mathematical formulations and structured probabilistic models to represent the behavior of a system in failure states~\cite{meera2022reliability,allan2013reliability}. By decomposing complex systems into simplified states or logical structures, these methods provide closed or semi-closed-form estimates of reliability indices with relatively low computational burden~\cite{myhre2021reliability}. This efficiency makes them attractive for preliminary planning studies, sensitivity analyses, and systems with limited structural complexity. Common analytical techniques applied to microgrid reliability evaluation include Markov models, Bayesian networks, reliability block diagrams (RBD), fault tree analysis (FTA), minimal cut set analysis, and failure mode and effects analysis (FMEA)~\cite{meera2022reliability,lopez2020reliability,momoh2017electric}.

Markov models provide a mathematically rigorous, state-based approach to evaluate microgrid reliability by representing systems as a set of discrete operational conditions ranging from normal operation through partial degradation to complete outage~\cite{hong2017optimal}. Transitions between these states are modeled as stochastic processes governed by probabilistic failure and repair rates, often assuming exponential distributions to satisfy the memoryless property of Markovian systems. Building on this framework, Lopez et al.~\cite{lopez2020reliability} demonstrated that each state represents a specific operational condition, and the model quantifies the transition likelihoods using a stochastic transition probability matrix. For example, a microgrid comprising a conventional backup generator and battery can be modeled using states corresponding to the generator availability and battery charge status. However, the primary challenge in applying Markov models to modern microgrids lies in their computational scalability. As the number of components and their respective operational modes increases, the state space grows exponentially, potentially reaching millions of states for even medium-sized microgrids.

Bayesian Network offers an alternative in this context by explicitly modeling conditional dependencies among components through probabilistic graphical structures~\cite{meera2022reliability}. Unlike Markov models, Bayesian networks capture causal relationships between components, such as the propagation of inverter failure to load curtailment in renewable-dominated microgrids. This capability enables a more efficient representation of complex interdependencies without exhaustive state enumeration. However, the construction of accurate conditional probability tables requires extensive empirical data, which is often unavailable for emerging microgrid technologies. Similarly, Reliability Block Diagrams (RBDs) provide a path-centric approach by representing component connectivity and redundancy through logical success and failure configurations~\cite{vcepin2011assessment}. Although RBDs offer superior computational efficiency and intuitive visualization, making them ideal for identifying topological vulnerabilities during preliminary planning, their application to dynamic microgrids is often restricted by their traditionally static nature. Specifically, RBDs frequently rely on the assumption of independent component failures and cannot inherently represent complex temporal dependencies, such as sequential state-of-charge variations in energy storage or the coordinated stochastic repair processes required for reliable operations.

\nomenclature[A]{RBD}{Reliability Block Diagram}
\nomenclature[A]{FTA}{Fault Tree Analysis}
\nomenclature[A]{FMEA}{Failure Mode and Effects Analysis}

Notably, Fault tree analysis (FTA) extends this logic by identifying combinations of component failures that lead to system-level failure using Boolean logic structures~\cite{baig2013reliability}. FTA maps multiple failure scenarios to a single top event, making it particularly effective for identifying common-cause failures and vulnerabilities in complex topologies. To transition from qualitative failure pathway identification to quantitative vulnerability assessment, minimal cut set analysis is frequently employed in conjunction with FTA~\cite{momoh2017electric}. This analysis represents the minimal combination of component failures whose simultaneous occurrence results in the top event. By identifying these sets, planners can isolate the most critical components from a reliability perspective, particularly those characterized by single-order cut sets wherein failure of a single component precipitates system failure. This identification enables the prioritization of reliability enhancement measures, such as strategic redundancy deployment or infrastructure hardening, at points of maximum system vulnerability.

Conversely, when a granular assessment of each component's failure modes is required, FMEA serves as a robust, inductive framework for reliability evaluation. Unlike the top-down deductive logic of FTA, FMEA systematically investigates each component in isolation to identify potential failure modes and their subsequent system-level impacts~\cite{momoh2017electric}.. In the context of reliability-oriented planning, this method allows for the evaluation of system-level impacts through a sequential examination of component states. Once a component failure is identified, its specific consequences for the end-user are quantified before transitioning to the next operational contingency.

Overall, it is important to note that analytical techniques have been invaluable for years, offering an average idea of how often things might go wrong~\cite{lopez2020reliability}. However, they exhibit fundamental limitations when applied to modern microgrids. Most notably, they struggle to represent low-probability, high-impact events such as extreme weather phenomena or extended renewable generation deficits. In addition, these methods face challenges in capturing nonlinear interactions among distributed generation, energy storage, and demand response, as well as time-varying operational strategies for modern microgrids~\cite{lopez2020reliability,ghanbarzadeh2025addressing}. As the scale and complexity of microgrids increase, these limitations reduce the fidelity of analytical reliability assessments and motivate the adoption of simulation-based or hybrid evaluation methods.

\subsubsection{Simulation-Based Methods}
\label{5.1.2}

Simulation-based reliability evaluation methods have become the dominant approach for assessing microgrid reliability, particularly in systems dominated by high renewable penetration, complex network topologies, and stochastic operating conditions~\cite{meera2022reliability}. Unlike analytical methods that rely on simplified mathematical abstractions, simulation-based methods explicitly model the stochastic behavior of system components, operational strategies, and temporal dependencies, allowing a realistic representation of microgrid operation under uncertainty~\cite{escalera2018survey,lopez2020reliability}.

MCS constitutes the most widely adopted framework within this paradigm. MCS evaluates system reliability by repeatedly simulating random realizations of component failures, repairs, load variations, and generation uncertainty over a defined time horizon, subsequently aggregating system performance metrics across numerous iterations~\cite{meera2022reliability}. This approach accommodates arbitrary probability distributions and complex interactions among distributed energy resources, storage systems, loads, and network configurations, making it particularly well-suited for microgrid reliability studies~\cite{ahmadi2026comprehensive}. For example, in a microgrid comprising photovoltaic generation, battery storage, and a backup diesel generator, a single MCS iteration might simulate a solar inverter failure occurring during peak demand with concurrent low battery state-of-charge, enabling quantification of the resulting load curtailment duration and magnitude across thousands of such scenarios.

\nomenclature[A]{MCS}{Monte Carlo Simulation}
\nomenclature[A]{SMCS}{Sequential Monte Carlo Simulation}
\nomenclature[A]{NSMCS}{Non-Sequential Monte Carlo Simulation}

MCS methods are broadly classified into sequential and non-sequential approaches~\cite{lopez2020reliability}. Sequential MCS preserves the chronological evolution of system states by explicitly modeling time-dependent component failure and repair processes~\cite{escalera2018survey}. This temporal fidelity proves essential for microgrids, where the reliability of the system depends on the "history" of operations such as whether a battery was charged yesterday to survive a power outage today. For example, while a static model might only look at average solar output, a sequential simulation can capture high-impact scenarios such as an inverter failure occurring exactly during a week-long \textit{Dunkelflaute} when battery reserves are already depleted. However, this chronological accuracy incurs substantial computational expense, as extended simulation horizons and large numbers of iterations are required to achieve statistical convergence.

In contrast, Non-Sequential Monte Carlo Simulation (NSMCS) evaluates the reliability of microgrids by treating each simulated state as an independent event~\cite{lopez2020reliability}, entirely decoupled from the chronological order of operations. By sampling component availability based on fixed probability distributions rather than a timeline, NSMCS significantly reduces the computational burden required for a study. This makes it an efficient tool for preliminary reliability assessments or for systems where the components do not have strong "memory" dependencies. However, this simplification limits its ability to accurately represent storage behavior, repair sequences, and prolonged outage events, which are increasingly important in reliability-oriented microgrid planning.

To overcome the limitations mentioned above, recent studies have explored strategies to enhance the computational efficiency of MCS while preserving modeling accuracy. Meera et al.~\cite{meera2022reliability} reported that machine learning techniques such as support vector machines and random forests can approximate system responses or classify failure states, thus accelerating the convergence of MCS. Such approaches prove particularly valuable when MCS is embedded within optimization frameworks, where repeated reliability evaluations would otherwise impose prohibitive computational demands.

In this context, Lopez-Prado et al.~\cite{lopez2020reliability} presented a comprehensive methodology for reliability evaluation using MCS, which could be considered as a state-of-the-art technique for microgrids. By incorporating a comprehensive set of variables, this methodology provides a holistic framework that captures the multifaceted nature of reliability in microgrids. However, while a flowchart was provided, the authors did not explain the methodology explicitly. A comprehensive breakdown of this methodology into several operational phases is provided below.

\begin{itemize}
  \item \textbf{Initialization Phase}
    \begin{enumerate}
      \item \textbf{Parameter Definition:} In this initial step, the essential parameters have been defined for each component, and denoted as \(N\), \(T\), \(\lambda_c\), \(\mu\), and TTR, which represent the number of iterations, time horizon, failure rate of components, repair rate, and time to repair for each component, respectively.
    \end{enumerate}

\item \textbf{Simulation Phase}
    \begin{enumerate}
      \item \textbf{Variable Initialization:} Variables \(n\) and \(t\) are set to 1 and 0, serving as counters for the number of iterations and time, respectively.
      \item \textbf{Stochastic Simulation:} A set of random numbers \(X\) is generated to represent the stochastic nature of component failures, each uniformly distributed over the interval \([0,1]\).
      \item \textbf{Time-To-Failure (TTF) Calculation:} The random numbers \(X\) are transformed into Time to Failure (TTF) for each component using the following equation:
\begin{equation}
    \text{TTF}= -\frac{1}{\lambda_c} \ln(X), \tag{1}
\end{equation}
thereby quantifying the expected lifespan of each component.
    \end{enumerate}

\item \textbf{Failure Identification}
    \begin{enumerate}
      \item \textbf{Component Selection:} In this step, the algorithm identifies the component $C$ with the minimum TTF and designates it as the component most likely to fail first. The time to failure of the component \(C\) is denoted as \(\text{TTF}_{c}\).
      \item \textbf{Spatial Localization:} In this part of the algorithm, it identifies where the component \(C\) is physically located within the overall layout of the system.
      \end{enumerate}

\item \textbf{Repair and Impact Assessment}
    \begin{enumerate}
        \item \textbf{Repair Time Estimation:} A new random number is generated and converted into Time To Repair (TTR) for component $C$ using the following equation:
        \begin{equation}
        \text{TTR} = -\frac{1}{\mu} \ln(X). \tag{2}
        \end{equation}
        \item \textbf{Load Point Analysis:} In this step, the following equation 
        \begin{equation}
        L_k(t) = w_h (h)\times w_m(m) \times P_{lk}, \tag{3}
        \end{equation}
        is used to calculate the electrical load at each load point in the system and identifies those adversely affected by the failure of component \(C\). Note that $L_k(t)$ represents the amount of electrical demand at a specific time $t$. The factors $w_h(h)$, $w_m(m)$, and $P_{lk}$ are used to adjust this demand based on the hour of the day, the month, and the peak load at that point $k$, respectively.
        \item \textbf{Distributed Generation (DG) Scenarios:} The algorithm generates various scenarios for Distributed Generation (DG) in this step, to assess how they can reduce the impact of any failure of component \(C\).
    \end{enumerate}
    \item \textbf{Iterative Update and Analysis}
    \begin{enumerate}
        \item \textbf{Failure Duration and Impact:} In this step, the algorithm calculates the duration of the failure and the unsupplied load at each affected load point.
        \item \textbf{TTF Recalculation:} A new TTF is calculated for component \(C\) for the subsequent iteration.
        \item \textbf{Failure Metrics:} The algorithm calculates the average failure rate and duration across all load points.
    \end{enumerate}
    \item \textbf{Convergence and Output}
    \begin{enumerate}
        \item \textbf{Iteration Check:} In this step, the algorithm checks whether $n<N$. If true, it returns to the stochastic simulation phase; otherwise, it proceeds to the final step.
        \item \textbf{Reliability Indices Calculation:} Reliability indices are calculated for the system, based on the aggregated data from all iterations.
        \item \textbf{Termination:} The algorithm concludes and the calculated reliability indices are compiled.
    \end{enumerate}
\end{itemize}

\begin{figure}[!t] 
    \centering 
    \includegraphics[width=1\linewidth]{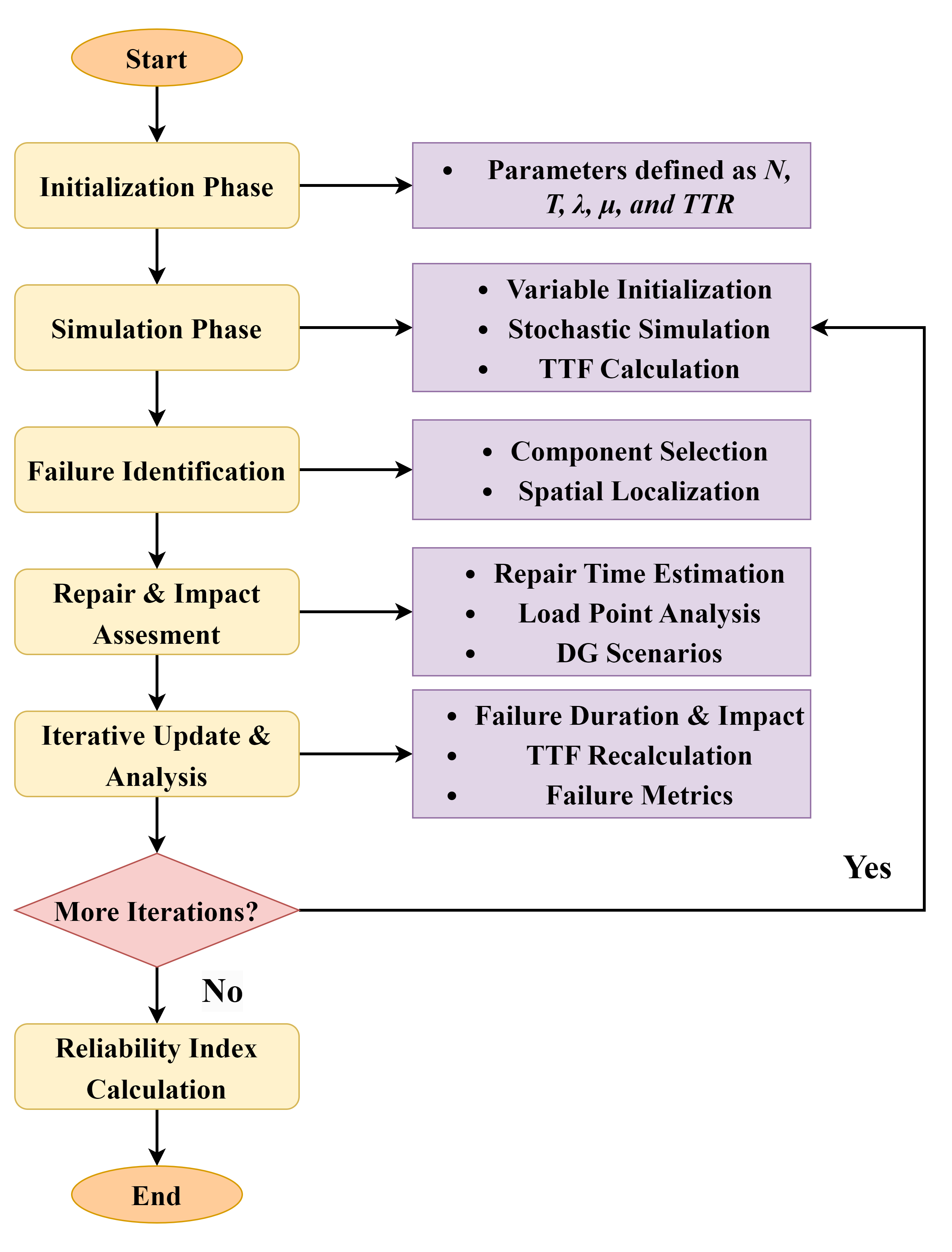} 
    \captionsetup{font=small}
    \captionsetup{justification=centering}
    \caption{A modified flowchart for reliability evaluation of distributed microgrids~\cite{lopez2020reliability}.} 
    \label{Fig.6}  
\end{figure}

Figure~\ref{Fig.6} presents a modified flowchart for the reliability evaluation of microgrids. 

Despite their advantages, simulation-based methods present notable challenges. Computational intensity remains a primary concern, particularly for large-scale microgrids, multi-year planning horizons, or studies incorporating extensive renewable variability and component degradation scenarios. Additionally, the stochastic nature of MCS introduces sampling uncertainty, necessitating careful specification of iteration counts and convergence criteria to ensure statistical validity. However, recognizing the complementary strengths of analytical and simulation-based approaches, hybrid reliability evaluation frameworks have recently emerged as a promising methodology for microgrid assessment~\cite{meera2022reliability}. 

\subsubsection{Hybrid Methods}
\label{5.1.3}

Hybrid reliability evaluation methods integrate analytical formulations with simulation-based techniques to overcome the individual limitations of each approach. While analytical methods offer computational efficiency and structural insight, they often rely on simplifying assumptions that limit their ability to capture stochastic behavior and temporal dependencies. Conversely, simulation-based methods such as Monte Carlo simulation provide high fidelity by explicitly modeling uncertainty and system dynamics, but at the cost of significant computational burden. Hybrid methods seek to balance these trade-offs by combining analytical tractability with probabilistic realism~\cite{patowary2019reliability}, making them suitable for modern, renewable-dominated microgrids.

In the reviewed literature, hybrid approaches are commonly implemented by embedding analytical reliability models within Monte Carlo simulation frameworks~\cite{lopez2020reliability,meera2022reliability}. For example, analytical techniques such as minimal cut set identification, fault tree analysis, or load-point reliability evaluation are first used to characterize fault impacts and system structure. Sequential or non-sequential Monte Carlo simulation is then applied to estimate component failure rates, outage durations, and system-level reliability indices under stochastic operating conditions. This structure enables efficient screening of critical contingencies while preserving the ability to capture time-dependent effects such as repair processes, load variation, and renewable intermittency.

Hybrid reliability evaluation has also been applied to assess islanding and load transfer strategies in microgrids~\cite{conti2014monte}. Analytical models are used to determine feasible islanding configurations and identify load points that can be supplied during component outages, while simulation quantifies the probabilistic frequency and duration of such events. Again, from an AI perspective, Ahmadi et al.~\cite{ahmadi2026comprehensive} proposed that machine learning techniques can enhance reliability evaluation methods by providing more accurate forecasting of renewable generation and load patterns, allowing probabilistic reliability evaluation with reduced computational burden compared to traditional Monte Carlo simulation. Their unified framework also suggests that fault detection capabilities should inform reliability assessment in real-time, with AI-based diagnostic output integrated into system adequacy evaluation rather than treated as separate functions. This represents a conceptual evolution from static reliability assessment conducted during planning phase to dynamic reliability monitoring during real-time operation. While concerns still remain regarding the robustness of AI-based reliability evaluation, particularly related to model uncertainty and the representativeness of training data, this shift suggests that future microgrids will not only be planned for reliability but will actively manage their own reliability margins through continuous, AI-driven state assessment.

Taken together, the preceding discussion of analytical, simulation-based, and hybrid reliability evaluation methods reveals both meaningful convergence and persistent disagreement across the reviewed literature. There is broad convergence that analytical methods, while computationally efficient, are insufficient for high-penetration renewable microgrids due to their structural inability to capture temporal dependencies, state-of-charge dynamics, and rare high-impact scarcity events. Similarly, there is growing consensus that sequential Monte Carlo simulation represents the most appropriate reference method for reliability evaluation in these systems, precisely because it preserves chronological system behavior and can represent storage depletion dynamics across multi-day shortage periods. The literature diverges, however, on what constitutes an acceptable computational trade-off for planning-stage reliability evaluation: some studies advocate full multi-year sequential simulation to capture interannual variability in extreme events, while others argue that non-sequential approaches or representative-period methods provide sufficient accuracy at manageable cost. This disagreement reflects a deeper unresolved tension between model fidelity and practical tractability that the field has not yet resolved through systematic benchmark comparisons across standardized microgrid contexts. The most significant gap, however, is not methodological but structural: the majority of reviewed studies treat planning and reliability assessment as sequential rather than co-optimized steps, with reliability evaluated post-hoc rather than embedded as a binding design driver within the optimization itself. Bridging this structural disconnect between evaluation and design represents one of the most consequential open challenges in reliability-oriented microgrid planning.

Following the same systematic data extraction and categorization protocol described in previous section, a quantitative distribution of reliability evaluation methods was derived from the reviewed literature. Studies were categorized according to their dominant reliability assessment approach, ensuring methodological consistency and comparability. 

\begin{table*}[!t]
\centering
\caption{Quantitative distribution of reliability evaluation methods in microgrid studies.}
\label{tab:reliability_method_distribution}
\begin{tabular}{lcc}
\toprule
\textbf{Method Category} & \textbf{Percentage (\%)} & \textbf{Characteristics and Applications} \\
\midrule
Simulation-based methods & 52.9 & Stochastic modeling of system behavior under uncertainty \\
Analytical methods       & 32.2 & Mathematical modeling of system states and component interactions \\
Hybrid methods           & 14.9 & Combined analytical and simulation-based reliability evaluation \\
\midrule
Total                    & 100.0 & -- \\
\bottomrule
\end{tabular}
\end{table*}

The distribution in Table~\ref{tab:reliability_method_distribution} indicates that simulation-based reliability evaluation methods dominate the literature, accounting for 52.9\% of the reported studies. This prevalence reflects the need to explicitly capture uncertainty, time-dependent behavior, and rare high-impact events in renewable-dominated microgrids, which are difficult to represent using purely analytical formulations. In contrast, analytical methods still remain relevant at 32.2\% due to their computational efficiency and interpretability, particularly for simplified systems and preliminary assessments. Hybrid methods, while comparatively less prevalent, signal an emerging research direction that seeks to balance computational tractability with modeling fidelity by integrating analytical insights with stochastic simulation. In particular, AI-enhanced reliability assessment also represents a nascent but promising area for computational acceleration. Hence, more research is required to systematically assess the robustness, applicability, and practical integration of AI-enhanced reliability assessment within established hybrid evaluation frameworks.

\subsection{Reliability Indices}
\label{5.2}

\begin{figure*}[!t] 
    \centering 
    \includegraphics[width=0.7\linewidth]{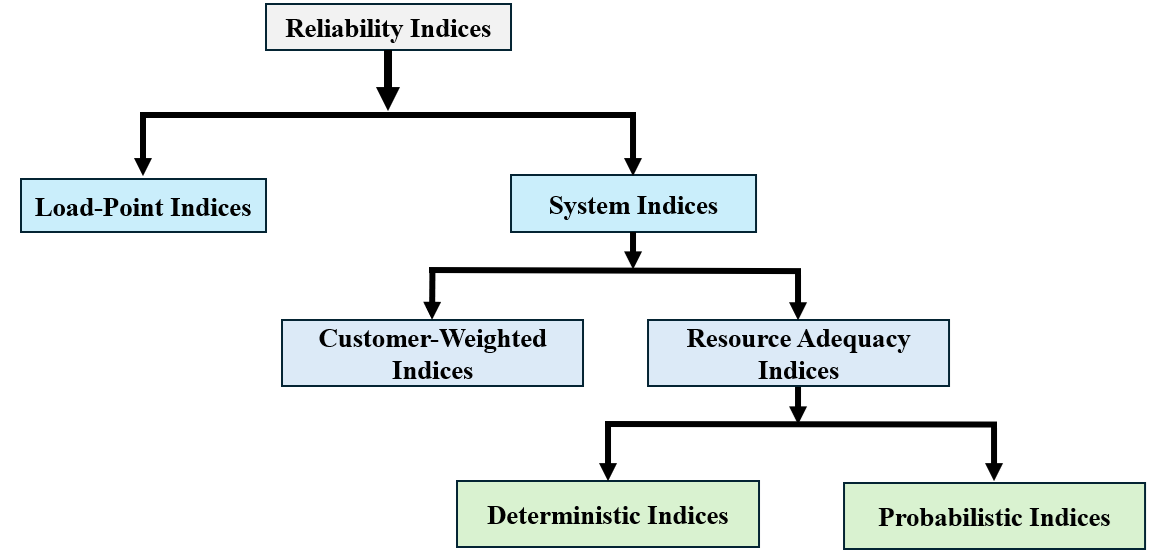} 
    \captionsetup{font=small}
    \caption{Classification of reliability indices.} 
    \label{Fig.7}  
\end{figure*}

While the previous section covered various methods used to evaluate the reliability of microgrids, such quantitative methods require indices that capture the various dimensions of system performance, ranging from individual load point behavior to aggregate system-wide impacts. The selection of reliability indices thus plays a fundamental role in planning decisions, as different metrics capture different aspects of system behavior and may lead to substantially different design outcomes. Although multiple categorization frameworks for reliability indices have been proposed in the   literature~\cite{sandelic2022reliability,meera2022reliability,lopez2020reliability}, such frameworks frequently conflate indices across distinct conceptual dimensions, leading to confusion in index selection and ambiguity in interpretation. Therefore, to improve structural clarity and analytical consistency, this review adopts a systematic hierarchical classification that distinguishes the spatial dimension (where reliability is evaluated) from the performance dimension (what aspect is quantified). As illustrated in Figure~\ref{Fig.7}, we first categorize indices spatially into Load-Point Indices and System Indices. System Indices are further sub-categorized into Customer-Weighted Indices (outcome-based) that quantify experienced service continuity from the end-user perspective, and Resource Adequacy Indices (capability-based) that assess generation sufficiency from either probabilistic or deterministic perspectives. This classification framework eliminates dimensional ambiguity and provides a clearer representation of how reliability indices are used in planning generation capacity and assessing supply continuity in microgrids.

\subsubsection{Load-point Indices}
\label{5.2.1}

Load-point indices quantify the reliability of power supply at individual nodes or customer connection points by aggregating the failure rates and repair times of the components serving the corresponding load point. These indices form the basis for system-level reliability assessment and are essential for identifying weak points in network topology~\cite{sandelic2022reliability,lopez2020reliability}. The common load-point indices are:

\begin{enumerate}
    \item \textbf{Average failure rate ($\lambda_i$)} - This index represents the expected number of interruptions occurring per year at a specific load point $i$. It is usually calculated as the sum of the failure rates of all components whose failure interrupts supply to that load point:
    \begin{equation}
    \lambda_i = \sum_{c=1}^{n} \lambda_c,\tag{4} 
    \end{equation}
    where $n$ is the total number of components whose failure results in an interruption at load point $i$, and $\lambda_c$ is the failure rate of component c (typically expressed in failures per year). It is important to note that this formulation follows the conventional assumption of a radial supply structure where each load point has a single supply path, such that any component failure causes interruption. However, microgrids frequently incorporate redundancy through multiple distributed generators, energy storage systems, or meshed network connections, enabling a load point to remain energized despite individual component failures. In such configurations, interruption occurs only when all alternative supply paths are simultaneously unavailable, requiring modification of the above formulation to account for parallel redundancy and coincident failure probabilities.

    \item \textbf{Average Annual Outage Time (Unavailability, $U_i$)} - This index represents the total expected duration of interruptions at a specific load point $i$ over a one-year period. This metric is generally expressed in hours per year and is calculated using the following formula:
    \begin{equation}
    U_i = \sum_{c=1}^{n} \lambda_c r_c, \tag{5}
    \end{equation}
    where $r_c$ denotes the average restoration time following the failure of component $c$. This index captures both the frequency of failures (through $\lambda_c$) and their duration (through $r_c$) , making it particularly relevant for assessing the cumulative impact of interruptions on customer service quality.

    \item \textbf{Average Outage Time Per Failure ($r_i$)} - This index represents the average duration of an interruption at a specific load point $i$. It is calculated as the ratio of the average annual outage time to the average failure rate, 
    \begin{equation}
    r_i= \frac{U_i}{\lambda_i}. \tag{6}
    \end{equation}
    This index provides insight into the restoration characteristics at each load point, reflecting both the nature of failures affecting that location and the effectiveness of repair and switching operations. Load points with high $r_i$ values may indicate areas where restoration is particularly challenging and may warrant targeted reliability improvement measures.
\end{enumerate}

The above-mentioned load-point indices serve as inputs to system-level calculations and provide the granularity necessary for a detailed reliability assessment. However, they do not directly capture the system-wide perspective required for aggregate planning decisions, motivating the development of system indices described in the following subsections.

\subsubsection{System Indices}
\label{5.2.2}

While load-point indices evaluate reliability at individual locations within the microgrid network, system indices aggregate load-point information to evaluate the overall reliability of the entire microgrid system~\cite{sandelic2022reliability}. These indices provide a comprehensive measure of the severity and system-wide impact of interruptions, thereby capturing the collective experience of all connected loads. In addition, system indices also play a critical role in benchmarking microgrid performance, facilitating comparisons among alternative designs, and supporting compliance with regulatory and planning requirements. As mentioned earlier, this review distinguishes system indices along two dimensions: Customer-Weighted Indices, which provide outcome-based quantification of service delivery from the customer perspective, and Resource Adequacy Indices, which offer capability-based assessment of generation-load balance. Together, these two classes of system indices provide complementary insights into microgrid reliability, addressing both the experienced quality of supply and the underlying adequacy of supplying resources.

\textbf{Customer-Weighted Indices (Performance-Based)}

Customer-weighted indices quantify reliability from the customer perspective by weighting interruption events according to the number of affected customers. These indices are typically used for historical performance assessment and are prominently featured in regulatory performance benchmarks worldwide~\cite{lopez2020reliability,meera2022reliability}.

\begin{enumerate}
    \item \textbf{System Average Interruption Frequency Index (SAIFI)}: This index measures the average frequency of sustained interruptions, defined as the expected number of sustained interruptions per customer per year. It is calculated as the ratio of the total number of customers interrupted to the total number of customers served~\cite{sandelic2022reliability}: 
          \begin{equation}
            \text{SAIFI} = \frac{\sum_{i} \lambda_i N_i}{N_T}, \tag{7}
          \end{equation}
    where \(N_i\) denotes the number of customers at load point $i$, and \(N_T\) is the total number of customers served by the system. 

    \item \textbf{System Average Interruption Duration Index (SAIDI)}: This index measures the total duration of sustained interruptions that the average customer experiences:
          \begin{equation} 
            \text{SAIDI} = \frac{\sum_{i} U_i N_i}{N_T}. \tag{8} 
          \end{equation}
    SAIDI is expressed in hours (or minutes) per customer per year and reflects the duration dimension of reliability. A lower SAIDI value means that, on average, customers experience shorter interruptions, which is indicative of a more reliable system.

    Together, SAIFI and SAIDI offer complementary insights into system reliability, as a network may experience rare but long duration outages, reflected by low SAIFI and high SAIDI, or frequent yet short interruptions, indicated by high SAIFI and low SAIDI. Therefore, both indices must be jointly evaluated to achieve a comprehensive assessment of reliability.

    \item \textbf{Momentary Average Interruption Frequency Index (MAIFI)}: This reliability index measures the average number of short-duration or brief interruptions experienced by each customer during a specified period. It is useful for capturing the frequency of short interruptions, which are often ignored by other indices like SAIFI and SAIDI. It is calculated as: 
          \begin{equation}
            \text{MAIFI} = \frac{\sum_{i} M_i N_i}{N_T}, \tag{9} 
          \end{equation}
    where \(M_i\) represents the number of momentary interruptions at load point $i$. The definition of momentary interruptions and the associated MAIFI can vary by system operators with thresholds ranging from less than one minute to up to five minutes. However interruption durations of less than one minute are particularly critical in renewable-dominated microgrids where low system inertia can lead to frequent sub-second frequency trips.

    \item \textbf{Customer Average Interruption Frequency Index \\(CAIFI)} - This index differs from \text{SAIFI} by considering only those customers who actually experienced an outage, excluding unaffected customers from the calculation:
    \begin{equation}
        \text{CAIFI} = \frac{\sum_{i} \lambda_i N_i}{N_{affected}},  \tag{10}
    \end{equation}
    where $N_{affected}$ is the total number of customers who experienced at least one outage during the reporting period. CAIFI provides insight into the reliability experienced by affected customers specifically, which may differ substantially from system-wide averages in networks with heterogeneous reliability performance. 

    \item \textbf{Customer Average Interruption Duration Index \\(CAIDI)} - This index represents the average duration of each outage event, calculated as the ratio of total customer interruption duration to total customer interruptions:
    \begin{equation} 
        \text{CAIDI} =  \frac{\sum_{i} U_i N_i}{\sum_{i} \lambda_i N_i} = \frac{SAIDI}{SAIFI}. \tag{11} 
    \end{equation}
    \text{CAIDI} reflects the average time required to restore service following an interruption and is expressed in hours (or minutes) per interruption.

   \item \textbf{Customers Experiencing Multiple Interruptions (CEMIn)} - This index characterizes the occurrence of multiple sustained interruptions affecting individual customers over a defined period:
    \begin{equation}
        \text{CEMIn} = \frac{N_{(k>n)}}{N_T},  \tag{12}
    \end{equation} 
    where ${N_{(k>n)}}$ is the number of customers experiencing more than n interruptions. In contrast to indices that aggregate and average interruption events across the full customer base, it explicitly identifies customers exposed to repeated service disruptions, thereby providing a more granular measure of reliability performance.
\end{enumerate}

\textbf{Resource Adequacy Indices}

While customer-weighted indices evaluate historical performance, resource adequacy indices assess the sufficiency of available generation and storage capacity to meet anticipated future demand. These indices address the fundamental  question of supply-demand balance that underpins microgrid planning and long-term reliability assessment. These indices are further distinguished by their treatment of uncertainty, with probabilistic indices explicitly capturing the stochastic behavior of generation and demand, whereas deterministic indices rely on fixed capacity margins derived from engineering judgment or historical performance detailed in the previous section.

\textbf{\itshape Deterministic Indices}

Deterministic indices evaluate the adequacy of the system using predefined capacity margins or reserve requirements to ensure supply sufficiency under assumed operating conditions. These indices do not explicitly model the uncertainty in generation availability or demand, resulting in relatively simple and transparent planning criteria~\cite{EPRI2023Metrics}.

\begin{enumerate}
    \item \textbf{Planning Reserve Margin (PRM)} - This metric is defined as the difference between the total installed generation capacity and the peak demand, divided by the peak demand, and typically expressed as a percentage:
    \begin{equation}
    \textbf{PRM} = \frac{C_{available} - P_{peak}}{P_{peak}} \times 100\%, \tag{13}
    \end{equation}
    where $C_{available}$ is the total available generation capacity and $P_{peak}$ is the forecasted peak demand. Typical PRM requirements in conventional power systems range from 6\% to 30\%~\cite{stenclik2024new}, providing a capacity buffer to accommodate forced outages, scheduled maintenance, and uncertainties in demand forecasts. However, Ghanbarzadeh et al.~\cite{ghanbarzadeh2025addressing} explicitly caution that PRM may be inadequate in systems with high penetrations of renewable energy resources, a limitation that is particularly pronounced in microgrids. In such systems, variable renewable generators cannot be assumed to contribute their full rated capacity during periods of peak demand. This is because the availability of renewable generation is inherently dependent on local weather conditions, which may not coincide with peak demand periods and can therefore lead to a mismatch between installed capacity and effective reliability contribution.

    \item \textbf{Energy Reserve Margin (ERM)} - This index extends PRM by focusing on energy sufficiency rather than just power capacity. It quantifies the surplus of available energy relative to total energy demand over a specified assessment horizon and is typically evaluated on an hourly basis, rather than being restricted to the peak load hour as in conventional PRM-based assessments. It is defined as:
    \begin{equation}
    \text{ERM} = \frac{E_{\text{avail}} - E_{\text{demand}}}{E_{\text{demand}}} \times 100, \tag{14}
    \end{equation}
    where $E_{avail}$ denotes the total available energy, and $E_{demand}$ represents the corresponding energy demand over the assessment period. By evaluating adequacy on an hourly basis rather than a single peak condition, ERM provides a more appropriate deterministic measure for microgrids with significant shares of variable and energy-limited resources. However, as a deterministic metric, ERM remains limited in its ability to capture uncertainty and temporal correlations between demand and renewable generation, which constrains its applicability in highly stochastic operating environments.

\end{enumerate}

\textbf{\itshape Probabilistic Indices}

Probabilistic indices evaluate the reliability of the system by explicitly representing uncertainty in demand and generation availability~\cite{ghanbarzadeh2025addressing}, offering a fundamentally different perspective from deterministic reserve-based metrics. Uncertainty is typically characterized by assigning probability distributions to stochastic system variables, allowing reliability to be evaluated over a wide range of possible operating conditions rather than a single deterministic scenario. Hence, these indices quantify the likelihood, frequency, and magnitude of generation shortfalls, providing a rigorous basis for microgrid planning under uncertainty.

\begin{enumerate}

    \item \textbf{Loss of Load Probability(LOLP)} - This index quantifies the likelihood that a microgrid is unable to meet demand (i.e., experiences a loss-of-load condition) at a given time due to insufficient net available supply under uncertainty in renewable generation, demand, and component availability. It is defined as the expected fraction of time steps in which net available supply is lower than the load demand, weighted by the probability of occurrence of uncertainty scenarios:
    \begin{equation}
        \textbf{LOLP} \;=\; \frac{1}{T}\sum_{t=1}^{T}\left(\sum_{s=1}^{S}\pi_s \,\mathbb{I}\!\left\{P_{\text{avail}}(t,s) < P_{\text{load}}(t,s)\right\}\right), \tag{15}
    \end{equation}
    where $t\in\{1,\dots,T\}$ indexes discrete time steps over the analysis horizon, $s\in\{1,\dots,S\}$ indexes uncertainty scenarios with probability $\pi_s$ (such that $\sum_{s=1}^{S}\pi_s=1$), and $\mathbb{I}\{\cdot\}$ is an indicator function equal to $1$ if the condition holds and $0$ otherwise. $P_{\text{avail}}(t,s)$ denotes the net available supply at time $t$ in scenario $s$ (e.g., renewable output after outages plus dispatchable resources and storage discharge, subject to power and state-of-charge constraints), and $P_{\text{load}}(t,s)$ denotes the load demand at time $t$ in scenario $s$. As a dimensionless metric, LOLP is typically expressed as a percentage, capturing the instantaneous risk of supply shortfall without regard to duration or magnitude.

    \item \textbf{Loss of Load Expectation (LOLE)} - This index quantifies the expected total duration over a specified planning horizon during which system demand exceeds net available supply. Unlike LOLP, which measures only the likelihood of a supply shortfall at an instant, LOLE explicitly incorporates the temporal persistence of inadequacy events by accumulating the time spent in loss-of-load conditions across the horizon:
    \begin{equation}
        \textbf{LOLE} \;=\; \sum_{t=1}^{T}\left(\sum_{s=1}^{S}\pi_s \,\mathbb{I}\!\left\{P_{\text{avail}}(t,s) < P_{\text{load}}(t,s)\right\}\right)\Delta t, \tag{16}
    \end{equation}
    where $\Delta t$ is the time-step duration (e.g., $\Delta t=1$~h). LOLE is typically expressed in hours per year or days per year and has been widely adopted as a benchmark reliability standard in generation planning. In many jurisdictions, this index is used to define acceptable levels of supply adequacy, with planning criteria commonly specified using threshold values such as 0.1 day/year (approximately 2.4 hours/year)~\cite{ghanbarzadeh2025addressing}.

    \item \textbf{Loss of Power Supply Probability (LPSP)} - This index reflects the extent to which available supply fails to meet demand, aggregated over a specified analysis period. It is defined as the ratio of the cumulative energy shortfall to the cumulative load demand, and is commonly used to assess supply adequacy:
    \begin{equation}
        \textbf{LPSP} \;=\; \frac{\sum_{t=1}^{T}\sum_{s=1}^{S}\pi_s \, P_{\text{shed}}(t,s)\,\Delta t}{\sum_{t=1}^{T}\sum_{s=1}^{S}\pi_s \, P_{\text{load}}(t,s)\,\Delta t}, \tag{17}
    \end{equation}
    where $P_{\text{shed}}(t,s)$ denotes the unserved load in scenario $s$ at time $t$, $P_{\text{load}}(t,s)$ denotes the load demand at time $t$ in scenario $s$, $\Delta t$ is the time step, and $T$ is the total number of time steps considered. LPSP varies between 0 and 1, with lower values indicating higher supply reliability. Owing to its explicit representation of energy shortfall, this index is particularly relevant for stand-alone microgrids where grid backup is unavailable and any generation shortfall directly impacts customers~\cite{thirunavukkarasu2023comprehensive}.

    \item \textbf{Expected Energy Not Supplied (EENS)} - This index represents the expected amount of electrical energy that cannot be delivered to consumers due to insufficient available generation capacity over a specified planning horizon. It is defined as the probability-weighted sum of curtailed (unserved) energy across scenarios and time steps:
    \begin{equation}
        \textbf{EENS} \;=\; \sum_{t=1}^{T}\sum_{s=1}^{S}\pi_s \, P_{\text{shed}}(t,s)\,\Delta t, \tag{18}
    \end{equation}
    where $P_{\text{shed}}(t,s)$ denotes the unserved load in scenario $s$ at time $t$, $\pi_s$ is the probability of scenario $s$, $\Delta t$ is the time-step duration, and $T$ is the total number of time steps considered. EENS is typically expressed in MWh/year and provides a direct link between reliability assessment and economic evaluation, as it can be readily combined with interruption cost or value-of-lost-load (VOLL) parameters. Hence, this index has received considerable attention in reliability-constrained planning and optimization models, where it enables explicit trade-offs between investment cost and expected reliability performance.

\end{enumerate}

Based on the above-mentioned discussion, it is evident that reliability assessment in microgrids relies on a structured set of indices that quantify system performance across multiple dimensions, including interruption occurrence, interruption duration, affected load or customer base, and energy not supplied. As reflected in the proposed hierarchical classification, reliability can be evaluated at individual load points or at the system level, and further interpreted from either the customer experience perspective or the standpoint of resource adequacy and supply sufficiency. Each class of indices captures a specific facet of system behavior, ranging from local vulnerability and service continuity to generation-demand balance under uncertainty, but no single index is capable of fully representing reliability in a comprehensive manner. This limitation arises from the inherently multi-dimensional nature of reliability, where frequency, duration, magnitude, and probability of supply shortfall must all be considered simultaneously. These challenges are further exacerbated in renewable-based microgrids, where the underlying system behavior differs significantly from the assumptions embedded in conventional reliability metrics. Consequently, a critical examination of the suitability of traditional reliability indices becomes essential in the context of modern renewable-rich microgrids.

Historically, reliability indices were developed for conventional power systems characterized by dispatchable generation, predictable demand patterns, and centralized network topologies. However, these underlying assumptions fundamentally break down in renewable-based microgrids, where system behavior is governed by decentralized configurations, inverter-dominated operation, and variable generation sources. These limitations become immediately apparent when examining load-point indices, which attribute interruptions primarily to component failures and associated repair processes. While such indices are effective for identifying structurally weak locations in radial distribution networks, they become less informative in renewable-based microgrids, where supply continuity is increasingly governed by energy availability rather than by component outages alone.

In renewable-based microgrids, a load point supplied by photovoltaic generation and battery storage may experience loss of service not because of component failure, but due to prolonged periods of low resource availability combined with depleted storage. Such events fall outside the scope of conventional load-point formulations, which implicitly assume interruptions to be discrete, component-triggered events independent of operational state. Consequently, load-point indices may underestimate reliability risks in microgrids where interruptions arise from correlated resource insufficiency rather than physical network failures.

Beyond adequacy-related limitations, conventional load-point indices also fail to capture the operationally dependent nature of component reliability in inverter-dominated microgrids. Failure rates of power electronic converters vary significantly with operating conditions such as thermal cycling, voltage stress, and environmental exposure \cite{peyghami2020overview,sandelic2022reliability}. For example, a photovoltaic inverter operating under highly variable irradiance conditions experiences markedly different stress profiles compared to steady-state operation, yet traditional load-point analysis assigns identical failure rates regardless of operating context. As a result, interruption frequency may be systematically underestimated at affected customer connection points. 

Similarly, restoration processes in microgrids are inherently state dependent. In systems with islanding capability, restoration following a feeder fault depends not only on repair actions but also on the availability of local generation and the state of charge of energy storage at the time of the event. A fault occurring during periods of high renewable output and sufficient storage may be mitigated through rapid islanding, whereas an identical fault under resource-scarce conditions may lead to prolonged outages. Such state-dependent restoration dynamics cannot be represented by static outage time parameters, limiting further the interpretability of load-point indices for reliability-oriented microgrid planning.

Notably, the limitations identified for load-point indices extend directly to system-wide reliability indices. This is due to the fact that system-wide indices are derived through aggregation of the same interruption frequencies and outage durations, thereby propagating the assumptions that limit the applicability of load-point indices in renewable-based microgrids. Among system-wide indices, customer-weighted metrics such as SAIFI, SAIDI, and CAIDI are the most widely adopted for utility benchmarking and regulatory compliance~\cite{kaddour2022impact}. However, their reliance on annual aggregation obscures seasonal and diurnal variations that are critical in renewable-based microgrids. A system may satisfy annual SAIDI or SAIFI targets while experiencing severe reliability degradation during specific periods when renewable generation is scarce and demand peaks.

In addition, customer-weighted indices implicitly assume homogeneous customer importance and do not distinguish between critical and non-critical loads. Renewable-based microgrids often serve heterogeneous load portfolios \cite{xu2016microgrids}, including essential services such as healthcare facilities, communication infrastructure, and emergency services alongside non-critical residential demand. In such contexts, a design that performs acceptably under customer-weighted averages may still fail to meet reliability requirements for high-priority loads if interruptions are concentrated on critical customers. Furthermore, these indices inadequately characterize microgrid-specific operational modes, particularly intentional and unintentional islanding. While transitions between grid-connected and islanded operation may introduce brief interruptions captured by indices such as MAIFI, these metrics cannot distinguish between planned islanding, which enhances resilience, and unplanned islanding, which indicates a reliability concern. Consequently, microgrids with advanced islanding capabilities may report higher momentary interruption frequencies than systems without islanding, inadvertently penalizing reliability-enhancing design features. Furthermore, these indices do not directly address the planning-oriented question of whether available generation and storage resources are sufficient to meet demand under uncertainty. This planning-oriented perspective is instead addressed by resource adequacy indices, which assess the balance between supply capability and demand.

Despite their widespread use in planning conventional power systems, resource adequacy indices face significant conceptual and practical challenges when applied to renewable-based microgrids. Deterministic resource adequacy indices such as the PRM implicitly assume that all installed capacity contributes equally to reliability, treating a megawatt of variable renewable generation as equivalent to a megawatt of dispatchable capacity. This assumption is invalid for weather-dependent resources whose availability may not coincide with demand peaks. Similarly, both PRM and ERM rely on static capacity representations that fail to capture the state-dependent availability of energy-limited resources such as battery storage. A battery may be fully available or entirely unavailable depending on prior dispatch history, a dependence that deterministic indices cannot represent.

In addition, deterministic indices (e.g., reserve-margin or threshold-based capacity criteria) do not quantify reliability risk, as they indicate only whether capacity exceeds a predefined threshold rather than the likelihood or severity of shortage events. A system marginally below the required reserve threshold may experience negligible shortage risk under low variability conditions, while a system comfortably exceeding the threshold may encounter substantial risk if generation and demand variability are high, yet deterministic indices provide no basis for distinguishing such cases. To address this, probabilistic indices explicitly model uncertainty, making them conceptually well-suited to renewable-dominated systems. However, the practical application of these indices to renewable-based microgrids also reveals significant limitations that must be acknowledged.

Probabilistic indices require detailed characterization of generation availability and demand uncertainty, which is challenging for renewable-based microgrids due to strong site-specificity and limited high-quality local datasets. Unlike conventional generators, whose forced outage can often be estimated from long operational histories~\cite{sergio2025weather}, renewable characterization depends on local meteorological records (and, for microgrids, behind-the-meter conditions) that may be unavailable, short in duration, or not representative of extreme events. Accurate assessment further requires joint modeling of demand and renewable availability, including temporal correlations. For instance, if peak demand occurs during evening hours when solar generation is unavailable, photovoltaic resources contribute little to reducing LOLP despite high average capacity factors. Capturing such effects requires chronological simulation or advanced probabilistic modeling that may exceed data availability and computational feasibility for microgrid planning.

Finally, threshold-based probabilistic indices such as LOLP and LOLE quantify the likelihood or expected duration of shortages but remain insensitive to shortage magnitude. Systems with identical LOLP or LOLE values may experience vastly different reliability consequences, as these indices treat minor deficits equivalently to severe supply interruptions. This limitation is particularly critical in renewable-based microgrids, where shortages often coincide with extended resource droughts and storage depletion. Even energy-based indices such as EENS, while useful for economic analysis, provide limited insight into the timing, spatial distribution, and criticality of shortages. Moreover, EENS estimates in microgrids are sensitive to the assumed operational policy (e.g., storage dispatch strategy, forecast quality, reserve allocation, and operational constraints), which can lead to optimistic reliability estimates, if perfect or fully optimal dispatch is implicitly assumed under real-world uncertainty.

Taken together, the preceding analysis demonstrates that traditional reliability indices, while well established in conventional power system planning, exhibit fundamental limitations when applied to renewable-based microgrids. Their underlying assumptions regarding  dispatchable generation, temporal independence, and homogeneous service requirements limit their ability to reflect the operational realities of decentralized, energy-limited, and inverter-dominated systems. Consequently, no single metric suffices to comprehensively characterize reliability performance in such systems. These limitations highlight the need to reconsider how reliability is evaluated and interpreted in microgrids, both in terms of the selection and combination of existing indices and in terms of more fundamental, multi-dimensional representations of reliability. Table~\ref{tab:recommendation matrix} consolidates the critical analysis of reliability indices presented in this section into a structured recommendation matrix, designed to guide practitioners in selecting appropriate metrics for specific planning contexts. In particular, it synthesizes the suitability of major reliability metrics against common planning questions in renewable-based microgrids. The purpose of this matrix is not to identify a universally superior index, but to clarify which metrics are most informative for a given planning objective, what dimension of reliability they capture, and what important limitations remain when they are used in isolation. Building on these insights, the following section discusses the findings of this review to articulate an integrated, reliability-oriented perspective on microgrid planning and design.

\begin{landscape}
\begin{table}[h]
\centering
\caption{Reliability index selection matrix for renewable-based microgrid planning}
\label{tab:recommendation matrix}

\renewcommand{\arraystretch}{1.6}
\setlength{\tabcolsep}{8pt}

\newcommand{\ccell}[2][\linewidth]{%
  \parbox[c]{#1}{\centering\strut #2\strut}%
}

\begin{tabular}{|
>{\centering\arraybackslash}m{3.5cm}|
>{\centering\arraybackslash}m{3.0cm}|
>{\centering\arraybackslash}m{3.7cm}|
>{\centering\arraybackslash}m{3.7cm}|
>{\centering\arraybackslash}m{4.2cm}|
>{\centering\arraybackslash}m{5.0cm}|}
\hline
 
\textbf{Decision Context}
& \textbf{Recommended Metric(s)}
& \textbf{What It Captures}
& \textbf{Key Strengths}
& \textbf{Key Limitations in Renewable MGs}
& \textbf{Practical Guidance} \\
\hline
 
\ccell{Topology \& Weak-Point Identification}
& \ccell{$\lambda_i$, $U_i$, $d_i$}
& \ccell{Identifies buses with the highest interruption frequency and outage duration}
& \ccell{Pinpoints structurally vulnerable locations; basis for all system-level index calculation}
& \ccell{Assumes component-induced faults only; misses adequacy-driven outages from renewable scarcity or storage depletion}
& \ccell{Apply to rank buses for reinforcement, redundancy, and protection placement; combine with adequacy indices to capture energy shortfalls.} \\
\hline
 
\multirow{2}{*}[-1.0cm]{\parbox{3.5cm}{\centering Generation \& Storage Adequacy}}
& \ccell{LOLP, LOLE, LPSP}
& \ccell{Captures the likelihood, duration, and fraction of time that supply falls short of demand under uncertainty}
& \ccell{Probabilistic; captures stochastic renewable and demand behavior; LOLE has established benchmarks, e.g., 0.1 day/yr}
& \ccell{Insensitive to shortage magnitude; requires detailed stochastic data; LPSP aggregation obscures when shortfalls occur}
& \ccell[4.6cm]{Combine LOLP/LOLE with EENS to capture occurrence and magnitude; apply LPSP for standalone or islanded sizing.} \\
\cline{2-6}
 
& \ccell{PRM, ERM}
& \ccell{Captures capacity (PRM) or energy (ERM) surplus above demand}
& \ccell{Simple and transparent; easy to communicate to stakeholders}
& \ccell{Deterministic; ignores variability, storage dynamics, and demand correlations}
& \ccell[4.6cm]{Use for preliminary screening only. Do not use as a standalone adequacy measure for renewable MGs; always supplement with probabilistic indices.} \\
\hline
 
\multirow{2}{*}[-1.0cm]{\parbox{3.5cm}{\centering Customer Service Quality}}
& \ccell{SAIFI + SAIDI}
& \ccell{Captures how often and how long the average customer is interrupted per year}
& \ccell{Widely adopted regulatory benchmarks; complementary pair covering interruption frequency and duration}
& \ccell{Annual averages mask seasonal and diurnal clustering; treat all customers as homogeneous; insensitive to unserved energy}
& \ccell[4.6cm]{Use for regulatory reporting and trend monitoring; supplement with CAIFI/CEMI$_n$ to reveal concentrated outage exposure and with EENS to reflect energy impact.} \\
\cline{2-6}
 
& \ccell{CAIFI, CAIDI, MAIFI, CEMI$_n$}
& \ccell{Captures outage exposure, restoration speed, momentary events, and repeat interruptions}
& \ccell{Reveals inequities and worst-served customers}
& \ccell{Thresholds may be arbitrary; does not distinguish critical vs non-critical loads}
& \ccell[4.6cm]{Use for equity-focused assessment and detailed outage characterization.} \\
\hline
 
\ccell{Cost--Reliability Trade-off \& Investment Justification}
& \ccell{EENS (with VOLL/CDF)}
& \ccell{Captures expected unserved energy (MWh/yr), which can be converted to monetary cost via VOLL or customer damage functions}
& \ccell{Captures shortage magnitude; enables direct economic valuation of reliability; well suited for constrained optimization}
& \ccell{A single scalar obscures timing and spatial distribution; sensitive to dispatch assumptions and VOLL estimation}
& \ccell[4.6cm]{Use as a central metric in reliability-constrained optimization; multiply by VOLL for interruption cost and pair with LOLE for investment appraisal.} \\
\hline

\ccell{Regulatory Compliance \& Benchmarking}
& \ccell{SAIFI + SAIDI + LOLE}
& \ccell{Captures customer-facing interruption performance (SAIFI, SAIDI) together with generation adequacy compliance (LOLE)}
& \ccell{Aligns with established utility reporting standards and regulatory planning benchmarks}
& \ccell{Designed for conventional centralized systems; annual aggregation masks seasonal vulnerability}
& \ccell[4.6cm]{Use as a basic set of metrics for regulatory reporting;; add EENS and load-specific analysis for renewable-based microgrids.} \\
\hline

\end{tabular}
\end{table}
\end{landscape}

\nomenclature[A]{PRM}{Planning Reserve Margin}
\nomenclature[A]{ERM}{Energy Reserve Margin}
\nomenclature[A]{SAIFI}{System Average Interruption Frequency Index}
\nomenclature[A]{SAIDI}{System Average Interruption Duration Index}
\nomenclature[A]{MAIFI}{Momentary Average Interruption Frequency Index}
\nomenclature[A]{CAIFI}{Customer Average Interruption Frequency Index}
\nomenclature[A]{CAIDI}{Customer Average Interruption Duration Index}
\nomenclature[A]{CEMIn}{Customers Experiencing Multiple Interruptions}
\nomenclature[A]{EENS}{Expected Energy Not Supplied}
\nomenclature[A]{LOLE}{Loss of Load Expectation}
\nomenclature[A]{LOLP}{Loss of Load Probability}
\nomenclature[A]{LPSP}{Loss of Power Supply Probability}

\section{Further Discussion}
\label{6}

This review provides a holistic examination of reliability-oriented microgrid design, demonstrating that reliability in renewable-based microgrids is an inherently cross-cutting outcome shaped by the interaction of planning assumptions, optimization formulations, operational flexibility mechanisms, and assessment frameworks. In contrast to prevailing approaches in the literature, where reliability is frequently treated as a secondary by-product of cost-optimal design, this study establishes it as a deliberate consequence of coordinated modeling strategies, planning methodologies, and evaluation criteria. Moreover, from a methodological standpoint, this study introduces a structured rapid review protocol combining PRISMA-based screening with AMSTAR 2 quality evaluation, tailored to the fast-evolving nature of microgrid research. Although rapid reviews are well established in the health sciences, their adoption in engineering remains limited. This work demonstrates that rigorous yet timely evidence synthesis is both feasible and valuable, positioning rapid reviews as a practical tool for researchers and decision-makers in dynamic domains that are rapidly evolving.

From the perspective of key planning factors, the consolidated evidence highlights that achieving high reliability in renewable-based microgrids is not contingent upon selecting a single optimal technology, but rather on the effective integration of multiple interdependent planning decisions. Reliability outcomes depend on coordinated choices across demand forecasting, generation and storage mix, network topology, and protection strategies, and the way reliability itself is assessed. These elements tightly interconnected as deficiencies in demand characterization, network design, or operational capability, can limit the effectiveness of appropriately sized generation and storage assets. Consequently, microgrid reliability should be conceptualized as an outcome of system-level coordination rather than as a function of isolated component performance.

A consistent finding across the literature is that reliability challenges in renewable-based microgrids are fundamentally time-dependent. Because both electrical demand and renewable generation exhibit strong temporal variability, planning frameworks must explicitly address uncertainty, forecast error, and extended periods of supply-demand imbalance. In this context, reliability is governed primarily by energy adequacy over time rather than by installed capacity alone. This adequacy–capacity distinction represents perhaps the most important conceptual shift when transitioning from conventional power system planning, where reserve margins and firm capacity provide reliable indicators of adequacy, to renewable-based microgrid planning, where reliability is instead governed by energy limitations, temporal correlation, and state dependent dynamics.

Energy storage emerges as a central flexibility resource, but its contribution to reliability is contingent upon operational dispatch strategies, state-of-charge management, and the system’s ability to sustain operation during prolonged low-renewable conditions. In particular, extended periods of \textit{Dunkelflaute}-type events, emerge from this review as a dominant yet systematically under-represented reliability risk. As discussed in section~\ref{3.1.2}, these events can persist for days or even weeks. Their impact is particularly severe when they coincide with storage depletion and demand peaks. Critically, \textit{Dunkelflaute} cannot be captured by looking at average renewable generation or short-term fluctuations. Instead, they require explicit modeling approaches that account for rare but prolonged scarcity conditions. This can be achieved through methods such as regime-based scenario construction, multi-year chronological simulations, or planning frameworks that explicitly consider system survivability under extreme conditions. Storage duration requirements should therefore be derived from the tail of the renewable scarcity distribution, rather than from average shortfall characteristics. This insight has direct implications for planning practice: tools that rely on typical meteorological year data or single-year simulations are likely to underestimate the true need for storage and backup capacity. As a result, systems designed using such approaches may appear reliable under normal conditions but remain vulnerable during prolonged renewable droughts.

Network topology and operational capabilities further shape whether available resources can effectively supply critical loads during contingency events. Features such as islanding, reconfiguration, and protection coordination determine the extent to which generation and storage assets remain usable under non-nominal conditions. These characteristics cannot be inferred from nameplate capacities or steady-state performance metrics and must be evaluated under stressed operating regimes that reflect realistic fault, islanding, and recovery scenarios. Together, these findings reinforce that reliability in renewable-based microgrids is governed by dynamic, state-dependent mechanisms that require explicit chronological representation.

Beyond technical considerations, this review underscores that reliability-oriented microgrid planning is situated within a broader socio-economic and institutional context. Microgrids are deployed to serve communities rather than merely electrical loads, and determinations of ``acceptable'' reliability are shaped by how outage impacts are valued, how critical services are prioritized, and what institutional and regulatory capacities exist to support implementation and long-term operation. Economic metrics such as the value of lost load, customer damage functions, and differentiated penalties for critical loads, translate interruption impacts into quantitative terms, while social factors, including affordability, equity, governance, and community acceptance, ultimately influence which reliability targets are feasible and sustainable. Collectively, these insights suggest that reliability-oriented microgrid design is best framed as a socio-technical, multi-criteria planning problem in which technical feasibility, economic efficiency, and service value must be considered jointly rather than in isolation.

Translating these interconnected planning considerations into implementable design decisions places optimization frameworks at the center of reliability-oriented microgrid planning. The reviewed literature reflects a clear shift away from purely deterministic sizing approaches toward multi-objective and uncertainty-aware formulations, acknowledging that reliable design must accommodate stochastic variability and competing objectives rather than average conditions alone. Within this landscape, heuristic and metaheuristic techniques dominate, followed by mathematical programming approaches and, more recently, AI-assisted methods. This methodological diversity reflects the increasing nonlinearity and combinatorial complexity of reliability-oriented planning problems, where flexible search strategies can offer advantages in exploring large solution spaces.

Importantly, the review indicates that the choice of optimization algorithm is rarely the primary determinant of planning outcomes. Instead, results are largely determined by how reliability is encoded within the problem formulation and whether candidate solutions remain feasible under uncertainty, islanded operation, and extended stress conditions. Advanced algorithms cannot compensate for simplified or poorly aligned reliability representations. Conversely, comparatively simple optimization methods can yield defensible designs when reliability constraints, chronological behavior, and operational feasibility are explicitly modeled. This finding underscores that methodological sophistication should be evaluated in terms of representational fidelity rather than algorithmic novelty alone.

A persistent limitation across the literature is the reliance on simplified reliability proxies, such as aggregated penalty terms or static reserve margins, even within advanced optimization frameworks. While these abstractions improve tractability, they obscure the time-coupled mechanisms—such as storage depletion, state-dependent operation, and prolonged renewable scarcity that ultimately govern reliability outcomes. Models that adopt explicit chronological representations, including scenario-based dispatch and storage dynamics, provide more accurate reliability characterization but incur higher computational cost. This trade-off has motivated the development of representative period selection, scenario reduction, and hybrid simulation-optimization workflows, which offer partial mitigation but remain unevenly applied and insufficiently standardized.

Beyond algorithmic formulation, the review also reveals that planning tools differ substantially in how directly reliability considerations can be embedded within design optimization. While many widely used platforms support cost-driven analysis with post-processing reliability assessment, fewer enable reliability constraints, outage scenarios, and network feasibility to be treated as binding design drivers within the optimization itself. Recent planning environments have begun to address this limitation by integrating multi-node network representations, explicit islanded operation modeling, and reliability-related constraints directly into optimization workflows. From a research perspective, this comparison underscores that the value of a planning tool lies not in its commercial maturity or computational speed, but in how transparently and defensibly it represents the mechanisms governing reliability. Consequently, optimization models and software should be evaluated not only on solution optimality, but on whether their reliability formulations are operationally meaningful and suitable for informing long-term investment and regulatory decisions in renewable-based microgrids.

Again, reliability-oriented planning is only as credible as the methods used to evaluate reliability under uncertain and evolving conditions. The reviewed literature shows a clear shift from purely analytical techniques toward simulation-based and hybrid reliability assessment, reflecting the growing complexity of renewable-based microgrids. While analytical methods remain useful for transparent interpretation and early-stage screening, they struggle to capture time-dependent behavior, nonlinear interactions, and rare but high-impact events such as extended renewable scarcity. Sequential Monte Carlo simulation is therefore widely regarded as essential, as it preserves chronological system behavior and can represent state-dependent mechanisms such as storage depletion during multi-day shortage events. However, its high computational burden continues to motivate research into acceleration strategies, including hybrid pipelines and AI-assisted surrogate models, which remain insufficiently validated for routine planning use.

These methodological choices directly influence how reliability outcomes are communicated through performance indices. As highlighted by the reviewed studies, indices determine which aspects of reliability are visible to planners and which remain hidden. Consequently, the selection and structuring of reliability indices is not merely a reporting step but an integral part of the planning process itself. This motivates the hierarchical classification adopted in this review, which separates the reliability indices by spatial scope and performance focus. The structure reduces dimensional ambiguity and provides a clearer basis for linking reliability metrics to generation planning and supply continuity in microgrids.

A central finding of this review is that conventional reliability indices, while well established and widely used, are poorly aligned with the operating realities of renewable-based microgrids. Developed primarily for centrally supplied systems with stable generation and weak time coupling, these indices can mask vulnerabilities associated with energy limitation, prolonged renewable shortfalls, and islanded operation when applied without adaptation. As a result, no single reliability index is sufficient to fully characterize microgrid reliability, and reliance on a single scalar metric may lead to incomplete or potentially misleading planning conclusions.

From this synthesis, two complementary insights for microgrid planning practice can be drawn. In the near term, reliability-oriented planning can adopt a structured and coordinated use of multiple, complementary reliability indices, rather than relying on any single metric. Interruption-based indices remain valuable for representing customer-level service experience, probabilistic adequacy measures provide insight into supply shortfall risk under renewable and demand uncertainty, and energy-based metrics support the economic valuation of unserved energy. However, the informative value of these indices depends critically on their joint interpretation, with explicit consideration of temporal characteristics (e.g., duration, sequencing, and timing of events) and differentiated service priorities. This coordinated use of existing indices constitutes a practical and immediately implementable improvement over prevailing planning practice, requiring no fundamental methodological innovation while substantially enhancing the transparency, interpretability, and defensibility of reliability assessments in microgrid design.

Looking ahead, an integrated, service-oriented reliability index tailored specifically to renewable-based microgrids is needed. Such indices should be chronological, state-aware, and spatially interpretable, explicitly reflecting storage dynamics, renewable-demand correlation, operating mode transitions, and differentiated service priorities. Importantly, these developments are not intended to replace near-term improvements to existing indices, but rather to build upon them, delineating an incremental pathway from current index-based practice toward more faithful reliability representations, suitable for optimization, investment justification, and regulatory decision-making in renewable-based microgrids.

Taken together, this review shows that advancing reliability-oriented microgrid planning requires moving beyond isolated choices of metrics, models, or algorithms toward an integrated decision-making logic. Reliability in renewable-based microgrids cannot be credibly ensured through post hoc evaluation or single-number indicators; instead, it must be consistently embedded across planning assumptions, optimization formulations, and assessment methods. A central conclusion of the review is that reliability outcomes are only meaningful when the underlying representations of uncertainty, operational feasibility, and service priorities are mutually consistent and aligned with the physical behavior of inverter-dominated, storage-dependent systems.

\begin{landscape}
\begin{figure}[p]
    \centering
    \includegraphics[width=0.97\paperwidth,height=0.97\paperheight,keepaspectratio]{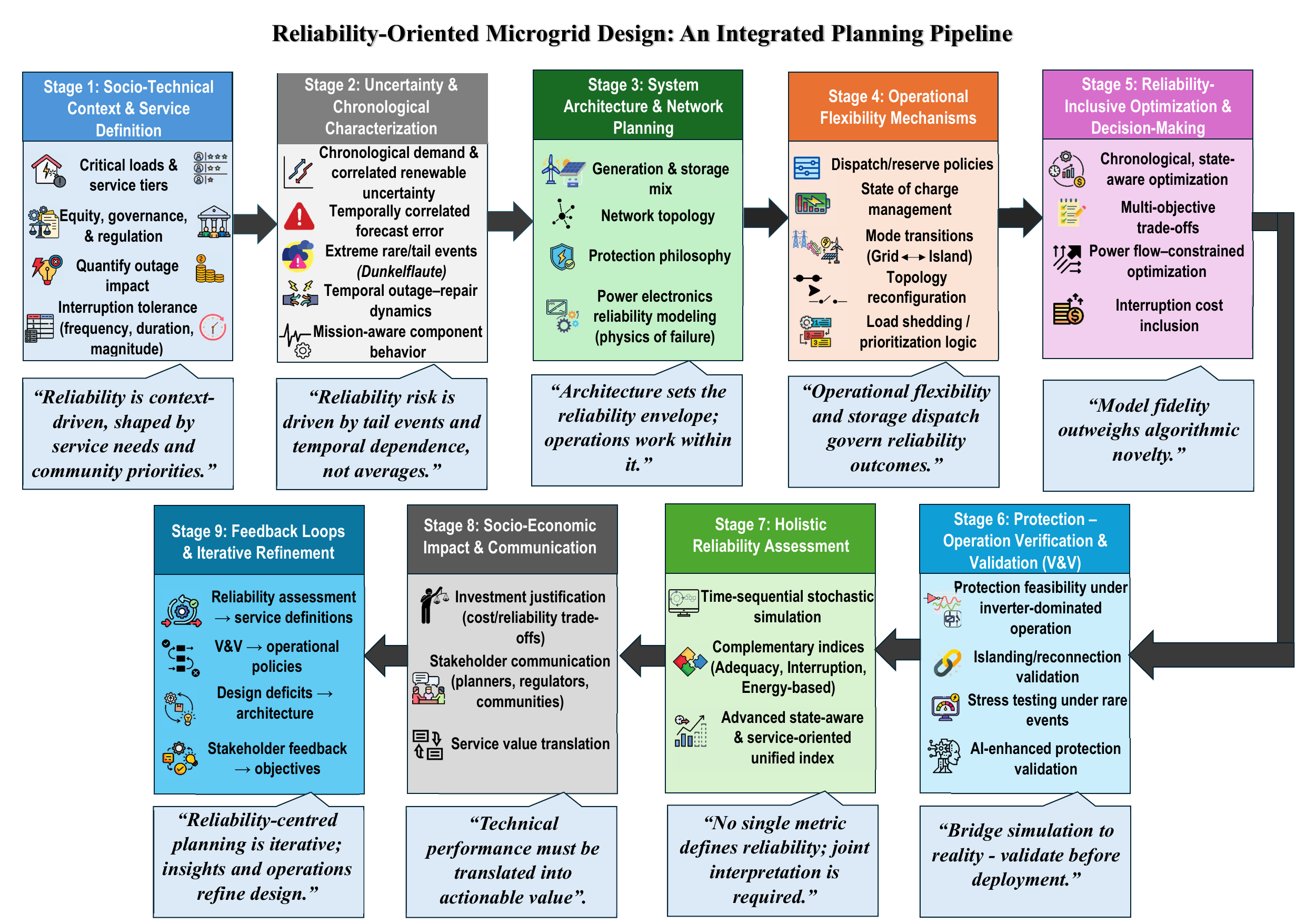}
    \caption{A unified reliability-oriented microgrid planning pipeline.}
    \label{fig:8}
\end{figure}
\end{landscape}

The reviewed evidence supports a unified planning perspective in which reliability-oriented microgrid design is best understood as a coordinated pipeline rather than a sequence of disconnected steps. This pipeline begins with defining service requirements and reliability objectives in decision-relevant terms that reflect load heterogeneity and interruption consequences, proceeds through uncertainty characterization that preserves chronology and adequacy-critical stress conditions, and integrates reliability directly into investment and operational optimization. Reliability evaluation then serves not as a detached reporting exercise, but as a mechanism for interpreting how design choices perform across frequency, duration, magnitude, and distribution of service interruptions. When implemented coherently, this integrated approach enables transparent trade-offs among reliability, cost, and operational flexibility, supporting defensible planning decisions rather than abstract performance claims. Figure~\ref{fig:8} illustrates a unified reliability-oriented microgrid planning pipeline, explicitly linking service requirements, chronological uncertainty representation, optimization formulation, and reliability evaluation into a coherent decision-making framework. In addition, to provide a formal algorithmic representation of the unified planning pipeline illustrated in Figure 8, Algorithm 1 specifies the inputs, outputs, decision logic, and feedback mechanisms at each stage. This algorithm formalizes the decision-making logic of the planning pipeline, making explicit the feedback loops between stages that are critical to ensuring convergent, reliability-compliant designs. Importantly, the algorithm highlights two key feedback mechanisms: (i) the protection verification loop (Stage 6 to Stage 3), which ensures that designs are operationally protected; and (ii) the reliability assessment loop (Stage 7 to Stage 5), which ensures that reliability targets are met before the design is finalized.

\begin{algorithm}[t]
\caption{Reliability-Oriented Microgrid Planning}
\label{alg:microgrid_pipeline}
\begin{algorithmic}[1]

\Require Load profiles $L(t)$, renewable resource data $R(t)$, candidate DER portfolio, network topology, reliability targets $(\mathrm{LOLE}_{\max}, \mathrm{EENS}_{\max})$, VOLL, budget $B_{\max}$
\Ensure Optimal DER sizing $x^*$, topology $\tau^*$, dispatch policy $\pi^*$

\Statex \textbf{Stage 1: Socio-Technical Context \& Service Definition}
\State Define critical load tiers and interruption tolerance thresholds
\State Quantify outage impact via VOLL and customer damage functions
\State Establish reliability targets from regulatory and community input

\Statex \textbf{Stage 2: Uncertainty \& Chronological Characterization}
\State Construct chronological demand--renewable scenarios $\mathcal{S} = \{s_1, \dots, s_n\}$
\State Identify tail events (e.g., \textit{Dunkelflaute} periods)
\State Assign scenario probabilities $\pi_s$, $\sum_s \pi_s = 1$

\Statex \textbf{Stage 3: System Architecture \& Network Planning}
\State Define candidate generation mix, ESS technologies, and network topologies
\State Incorporate power flow and power electronics reliability models
\State Define protection philosophy and candidate device locations

\Statex \textbf{Stage 4: Operational Flexibility}
\State Define dispatch and reserve policies; SoC management rules
\State Specify mode transitions (grid-connected/islanded)
\State Define load shedding priority and reconfiguration logic

\Statex \textbf{Stage 5: Reliability-Inclusive Optimization}
\State Formulate objective: $\min (cost)$ \& $\max (reliability)$, subject to constraints
\State Select technique:MILP for discrete decisions; stochastic/chance-constrained for uncertainty
\State Solve to obtain $(x^*, \tau^*, \pi^*)$

\Statex \textbf{Stage 6: Protection \& Operational Verification}
\State Verify protection feasibility under inverter-dominated conditions
\State Validate islanding and reconnection transitions
\State Perform stress testing under extreme scenarios
\If{verification fails}
    \State Return to Stage 3 with revised constraints
\EndIf

\Statex \textbf{Stage 7: Holistic Reliability Assessment}
\State Evaluate via time-sequential stochastic simulation
\State Compute indices: SAIFI, SAIDI, LOLE, LOLP, EENS
\If{LOLE $>$ $\mathrm{LOLE}_{\max}$ or EENS $>$ $\mathrm{EENS}_{\max}$}
    \State Return to Stage 5 with tightened constraints
\EndIf

\Statex \textbf{Stage 8: Socio-Economic Impact Assessment}
\State Compute cost--reliability trade-offs
\State Translate results for stakeholder communication

\Statex \textbf{Stage 9: Iterative Refinement}
\State Feed reliability insights into service definitions
\State Update operational policies and system architecture
\State Repeat until convergence

\State \Return $(x^*, \tau^*, \pi^*)$
\end{algorithmic}
\end{algorithm}

Finally, this review identifies several high-priority gaps that limit reliable translation from planning designs to real-world operation. These include insufficient treatment of time-coupled and correlated uncertainty, weak integration of inverter-dominated protection and mode transitions, fragmented use of reliability metrics with limited linkage to design decisions, and under-representation of service value and critical-load differentiation in planning objectives. Addressing these gaps does not require abandoning existing methods, but rather reorienting them around reliability as an organizing principle. From this perspective, future research should prioritize planning frameworks that are chronological and state-aware, employ complementary reliability measures rather than single scalar metrics, embed operational feasibility directly into optimization, and report results in ways that remain interpretable to planners, regulators, and communities. Such an approach provides a practical pathway toward microgrid designs that are not only cost-effective and sustainable, but demonstrably reliable under the full range of conditions faced by modern renewable-based systems.

\section{Findings, Recommendations, and Conclusion}
\label{7}

Based on the discussion above, the key findings, and recommendations are presented below.

Findings:
\begin{itemize}
  \item Microgrid reliability is a system-level outcome shaped by coordinated planning, not an incidental result of cost optimization.
  \item Reliability in renewable-based microgrids is governed by time-dependent energy adequacy, not installed capacity alone.
    \item Prolonged low-renewable events (\textit{Dunkelflaute}) dominate reliability risk and cannot be inferred from average resource performance.
  \item Energy storage contributes to reliability only when dispatch, state-of-charge, and duration are explicitly modeled.
  \item Power electronics reliability is mission-profile dependent; constant failure rate assumptions underestimate long-term reliability degradation.
  \item Network topology and islanding capabilities determine whether resources can supply critical loads during contingencies.
  \item Simplified reliability proxies obscure time-coupled mechanisms critical to renewable-based microgrid performance.
  \item Optimization outcomes depend more on reliability formulation than on the choice of algorithm.
  \item Chronological and uncertainty-aware models yield more credible reliability insights than static approaches.
  \item Conventional reliability indices are poorly aligned with inverter-dominated, storage-dependent renewable-based microgrids.
  \item Legacy planning tools lack outage modeling, VOLL integration, and power flow analysis critical for reliability-oriented design.
  \item No single reliability index can adequately characterize microgrid performance under uncertainty.
  \item Microgrid reliability targets are shaped by socio-economic context, service value, and critical load priorities.
\end{itemize}

Despite the breadth of the findings presented above, several limitations of this review warrant acknowledgment. First, although the adopted rapid review methodology enabled a timely synthesis of a fast-evolving field, it necessarily involved a more focused scope than a full systematic review. The streamlined screening process, while structured through PRISMA-based procedures and AMSTAR2 quality assessment, may not achieve exhaustive coverage, and some niche or interdisciplinary studies may have been excluded. Second, the synthesis corpus was intentionally centered on review-level literature, which improves interpretive breadth and enables efficient mapping of the research landscape, but may under-represent recent primary technical studies that have not yet been absorbed into published reviews. The inclusion of 59 primary studies through bidirectional snowballing mostly mitigates this limitation, but does not ensure complete coverage. Third, restricting the search to English-language publications in major databases may have excluded relevant non-English or regional studies. Fourth, substantial heterogeneity across studies in terms of modeling approaches, objectives, uncertainty treatment, and reliability metrics precluded formal meta-analysis. Therefore, the quantitative results presented reflect frequency-based methodological trends rather than statistically pooled estimates. Finally, while the adapted AMSTAR2 framework enhances transparency in quality appraisal, it does not directly measure the technical merit of engineering studies, given the differences between engineering review practices and standard systematic review protocols. These limitations do not undermine the main findings but define the boundaries within which the results should be interpreted.

  \clearpage
\clearpage

Building on the findings and evidentiary boundaries outlined above, several priority research directions emerge from the gaps identified in this review. The most fundamental of these is the need for reliability metrics suited to renewable-based microgrids. As demonstrated in section 5.2, existing indices fail to capture interruption frequency, duration, magnitude, timing, and spatial distribution under the stochastic, state-dependent conditions of inverter-dominated systems. Future work should therefore advance state-aware, service-oriented metrics that incorporate chronological system evolution, storage state dependence, and critical-load heterogeneity, moving beyond aggregated indices toward multi-dimensional representations that support investment and regulatory decision-making. Closely related is the need for integrated planning frameworks that co-represent key reliability drivers currently treated in isolation. Prolonged renewable scarcity, mission-profile-dependent component degradation, and protection constraints under inverter-dominated operation each exert first-order impacts on reliability yet are rarely captured within a unified optimization framework. Their joint representation is essential to improving the realism and credibility of planning outcomes.

The trade-off between model fidelity and computational tractability remains a central challenge. This review shows that formulation quality has a greater influence on planning credibility than solver sophistication, yet high-fidelity models become computationally prohibitive for long-term, uncertainty-rich analyses. Promising approaches, including decomposition methods, representative-period selection with fidelity guarantees, and AI-based surrogate models, offer pathways to scalability but require rigorous validation for reliability-critical applications. In parallel, the translation of technical reliability into stakeholder-relevant decision metrics remains underdeveloped. Although value of lost load and customer damage functions are recognized as essential, their application in microgrid planning is inconsistent. Standardized frameworks are needed to express reliability in terms of differentiated service levels, context-specific interruption valuation, and socio-technical criteria, enabling transparent communication of cost–reliability trade-offs.

Finally, stronger empirical grounding is required to bridge the gap between modeled and realized reliability. The literature remains dominated by simulation-based studies with limited validation against operational systems. Priority areas include long-term data from deployed microgrids, lifecycle reliability data for power electronic components, and systematic comparisons between predicted and observed performance. Without such validation, the credibility and practical applicability of planning models will remain limited. Based on these directions, a set of actionable recommendations is presented in the following subsection.

Recommendations:

\begin{itemize}
    \item Treat reliability as a governing design objective rather than a secondary outcome of cost minimization.
    \item Frame microgrid reliability as a socio-technical, multi-criteria planning problem.
    \item Explicitly model chronological demand, renewable variability, and correlated uncertainty in planning frameworks.
    \item Represent \textit{Dunkelflaute} events when determining storage duration and backup strategies.
    \item Prioritize representational fidelity over algorithmic complexity when selecting optimization methods.
    \item Represent critical-load heterogeneity using differentiated penalties, priorities, or service-level targets.
    \item Incorporate mission-profile or stress-aware converter reliability instead of fixed component failure rates.
    \item Adopt planning tools that support outage modeling, VOLL integration, and power flow analysis, rather than relying solely on techno-economic sizing tools.
    \item Use complementary reliability indices jointly rather than relying on any single aggregated metric.
    \item Develop a state-aware and service-oriented new reliability index tailored to renewable-based microgrids.

\end{itemize}

Finally, as there is a growing emphasis on incorporating fully renewable energy sources into microgrids, a critical question remains unanswered: can these fully renewable-based systems match the reliability of traditional systems? This represents a significant knowledge gap that future research must address.

\section*{Acknowledgement}
This work was supported in part by the RACE for 2030 CRC under Grant 21.PhD.N3.000006 and the Australian Research Council (ARC) Discovery Early Career Researcher Award (DECRA) under Grant DE230100046.

We dedicate this paper to the memory of our esteemed colleague, Professor Ariel Liebman, whose contributions will always be remembered.

\section*{Declaration of generative AI use}

During the preparation of this work the authors used ChatGPT to assist with writing clarity and language polishing. After using this tool, the author reviewed and edited the content as needed and takes full responsibility for the content of the published article.

\section*{Declaration of Competing Interest}
The authors declare that they have no known competing interests or personal relationships that could have appeared to influence the work reported in this paper.


\bibliographystyle{unsrt}

\bibliography{reference}
\clearpage
\onecolumn
\appendix
\begin{landscape}
\begin{center}
\textbf{Appendix A}
\end{center}
\begin{center}
\textbf{Table A1: Search String and Filters Applied to the Academic Databases.}\\
\end{center}

\begin{longtable}{|p{0.15\textwidth}|p{0.25\textwidth}|p{0.45\textwidth}|}

\hline
\textbf{Database} & \textbf{Search String} & \textbf{Filters Applied} \\
\hline
Scopus &  TITLE-ABS-KEY (``Microgrid" OR ``microgrid*" OR ``distributed*" OR ``DER" OR ``100\% Renewable") AND TITLE-ABS-KEY (reliability* OR reliable OR "Reliability" OR ``Reliability index" OR ``reliability management") AND TITLE-ABS-KEY (``Design" OR ``optimisation" OR ``optimisation*" OR ``Optimization" OR ``Planning") Edit Save Set alert & 
(LIMIT-TO (PUBYEAR, 2025) OR (LIMIT-TO (PUBYEAR, 2024) OR (LIMIT-TO (PUBYEAR, 2023) OR LIMIT-TO (PUBYEAR, 2022) OR LIMIT-TO (PUBYEAR, 2021) OR LIMIT-TO (PUBYEAR, 2020) OR LIMIT-TO (PUBYEAR, 2019) OR LIMIT-TO (PUBYEAR, 2018) OR LIMIT-TO (PUBYEAR, 2017) OR LIMIT-TO (PUBYEAR, 2016) OR LIMIT-TO (PUBYEAR, 2015) OR LIMIT-TO (PUBYEAR, 2014)) AND (LIMIT-TO (DOCTYPE, ``review")) AND (LIMIT-TO (SUBJAREA, ``ENGI") OR LIMIT-TO (SUBJAREA, ``ENER")) AND (LIMIT-TO (EXACTKEYWORD, ``Reliability") OR LIMIT-TO (EXACTKEYWORD, ``Optimization") OR LIMIT-TO (EXACTKEYWORD, ``Distributed Power Generation") OR LIMIT-TO (EXACTKEYWORD, ``Electric Power Distribution Networks") OR LIMIT-TO (EXACTKEYWORD, ``Microgrid") OR LIMIT-TO (EXACTKEYWORD, ``Reliability Analysis") OR LIMIT-TO (EXACTKEYWORD, ``Microgrids")) \\
\hline
Web of Science & {(Microgrid\* AND (Reliability OR Design OR Planning OR Optimisation)) (Topic) OR (Microgrid AND Reliability) (Title) AND (Microgrid\* AND Reliability OR Design OR Planning OR Optimisation) (Abstract)} & {Review Article (Document Types) and 2025 or 2024 or 2023 or 2022 or 2021 or 2020 or 2019 or 2018 or 2017 or 2016 or 2015 or 2014 (Publication Years) and Energy Fuels or Green Sustainable Science Technology or Engineering Electrical Electronic (Web of Science Categories) and Elsevier or Mdpi or IEEE or Wiley or Springer Nature or Taylor \& Francis (Publishers) and RENEWABLE SUSTAINABLE ENERGY REVIEWS or ENERGIES or IEEE ACCESS or APPLIED ENERGY or INTERNATIONAL JOURNAL OF ELECTRICAL POWER ENERGY SYSTEMS or ELECTRIC POWER SYSTEMS RESEARCH or INTERNATIONAL JOURNAL OF ENERGY RESEARCH or JOURNAL OF CLEANER PRODUCTION or JOURNAL OF ENERGY STORAGE or SUSTAINABILITY or JOURNAL OF MODERN POWER SYSTEMS AND CLEAN ENERGY or WILEY INTERDISCIPLINARY REVIEWS ENERGY AND ENVIRONMENT or SUSTAINABLE ENERGY TECHNOLOGIES AND ASSESSMENTS or ENERGY or RENEWABLE ENERGY FOCUS or SUSTAINABLE CITIES AND SOCIETY or ENERGY AND BUILDINGS or ENERGY FOR SUSTAINABLE DEVELOPMENT or SOLAR ENERGY or RENEWABLE ENERGY or ENERGY REPORTS or CLEAN TECHNOLOGIES or CLEAN TECHNOLOGIES AND ENVIRONMENTAL POLICY or ENERGY ENVIRONMENT AND SUSTAINABILITY or SUSTAINABLE ENERGY GRIDS NETWORKS or TECHNOLOGY AND ECONOMICS OF SMART GRIDS AND SUSTAINABLE ENERGY or IEEE TRANSACTIONS ON SMART GRID (Publication Titles)}\\
\hline
IEEE & {((``Document Title": ``Microgrid" AND (``Document Title": ``Reliability" OR ``Document Title": ``Design" OR ``Document Title": ``Planning" OR ``Document Title": ``Optimization" OR ``Document Title": ``Distributed" OR ``All Metadata": ``Network Topology" OR ``All Metadata": ``Reliability" OR ``All Metadata": ``Design" OR ``All Metadata": ``Planning") AND ``Document Title": ``Review") AND (``Abstract": ``Microgrid" AND ``Abstract": ``Reliability" OR ``Abstract": ``Design" OR ``Abstract": ``Planning"))} &  Filters Applied: Conferences, Journals, Early Access Articles, Review Article, IEEE, SGEPRI, 2014 - 2025\\
\hline
Compendex & {((((((Microgrid OR microgrid\*) AND (Reliability OR Design OR Planning OR Optimisation OR Network Topology)) WN KY) AND (((Microgrid) AND (Reliability OR Design OR Planning OR Network Topology OR Optimisation) AND (Review)))} &  Filters Applied: Review Papers, 2014 - 2025\\
\hline

\end{longtable}
\end{landscape}

\end{document}